\documentclass[12pt]{article}

\RequirePackage{amsthm,amsmath,amsfonts,amssymb}
\RequirePackage{natbib}
\RequirePackage[colorlinks,citecolor=blue,urlcolor=blue]{hyperref}
\RequirePackage{graphicx}
 \usepackage{comment}
\usepackage{dblfloatfix}
\usepackage{xcolor}
\usepackage{multirow}
\usepackage[margin=1in]{geometry}

\newcommand{\bGamma}{{\mathbf\Gamma}}
\newcommand{\bLambda}{{\mathbf\Lambda}}
\newcommand{\bSigma}{{\mathbf\Sigma}}

\def\e{{\mathbf e}}

\def\v{{\mathbf v}}

\def\x{{\mathbf x}}
\def\y{{\mathbf y}}
\def\z{{\mathbf z}}
\def\0{{\mathbf 0}}

\def\A{{\mathbf A}}
\def\B{{\mathbf B}}

\def\D{{\mathbf D}}

\def\F{{\mathbf F}}
\def\H{{\mathbf H}}
\def\L{{\mathbf L}}
\def\R{{\mathbf R}}
\def\S{{\mathbf S}}
\def\I{{\mathbf I}}
\def\M{{\mathbf M}}
\def\Q{{\mathbf Q}}
\def\T{{\mathbf T}}
\def\U{{\mathbf U}}
\def\X{{\mathbf X}}
\def\U{{\mathbf U}}
\def\V{{\mathbf V}}
\def\W{{\mathbf W}}
\def\X{{\mathbf X}}
\def\Y{{\mathbf Y}}
\def\Z {{\mathbf Z}}

\newcommand\blfootnote[1]{%
  \begingroup
  \renewcommand\thefootnote{}\footnote{#1}%
  \addtocounter{footnote}{-1}%
  \endgroup
}

\def\eop{\hfill $\Box$}
\def\tr{\mbox{tr}}
\def\Diag{\mbox{Diag}}

\newtheorem{theorem}{Theorem}
\newtheorem{lemma}{Lemma}

\newtheorem{assumption}{Assumption}

\newcommand{\change}[1]{{\leavevmode\color{red}{#1}}}

\begin{document}
\title{Limiting laws and consistent estimation criteria for fixed and diverging number of spiked eigenvalues}

\date{}
\author{}
\author{%
{Jianwei Hu$^{1^*}$, Jingfei Zhang$^{2^*}$, Jianhua Guo$^3$ and Ji Zhu$^4$}
\vspace{1.6mm}\\
\fontsize{11}{10}\selectfont\itshape $^1$\,Department of Statistics and Key Lab NAA--MOE, Central China Normal University, Wuhan, China. \\
\fontsize{11}{10}\selectfont\itshape $^2$\,Goizueta Business School, Emory University, Atlanta, GA, USA. \\
\fontsize{11}{10}\selectfont\itshape $^3$\,School of Mathematics and Statistics, Beijing Technology and Business University, Beijing, China. \\
\fontsize{11}{10}\selectfont\itshape $^4$\,Department of Statistics, University of Michigan, Ann Arbor, MI, USA. \\
\blfootnote{The first two authors contributed equally to this work.}}

\renewcommand{\baselinestretch}{1.35}
\maketitle

\vspace{-0.5in}
\begin{abstract}

In this paper, we study limiting laws and consistent estimation criteria for the extreme eigenvalues in a spiked covariance model of dimension $p$ with the number of spikes $k$.
Allowing both $p$ and $k$ to diverge, we derive limiting distributions of the spiked sample eigenvalues using random matrix theory techniques. Notably, our results are established under a general spiked covariance model, where the bulk eigenvalues are allowed to differ, and the spiked eigenvalues need not be uniformly upper bounded or tending to infinity, as have been assumed in the existing literature. Based on the above derived results, we formulate general estimation criteria that can consistently estimate $k$, while $k$ can be fixed or grow at an order of $k=o(n^{1/3})$. Our results are established under both Gaussian distributions and general distributions with finite fourth moments, with different growth rate conditions on $k$. The effectiveness of the proposed estimation criteria is illustrated through simulation studies and applications to three real-world data sets.
\end{abstract}


\newpage
\baselineskip=26.5pt

\section{Introduction}
\label{section:introduction}

Principal component analysis (PCA) has been widely used for a range of purposes including dimension reduction, data visualization, and downstream supervised or unsupervised learning \citep{abdi2010principal}.  In PCA, the leading eigenvalues and eigenvectors of the population covariance matrix are estimated using their empirical counterparts. Hence, it is important to understand properties of the sample covariance matrix, and to investigate how to determine the number of significant components using the sample covariance matrix.

Consider a sample of size $n$ denoted as $\x_1,\ldots,\x_n$ from a $p$-dimensional multivariate distribution with mean $\0$ and covariance $\bSigma_{p\times p}$. Let $\lambda_1\geq\cdots\geq\lambda_{p}$ denote the eigenvalues of $\bSigma$.
In our work, we focus on the setting with $\lambda_1\geq\cdots\geq\lambda_k>\lambda_{k+1}\geq\cdots\geq\lambda_{p}$, typically referred to as the general spiked covariance model \citep{Ke:Ma:Lin:2023}. The number, $k$, is referred to as the number of signals or the number of spikes. This model is a generalization of the standard spiked covariance model \citep{Johnstone:2001} that assumes $\lambda_1\geq\cdots\geq\lambda_k>\lambda_{k+1}=\cdots=\lambda_{p}$. Let $\S_n=\sum_{i=1}^n\x_i\x_i^\top/n$ be the sample covariance matrix and denote its eigenvalues by $d_1\geq\cdots\geq d_{p}$.
In this work, we study limiting laws of $d_1,\ldots, d_{k}$, referred to as spiked sample eigenvalues, and derive consistent criteria for estimating the number of spikes, $k$, in both fixed and diverging $p$ regimes. The three main focuses of this work are detailed as below.

As a first focus of this paper, we aim to study the limiting laws of $d_1,\ldots, d_{k}$ when $p\rightarrow\infty$. Under this regime, many efforts have been made to understand the behaviors of the eigenvalues and eigenvectors of $\S_n$; see, for example, \cite{Jonsson:1982}, \cite{Johansson:1998},  \cite{Bai:Silverstein:2004}, \cite{Paul:2007}, \cite{Nadakuditi:Edelman:2008}, \cite{Ma:2012}, \cite{Onatski:Moreira:Hallin:2013}, \cite{Wang:Yao:2013}, \cite{Wang:Silverstein:Yao:2014}, \cite{Passemier:Mckay:Chen:2015}, \cite{Wang:Fan:2017}, \cite{Bai:Choi:Fujikoshi:2018}, \cite{Cai:Han:Pan:2020}, \cite{Li:Han:Yao:2020}, \cite{Jiang:Bai:2021} and \cite{Zhang:Zheng:Pan:Zhong:2022}.
Assuming that $n,p\to\infty$, $p/n\rightarrow c>0$, we establish the limit of the spiked sample eigenvalue $d_i$ 
for $1\le i\le k$ under regularity conditions.
Notably, our result allows $k$ to diverge at the rate of $o(n^{1/3})$ for Gaussian distributions and $o(n^{1/4})$ for general distributions with finite fourth moments. This is improved over what has been established in the literature, such as $k=o(n^{1/6})$ in \cite{Cai:Han:Pan:2020}. Moreover, we show that $d_i$, $1\le i\le k$, is $\sqrt{n/k}$-consistent and demonstrate its asymptotic normality. An important improvement of our work over most known results in the literature is that our results hold even when the number of the spikes $k$ diverges. 
Furthermore, our results allow for general spiked eigenvalues $\lambda_1,\ldots,\lambda_k$, as we do not require them to be uniformly upper bounded or tending to infinity, and for general bulk eigenvalues $\lambda_{k+1},\ldots,\lambda_{p}$ that may differ from each other.

As a second focus of this paper, we aim to investigate consistent information criteria for estimating $k$ under a fixed $p$. Under this setting, commonly considered criteria for estimating the number of spikes $k$ include AIC \citep{Wax:Kailath:1985,Zhao:Krishnaiah:Bai:1986}, with a penalty term $v(k,p)$, and BIC \citep{Wax:Kailath:1985,Zhao:Krishnaiah:Bai:1986}, with a penalty term $v(k,p)\log n/2$, where $v(k,p)$ denotes the number of free parameters under $k$.
An important work by \citet{Zhao:Krishnaiah:Bai:1986} showed that any penalty term $C_nv(k,p)$ may lead to an asymptotically consistent estimator, as long as $C_n/\log\log n\to\infty$ and $C_n/n\to0$. Under general distributions with finite fourth moments and a standard spiked covariance model, we give a non-asymptotic result and show that consistency can be ensured as long as $C_n\to\infty$ and $C_n/n\to0$, which much relaxes the condition in \cite{Zhao:Krishnaiah:Bai:1986} and broadens the class of consistent estimation criteria. For example, based on our result, both $C_n=\log\log n$ and $C_n=(\log\log n)^{1/2}$ lead to consistent estimation criteria. As discussed in \citet{Kritchman:Nadler:2009} and \citet{Nadler:2010}, although BIC is asymptotically consistent, the penalty term with $C_n=\log n/2$ can be too large in finite sample and lead to underestimation, especially when the signal to noise ratio (SNR) is low. Hence, our result offers theoretical guarantee for weaker penalty terms compared to BIC. 

As a third focus of this paper, we aim to investigate consistent information criteria for $k$ when $p\rightarrow\infty$. In this regime, important progresses have been made recently on estimating the number of spikes in PCA or the number of factors in factor analysis, using methods such as information criterion \citep{Bai:Ng:2002,Bai:Choi:Fujikoshi:2018,chen2022determining}, identifying eigengap \citep{Onatski:2009,Passemier:Yao:2014,Cai:Han:Pan:2020}, hypothesis testing \citep{ke2016detecting,bao2022statistical} and eigenvalue thresholding \citep{buja1992remarks,saccenti2017considering, Dobriban:Owen:2019, Dobriban:2020, Jiang:Bai:2021,Fan:Guo:Zheng:2022,Ke:Ma:Lin:2023}.
Most existing works focused on a fixed $k$ as $n,p\rightarrow\infty$.
\citet{ke2016detecting} allowed a divergent $k$ but only for testing $H_0: \bSigma=\I$ against a spiked covariance matrix as the alternative.
Some works also required the spiked eigenvalues $\lambda_1,\ldots,\lambda_k$ to diverge \citep{Bai:Ng:2002,Cai:Han:Pan:2020} or be uniformly upper bounded \citep{Bai:Choi:Fujikoshi:2018,Jiang:Bai:2021}, while some required the population eigenvectors to satisfy certain forms of ``delocalization'' conditions \citep{Dobriban:Owen:2019, Dobriban:2020} or sparsity conditions \citep{Fan:Guo:Zheng:2022}.

In our work, we propose consistent general information criteria that allow for general spiked eigenvalues, in that they need not be divergent or uniformly upper bounded, a fixed $k$ or $k=o(n^{1/3})$, and without requiring delocalization or sparsity conditions on the eigenvectors.
Specifically, under the standard spiked covariance model with general distributions with finite fourth moments, we propose a general information criterion (GIC) that
includes AIC and BIC as special cases,
and is shown to achieve consistency under more general conditions compared to those established in \citet{Bai:Choi:Fujikoshi:2018} for AIC and BIC.
Furthermore, under the general spiked covariance model, motivated by the adjusted eigenvalue thresholding procedure in \citet{Fan:Guo:Zheng:2022}, we propose an adjusted GIC that is consistent under general distributions with finite moments of order greater than 4. In Section \ref{section:experiments}, we compare our proposed estimation criteria with many existing methods including those mentioned above, and demonstrate their efficacy.

\textbf{Main contributions.} In summary, our work makes three key contributions. First, for diverging $p$ with $p/n\rightarrow c>0$, we derive limits for $d_1,\ldots,d_k$ that are $\sqrt{n/k}$-consistent and establish their asymptotic normality. Our results accommodate an increasing $k$ up to $ o(n^{1/3})$ or $ o(n^{1/4})$, with general spiked eigenvalues and bulk eigenvalues. Second, for fixed $p$, we propose a general estimation criterion that consistently estimates $k$, achieving consistency under weaker penalty terms than the BIC. Third, under a standard spiked covariance model with diverging $p$ and $p/n\rightarrow c>0$, we introduce a consistent and unifying GIC that adapts to the signal-to-noise ratio, allowing for $k=o(n^{1/3})$. This approach is extended to general spiked covariance models with an adjusted GIC.

The rest of the paper is organized as follows.
Section \ref{section:asy} establishes limiting laws of spiked sample eigenvalues and Section \ref{sec:ic} proposes consistent estimation criteria under standard and general spiked covariance models with fixed and diverging $p,k$.
Section \ref{section:experiments} presents simulation studies that compare with several other methods and Section \ref{section:empirical} analyzes three real-world data sets. A short discussion section concludes the article. All proofs are provided in the Supplementary Materials.

\section{Asymptotics of spiked sample eigenvalues}
\label{section:asy}
We start with introducing the notation. For a matrix $\X\in\mathbb{R}^{p\times p}$, we use $\Vert\X\Vert$, $\vert\X\vert$ and $\tr(\X)$ to denote the spectral norm, determinant and trace of $\X$, respectively. We use $\Diag\{d_1,\cdots,d_p\}$ to denote a $p \times p$ diagonal matrix with diagonal elements $d_1,\cdots,d_p$ and $\Diag(\X)$ to denote a diagonal matrix with diagonal elements $X_{11},\ldots,X_{pp}$.

Let $\X=[\x_1,\ldots,\x_n]^\top$ denote an independent and identically distributed  sample of size $n$ from a $p$-dimensional multivariate distribution with mean $\0$ and covariance $\bSigma_{p\times p}$. Let $\lambda_1\geq\cdots\geq\lambda_k\geq\lambda_{k+1}\geq\cdots\geq\lambda_{p}$ denote the eigenvalues of $\bSigma$. By spectral decomposition, we may express $\bSigma$ as
\begin{equation}\label{eqn:de}
\bSigma=\bGamma\bLambda\bGamma^\top,\quad \bLambda= \begin{pmatrix}
 \bLambda_1 \,\,& \0 \\
  \0 \,\,&  \bLambda_2
\end{pmatrix},
\end{equation}
where $\bLambda_1=\Diag\{\lambda_1,\ldots,\lambda_k\}$, $\bLambda_2=\Diag\{\lambda_{k+1},\ldots,\lambda_p\}$, and $\bGamma=(\bGamma_1, \bGamma_2)$ is a $p\times p$ orthogonal matrix collecting the eigenvectors of $\bSigma$. We make the following assumptions.

\begin{assumption}\label{ass1}
Write $\x_i=\bSigma^{\frac{1}{2}}\y_i$ and $\Y=[\y_1,\ldots,\y_n]^\top\in\mathbb{R}^{n\times p}$, where $Y_{ij}$'s are independent and identically distributed random variables with $\mathbb{E}(Y_{ij})=0$ and $\mathrm{Var}(Y_{ij})=1$.
\end{assumption}

\begin{assumption}\label{ass2} $\lambda_1\geq\cdots\geq\lambda_k>\lambda_{k+1}\geq\cdots\geq\lambda_{p}$.
\end{assumption}

\begin{assumption}\label{ass3}
$n,p\to\infty$ and $c_n=(p-k)/n\rightarrow c>0$.
\end{assumption}

The spiked covariance model under Assumption \ref{ass2} is usually referred to as the general spiked covariance model \citep{Ke:Ma:Lin:2023}. Its special case when $\lambda_1\geq\cdots\geq\lambda_k>\lambda_{k+1}=\cdots=\lambda_{p}$ is often called the standard spiked covariance model \citep{Johnstone:2001}.

Define $\Z=[\z_1,\ldots,\z_n]^\top\in\mathbb{R}^{n\times p}$, where $\z_i=\bGamma^\top\x_i$, and we have $\frac{1}{n}\Z^\top\Z$ and $\frac{1}{n}\X^\top\X$ have the same eigenvalues. In the subsequent development, we focus on the analysis of $\Z$, as what has been done in the literature \citep{Paul:2007,Bai:Yao:2008,Wang:Fan:2017,fan2018large,Cai:Han:Pan:2020,Jiang:Bai:2021,Fan:Guo:Zheng:2022,Zhang:Zheng:Pan:Zhong:2022}.
We denote the sample covariance matrix as $\S_n=\frac{1}{n}\Z^\top\Z$ and the eigenvalues of $\S_n$ as $d_1\geq\cdots\geq d_p$. Write $\Z=[\Z_1,\Z_2]$ with $\Z_1\in\mathbb{R}^{n\times k}$ and $\Z_2\in\mathbb{R}^{n\times (p-k)}$. It is seen that $\Z=\Y\bGamma\bLambda^{\frac{1}{2}}$, $\Z_1=\Y\bGamma_1\bLambda_1^{\frac{1}{2}}$ and $\Z_2=\Y\bGamma_2\bLambda_2^{\frac{1}{2}}$.

Next, we state some known results on the empirical distribution of sample eigenvalues.
Let $\{\beta_j\}_{k+1\le j\le p}$ be the eigenvalues of $\S_{22}=\Z_2^\top\Z_2/n$. Then the empirical spectral distribution of $\S_{22}$ is defined as $F_n(x)=\frac{1}{p-k}\sum_{j=k+1}^{p}I_{(-\infty,x]}(\beta_j)$, where $I_\mathcal{A}(\cdot)$ is the indicator function of $\mathcal{A}$.  Let $H_n(x)=\frac{1}{p-k}\sum_{j=k+1}^{p}I_{(-\infty,x]}(\lambda_j)$ be the empirical spectral distribution of $\bLambda_{2}$, and assume that $H_n(x)$ converges weakly to a limiting spectral distribution $H(x)$ as $p\rightarrow\infty$. It is known that the empirical spectral distribution $F_n(x)$ converges weakly to $F^{c,H}$, a nonrandom c.d.f. that depends on $c$ and $H(x)$, as $p\rightarrow\infty$ \citep{Silverstein:1995,Bai:Silverstein:2004}.
Denote by $\mathcal{F}_{H}$ the support for $H(x)$ on $\mathbb{R}$.  For $x\notin\mathcal{F}_{H}$ and $x\neq 0$, we define
$$
\psi(x)=x+c\int\frac{xt}{x-t}dH(t).
$$
Correspondingly, the derivative of $\psi(x)$ is
$\psi'(x)=1-c\int\frac{t^2}{(x-t)^2}dH(t)$.
Now we additionally make the following assumption.
\begin{assumption}\label{ass4}
$\psi'(\lambda_k)>0$ and $\lambda_{k+1}$ is bounded.
\end{assumption}
Note that a sufficient condition for $\psi'(\lambda_k)>0$ is $\lambda_k/\lambda_{k+1}>1+\sqrt{c}$. The $i$-th largest eigenvalue $\lambda_i$ is said to be a distant spiked eigenvalue if $\psi'(\lambda_i)>0$ \citep{Bai:Yao:2012}.

Define $F^{c_n,H_n}(x)$ from $F^{c,H}(x)$ with $c$ and $H$ replaced by $c_{n}$ and $H_n$, respectively. Similar to \citet{Bai:Silverstein:2004,Jiang:Bai:2021}, we use $c_n$ and $H_n(x)$ rather than $c$ and $H(x)$ as the convergence of $c_n \rightarrow c$ and that of $H_n(x)\rightarrow H(x)$ may be arbitrarily slow. Let
$$
\psi_n(x)=x+c_n\int\frac{xt}{x-t}dH_n(t).
$$
Define the following Stieltjes transformations,
\[
s(d)=\int\frac{1}{x-d}dF^{c_n,H_n}(x), \,\, \underline{s}(d)=-\frac{1-c_n}{d}+c_ns(d),
\]
\[
s_n(d)=\frac{1}{p-k}\sum_{j=k+1}^{p}\frac{1}{\beta_j-d}=\int\frac{1}{x-d}dF_{n}(x), \,\, \underline{s}_n(d)=-\frac{1-c_n}{d}+c_ns_n(d),
\]
and
\[
m_1(d)=\int\frac{x}{d-x}dF^{c_n,H_n}(x),\,\,m_2(d)=\int\frac{x^2}{(d-x)^2}dF^{c_n,H_n}(x), \,\,m_3(d)=\int\frac{x}{(d-x)^2}dF^{c_n,H_n}(x).
\]
It is seen that $m_3(d)=-m_1'(d)$. Define \[
\hat{s}_n(d)=\frac{1}{p-k}\sum_{j=k+1}^{p}\frac{1}{d_j-d}, \,\, \hat{\underline{s}}_n(d)=-\frac{1-c_n}{d}+c_n\hat{s}_n(d).
\]


Our first result gives the limits for the spiked sample eigenvalues $d_1\geq\cdots\geq d_k$ and the consistent estimators for the spiked eigenvalues $\lambda_1\geq\cdots\geq \lambda_k$, which will be used to establish the consistency of GIC and AGIC, respectively.
\begin{theorem}\label{thm1}
Suppose Assumptions \ref{ass1}-\ref{ass4} and one of the following two conditions hold:
\begin{enumerate}
\item[(i)] $Y_{ij}$ is Gaussian and $k=o(n^{1/3})$,
\item[(ii)] $\mathbb{E}Y_{ij}^4<\infty$ and $k=o(n^{1/4})$.
\end{enumerate}
Then, for all $1\leq i\leq k$, we have
\begin{equation}\label{eq:thm1}
\frac{d_i}{\psi_n(\lambda_i)}\stackrel{p}{\longrightarrow}1,
\end{equation}
\begin{equation}\label{eq:thm01}
\frac{-\hat{\underline{s}}_n^{-1}(d_i)}{\lambda_i}\stackrel{p}{\longrightarrow}1,
\end{equation}
where $\hat{\underline{s}}_n^{-1}(d_i)=1/\hat{\underline{s}}_n(d_i)$.
\end{theorem}

When $\lambda_i\rightarrow\infty$, \eqref{eq:thm1} can be alternatively written as $d_i/\lambda_i\stackrel{p}{\longrightarrow}1$, as $\lambda_i/\psi_n(\lambda_i)\rightarrow1$.
Several existing results on high-dimensional PCA are relevant to Theorem \ref{thm1}.
Assuming that $n,p\to\infty$, $p/n\rightarrow c>0$ and $k$ is fixed, \cite{Baik:Silverstein:2006} established \eqref{eq:thm1} under the condition that $\lambda_1$ is bounded, \citet{Bai:Choi:Fujikoshi:2018} established \eqref{eq:thm1} under the condition that $\lambda_1$ is bounded or $\lambda_k\rightarrow\infty$, and \citet{Bai:Ding:2012} established \eqref{eq:thm01} under the condition that $\lambda_1$ is bounded.
Assuming $\lambda_{k+1}\ge\cdots\ge\lambda_{p}$ are uniformly bounded and $\lambda_k\rightarrow\infty$, \citet{Cai:Han:Pan:2020} established the following result
$$
\frac{d_i}{\lambda_i}-1=\mathcal{O}_p\left(\frac{p}{n\lambda_i}+\frac{1}{\lambda_i}+\frac{k^4}{n}\right),
$$
under the condition that $k=\min\{o(n^{1/6}), o(\lambda_k^{1/2})\}$.
In comparison, our result in Theorem \ref{thm1} allows $k$ to diverge, at a faster rate of $o(n^{1/3})$ and also relaxes the conditions on the spiked eigenvalues $\lambda_1,\ldots,\lambda_k$, in that we only require $\psi'(\lambda_k)>0$ and $\lambda_{k+1}$ is bounded. In Section \ref{sec:ic}, \eqref{eq:thm1} will be used to establish the consistency of GIC for the standard spiked covariance model and \eqref{eq:thm01} will be used to establish the consistency of AGIC for the general spiked covariance model.

Our next result shows that the spiked sample eigenvalues $d_1,\ldots,d_k$ are $\sqrt{n/k}$-consistent.
Recall $\S_{22}=\frac{1}{n}\Z_2^\top\Z_2$.
By the definition of eigenvalues, each $d_i$ solves the equation
\[
0=|d \I-\S_n|=|d \I-\S_{22}| \cdot |d \I-\mathbf K_n(d)|,
\]
where the $k\times k$ random matrix $\mathbf K_n(d)$ is defined as
\[
\mathbf K_n(d)=\frac{1}{n}\Z_1^\top(\I+\H_n(d))\Z_1,
\]
\[
\H_n(d)=\frac{1}{n}\Z_2(d \I-\S_{22})^{-1}\Z_2^\top.
\]
As the spectrum of $\S_{22}$ goes inside the support of $F^{c,H}$ in probability \citep{Bai:Yao:2012} and Theorem \ref{thm1} shows that $d_1,\ldots,d_k$ go outside the support of $F^{c,H}$ in probability, we can deduce 
that $d_1,\ldots,d_k$ solve the determinant equation
\[
\left| d\I-\mathbf  K_n(d)\right|=0.
\]
Define the random matrix $\A(d)$ as
\begin{eqnarray*}
\A_n(d)&=&\bLambda_{1}^{-\frac{1}{2}}(d\I-\mathbf K_n(d))\bLambda_{1}^{-\frac{1}{2}}\\
&=&\Diag\left\{\frac{d-\frac{1}{n}\tr(\I+\H_n(d))\lambda_1}{\lambda_1},\ldots,\frac{d-\frac{1}{n}\tr(\I+\H_n(d))\lambda_k}{\lambda_k}\right\}+\M_n(d),
\end{eqnarray*}
where $\M_n(d)=\frac{1}{n}\tr(\I+\H_n(d))\I-\bLambda_{1}^{-\frac{1}{2}}\mathbf K_n(d)\bLambda_{1}^{-\frac{1}{2}}$.
It is seen from the definition of $\A_n(d)$ that $\left| d\I-\mathbf  K_n(d)\right|=0$ if and only if $\vert \A_n(d)\vert=0$.
In our proof, we show that the solutions to $\vert \A_n(d)\vert=0$ tend to the solutions to $
d/\lambda_i-\frac{1}{n}\tr(\I+\H_n(d))=0$, $1\le i\le k$.

The following lemma gives the order at which the $i$-th diagonal element of the random matrix $\A_n(d_i)$, $1\leq i\leq k$, converges to 0.
\begin{lemma}
\label{lemma:diagonal}
Suppose Assumptions \ref{ass1}-\ref{ass4} hold, and there exists a positive constant $c_0$ not depending on $n$ such that $\lambda_{i}/\lambda_{i+1}\geq c_0>1$ for $1\le i\le k-1$. Then, for $1\le i\le k$,
\begin{enumerate}
\item[(i)] when $Y_{ij}$ is Gaussian and $k=o(n^{1/3})$, it holds that $\A_n(d_i)_{ii}=\mathcal{O}_p(k^{3/2}/n)$,
\item[(ii)] when $\mathbb{E}Y_{ij}^4<\infty$ and $k=o(n^{1/4})$, it holds that $\A_n(d_i)_{ii}=\mathcal{O}_p(k^{5/2}/n)$.
\end{enumerate}
Here $\A_n(d_i)_{ii}$ denotes the $i$-th diagonal element of $\A_n(d_i)$.
\end{lemma}

In Lemma \ref{lemma:diagonal}, we establish the order at which the {$i$-th} diagonal element of the random matrix $\A_n(d_i)$, $1\leq i\leq k$, converges to 0, which depends on both $k$ and $n$. This result requires a nontrivial proof and is crucial in establishing the convergence rate of $d_i$ in Theorem \ref{thm2} and the asymptotic normality of $d_i$ in Theorem \ref{thm3}.
The assumption on $\lambda_{i}/\lambda_{i+1}$'s requires the spiked eigenvalues to be well-separated. Similar assumptions have been made in \citet{Wang:Fan:2017,Cai:Han:Pan:2020,Zhang:Zheng:Pan:Zhong:2022}.

\begin{theorem}\label{thm2}
Suppose Assumptions \ref{ass1}-\ref{ass4} hold, and there exists a positive constant $c_0$ not depending on $n$ such that $\lambda_{i}/\lambda_{i+1}\geq c_0>1$ for $1\le i\le k-1$. Then, for $1\le i\le k$,
\begin{enumerate}
\item[(i)] when $Y_{ij}$ is Gaussian and $k=o(n^{1/3})$, it holds that
\begin{equation*}
\frac{d_i-\psi_n(\lambda_i)}{\lambda_i}=\mathcal{O}_p(\sqrt{k/n}),
\end{equation*}
\[
\frac{-\hat{\underline{s}}_n^{-1}(d_i)-\lambda_i}{\lambda_i}=\mathcal{O}_p(\sqrt{k/n}),
\]
\item[(ii)] when $\mathbb{E}Y_{ij}^4<\infty$ and $k=o(n^{1/4})$, it holds that
\begin{equation}\label{eq:thm3}
\frac{d_i-\psi_n(\lambda_i)}{\lambda_i}=\mathcal{O}_p(k/\sqrt{n}),
\end{equation}
\begin{equation}\label{eq:thm03}
\frac{-\hat{\underline{s}}_n^{-1}(d_i)-\lambda_i}{\lambda_i}=\mathcal{O}_p(k/\sqrt{n}).
\end{equation}
\end{enumerate}
\end{theorem}

It was shown in \citet{Cai:Han:Pan:2020} that
\begin{equation}\label{eq:cai}
\frac{d_i}{\lambda_i}-1=\mathcal{O}_p\left(\frac{p}{n\lambda_i}+\frac{1}{\lambda_i}+\frac{k^4}{n}\right),
\end{equation}
assuming $\lambda_{k+1}\ge\cdots\ge\lambda_{p}$ are uniformly bounded, $\lambda_k\rightarrow\infty$ and $k=\min\{o(n^{1/6}), o(\lambda_k^{1/2})\}$. The result in \eqref{eq:thm3} can be written as
$$
\frac{d_i-\psi_n(\lambda_i)}{\lambda_i}=\frac{d_i}{\lambda_i}-1-c_n\int\frac{t}{\lambda_i-t}dH_n(t)=\mathcal{O}_p(k/\sqrt{n}).
$$
In comparison, our result allows $k$ to diverge, at a faster rate of $o(n^{1/4})$ and also relaxes the conditions on the spiked eigenvalues $\lambda_1,\ldots,\lambda_k$, in that we only require $\psi'(\lambda_k)>0$ and $\lambda_{k+1}$ is bounded.
The two convergence rates in \eqref{eq:thm3} and \eqref{eq:cai} are not directly comparable due to the involvement of $\lambda_i\rightarrow\infty$ in \eqref{eq:cai}. However, if we assume, for example, $\lambda_k=\mathcal{O}(k^\alpha)$, $\alpha>2$, and $k=o(n^{1/\beta})$, $\beta=2+2\alpha$, then the rate in \eqref{eq:thm3} improves over that in \eqref{eq:cai}.

\begin{theorem}\label{thm3}
Suppose Assumptions \ref{ass1}-\ref{ass4} hold, and there exists a positive constant $c_0$ not depending on $n$ such that $\lambda_{i}/\lambda_{i+1}\geq c_0>1$ for $1\le i\le k-1$.
Assume one of the following two conditions holds:
\begin{enumerate}
\item[(i)] $Y_{ij}$ is Gaussian and $k=o(n^{1/3})$,
\item[(ii)] $\mathbb{E}Y_{ij}^4=\nu_4<\infty$ and $k=o(n^{1/5})$.
\end{enumerate}
Then, for $1\le i\le k$,
\begin{equation}\label{eq:thm4}
\frac{\sqrt{n}(d_i-\psi_n(\lambda_i))}{\lambda_i\sigma_{\lambda_i}}\stackrel{d}{\longrightarrow}N(0,1),
\end{equation}
where under (i) $\sigma_{\lambda_i}^2=2\psi_n'(\lambda_i)$ and under (ii)
$\sigma_{\lambda_i}^2=(\nu_4-3)[\psi_n'(\lambda_i)]^2\sum_{j=1}^p\Gamma_{ji}^4+2\psi_n'(\lambda_i)$ with $\bGamma$ defined as in \eqref{eqn:de}, and
\begin{equation}\label{eq:thm04}
\frac{\sqrt{n}(-\hat{\underline{s}}_n^{-1}(d_i)-\lambda_i)}{\lambda_i\tau_{\lambda_i}}\stackrel{d}{\longrightarrow}N(0,1),
\end{equation}
where under (i) $\tau_{\lambda_i}^2=2/\psi_n'(\lambda_i)$ and under (ii) $\tau_{\lambda_i}^2=(\nu_4-3)\sum_{j=1}^p\Gamma_{ji}^4+2/\psi_n'(\lambda_i)$.
\end{theorem}

Assuming that $k$ is fixed and $\lambda_1$ is bounded, \cite{Paul:2007} established \eqref{eq:thm4} under the assumption that $Y_{ij}$ is Gaussian, \citet{Bai:Yao:2008} established \eqref{eq:thm4} under the assumption that $\mathbb{E}Y_{ij}^4=\nu_4<\infty$ and a block diagonal $\bSigma$, and \cite{Zhang:Zheng:Pan:Zhong:2022} established \eqref{eq:thm4} under the assumption that $\mathbb{E}Y_{ij}^4=\nu_4<\infty$ and a general $\bSigma$. Assuming that $\lambda_k\rightarrow\infty$ and $k=\min\{o(n^{1/6}), o(\lambda_k^{1/2})\}$, \citet{Cai:Han:Pan:2020} established asymptotic normality of spiked eigenvalues under the assumption that $\mathbb{E}Y_{ij}^4=\nu_4<\infty$ and a general $\bSigma$. Compared with most existing results, one important contribution of Theorem \ref{thm3} is that it allows $k$ to diverge, and it does not require the spiked eigenvalues to be bounded.

In our work, Theorem \ref{thm2} is an important intermediate result in establishing Theorem \ref{thm3}. While Theorem \ref{thm3} provides a stronger result when compared to Theorem \ref{thm2}, it also imposes stronger assumptions on $k$. Specifically, the conditions $k = o(n^{1/3})$ and $k = o(n^{1/4})$ in Theorem \ref{thm2} are less restrictive than the requirements $k = o(n^{1/3})$ and $k = o(n^{1/5})$ in Theorem \ref{thm3}. We also note that Theorems \ref{thm2} and \ref{thm3} are of independent interest and are not applied in Section 3.

In what follows, we generalize Lemma \ref{lemma:diagonal}, Theorem \ref{thm2} and Theorem \ref{thm3} to the case of multiple spikes. Let $S(i)=\{1\leq j\leq k: \lambda_j=\lambda_i\}$ and let $|S(i)|$ be the cardinality of $S(i)$.

\renewcommand{\thelemma}{\arabic{lemma}$^*$}
\setcounter{lemma}{0}
\begin{lemma}
\label{lemma:diagonal11}
Suppose Assumptions \ref{ass1}-\ref{ass4} hold, and $|S(i)|<\infty$ for $1\le i\le k$. Then, for $1\le i\le k$,
\begin{enumerate}
\item[(i)] when $Y_{ij}$ is Gaussian and $k=o(n^{1/3})$, it holds that $$|\A_n(d_i)_{{S(i)\times S(i)}}|=\mathcal{O}_p((\frac{k^{2}}{n})^{|S(i)|}),$$
\item[(ii)] when $\mathbb{E}Y_{ij}^4<\infty$ and $k=o(n^{1/4})$, it holds that $$|\A_n(d_i)_{{S(i)\times S(i)}}|=\mathcal{O}_p((\frac{k^{3}}{n})^{|S(i)|}).$$
\end{enumerate}
Here $|\A_n(d_i)_{{S(i)\times S(i)}}|$ denotes the determinant of the diagonal block of $\A_n(d_i)_{S(i)\times S(i)}$.
\end{lemma}

\renewcommand{\thetheorem}{\arabic{theorem}$^*$}
\setcounter{theorem}{1}
\begin{theorem}\label{thm222}
Suppose Assumptions \ref{ass1}-\ref{ass4} hold, and $|S(i)|<\infty$ for $1\le i\le k$. Then, for $1\le i\le k$,
\begin{enumerate}
\item[(i)] when $Y_{ij}$ is Gaussian and $k=o(n^{1/3})$, it holds that
\begin{equation*}
\frac{d_i-\psi_n(\lambda_i)}{\lambda_i}=\mathcal{O}_p(\sqrt{k/n}),
\end{equation*}
\[
\frac{-\hat{\underline{s}}_n^{-1}(d_i)-\lambda_i}{\lambda_i}=\mathcal{O}_p(\sqrt{k/n}),
\]
\item[(ii)] when $\mathbb{E}Y_{ij}^4<\infty$ and $k=o(n^{1/4})$, it holds that
\begin{equation*}
\frac{d_i-\psi_n(\lambda_i)}{\lambda_i}=\mathcal{O}_p(k/\sqrt{n}),
\end{equation*}
\[
\frac{-\hat{\underline{s}}_n^{-1}(d_i)-\lambda_i}{\lambda_i}=\mathcal{O}_p(k/\sqrt{n}).
\]
\end{enumerate}
\end{theorem}

\renewcommand{\thetheorem}{\arabic{theorem}$^*$}
\setcounter{theorem}{2}
\begin{theorem}\label{eq:thm333}
Suppose Assumptions \ref{ass1}-\ref{ass4} hold, and $|S(i)|<\infty$ for $1\le i\le k$.
Assume one of the following two conditions holds:
\begin{enumerate}
\item[(i)] $Y_{ij}$ is Gaussian and $k=o(n^{1/4})$,
\item[(ii)] $\mathbb{E}Y_{ij}^4=\nu_4<\infty$ and $k=o(n^{1/6})$.
\end{enumerate}
Then, for $1\le i\le k$, the joint distribution of
\begin{equation*}
\frac{\sqrt{n}(d_i-\psi_n(\lambda_i))}{\lambda_i},\,\,\,\, i\in S(i),
\end{equation*}
converges to the joint distribution of the $|S(i)|$ eigenvalues of the Gaussian random matrix
\[
\mathrm{T}(\lambda_i)=\frac{\lambda_i\psi'(\lambda_i)}{\psi(\lambda_i)}\left[\bLambda_{1}^{-\frac{1}{2}}\mathbf R(\psi(\lambda_i))\bLambda_{1}^{-\frac{1}{2}}\right]_{S(i)\times S(i)},
\]
and the joint distribution of
\begin{equation*}
\frac{\sqrt{n}(-\hat{\underline{s}}_n^{-1}(d_i)-\lambda_i)}{\lambda_i},\,\,\,\, i\in S(i),
\end{equation*}
converges to the joint distribution of the $|S(i)|$ eigenvalues of the Gaussian random matrix
\[
\mathrm{\tilde{T}}(\lambda_i)=\frac{\lambda_i}{\psi(\lambda_i)}\left[\bLambda_{1}^{-\frac{1}{2}}\mathbf R(\psi(\lambda_i))\bLambda_{1}^{-\frac{1}{2}}\right]_{S(i)\times S(i)},
\]
where $\mathbf R(\psi(\lambda_i))$ is the Gaussian matrix limit of the random matrix $\mathbf R_n(\psi_n(\lambda_i))$, and $\mathbf R_n(\psi_n(\lambda_i))=\{\Z_1^\top(\I+\H_n(\psi_n(\lambda_i)))\Z_1-\tr(\I+\H_n(\psi_n(\lambda_i)))\bLambda_{1}\}/\sqrt{n}$.
\end{theorem}
Note that $\mathrm{T}(\lambda_i)$ is a  Gaussian random variable with mean 0, $1\leq i\leq k$, and the explicit form of its variance is given in \citet{Bai:Yao:2008}. The above result matches with the result in \citet{Bai:Yao:2008}, which was shown under the fixed $k$ case.

\section{Consistent information criteria for estimating $k$}
\label{sec:ic}
In Section \ref{section:ILP}, we first discuss consistent information criterion under the standard spiked covariance model with a fixed $p$, then generalize to the standard spiked covariance model with a diverging $p$ in Section \ref{sec:GIC} and to the general spiked covariance model with a diverging $p$ in Section \ref{sec:AGIC}.

\subsection{Standard spiked covariance model with a fixed $p$}\label{section:ILP}

Under the standard spiked covariance model, the information criteria, such as AIC and BIC, are scale invariant \citep{Bai:Choi:Fujikoshi:2018}. Hence, without loss of generality, we make the following assumption.
\renewcommand{\theassumption}{\arabic{assumption}$^*$}
\setcounter{assumption}{1}
\begin{assumption}\label{ass22} $\lambda_1\geq\cdots\geq\lambda_k>\lambda_{k+1}=\cdots=\lambda_{p}=\lambda=1$.
\end{assumption}

Under Assumption \ref{ass22}, the signal-to-noise ratio (SNR) can be calculated as
\[
\text{SNR}=\frac{\lambda_k-\lambda}{\lambda}=\lambda_k-1.
\]

Correspondingly, under Assumption \ref{ass2}, the signal-to-noise ratio (SNR) is
\[
\text{SNR}=\frac{\lambda_k-\lambda_{k+1}}{\lambda_{k+1}}=\lambda_k/\lambda_{k+1}-1.
\]

It is known that under the Gaussian assumption, the log-likelihood of the spiked covariance model, defined through Assumptions \ref{ass1} and \ref{ass22}, can be written as a function of sample eigenvalues $d_1,\ldots,d_k$ and $\bar{d}_{k+1}$ \citep{Anderson:2003,Wax:Kailath:1985,Zhao:Krishnaiah:Bai:1986}, given as
\begin{equation}\label{eq:L}
\log L_{k}=-\frac{n}{2}\left\{\sum_{i=1}^{k}\log d_i+(p-k)\log\bar{d}_{k+1}\right\},
\end{equation}
where $\bar{d}_{k+1}$ is defined as
\begin{equation}\label{eq:dbar}
\bar{d}_{k+1}=\frac{1}{p-k}\sum_{i=k+1}^pd_i.
\end{equation}
The AIC criterion, based on \citet{Akaike:1974}, estimates $k$ by minimizing
\[
\text{AIC}(k')=-\log L_{k'}+k'(p-k'/2+1/2),
\]
and BIC (or MDL), based on \citet{Schwarz:1978}, estimates $k$ by minimizing
\[
\text{BIC}(k')=-\log L_{k'}+k'(p-k'/2+1/2)\log n/2.
\]
Under the fixed $p$ regime, several works have pointed out that AIC is not consistent, i.e.,
$$
\lim_{n\rightarrow\infty}\mathbb{P}\left[\arg\min_{k'} \text{AIC}(k')=k\right]<1,
$$
while BIC is consistent, i.e.,
$$
\lim_{n\rightarrow\infty}\mathbb{P}\left[\arg\min_{k'} \text{BIC}(k')=k\right]=1.
$$
See \citet{Wax:Kailath:1985} and \citet{Zhao:Krishnaiah:Bai:1986} for more details. Note that BIC is a special case of a more general estimation criterion proposed in \citet{Zhao:Krishnaiah:Bai:1986}, defined as
\begin{equation}
\label{eq:gic}
\mathrm{GIC}(k')=-\log\tilde L_{k'}+k'(p-k'/2+1/2)C_n,
\end{equation}
where $\log\tilde L_{k'}$ approximates $\log L_{k'}$ via a quadratic expansion and is defined as
\begin{equation}\label{eq:tildeL}
\log\tilde L_{k'}=-\frac{n}{2}\left\{\sum_{i=1}^{k'}\log d_i+\sum_{i=k'+1}^p(d_i-1)\right\}.
\end{equation}
The two terms $\log\tilde L_k$ and $\log L_{k}$ have the same distribution asymptotically (see Lemma \ref{lemma:4}).
Thus, it is seen that when $C_n=\log n/2$, the estimation criterion $\mathrm{GIC}(k')$ is asymptotically equivalent to $\text{BIC}(k')$.
\citet{Zhao:Krishnaiah:Bai:1986} showed that GIC in \eqref{eq:gic} is consistent as long as $C_n/\log\log n\to\infty$ and $C_n/n\to0$.
In what follows, we show that these conditions on $C_n$ can be further improved and  GIC is consistent under general distributions with finite fourth moments.

\renewcommand{\thetheorem}{\arabic{theorem}}
\setcounter{theorem}{3}
\begin{theorem}
\label{thm4}
Suppose that Assumptions \ref{ass1}, \ref{ass22} hold, $\mathbb{E}Y_{ij}^4<\infty$ and $p$ is fixed while $n\rightarrow\infty$. Let $\emph{GIC}(k)$ be defined as in \eqref{eq:gic}. There exists a constant $\epsilon_k>0$ that depends on $\lambda_k$ such that, for $k'<k$, it holds with probability at least $1-c_1/n$ that
$$
\mathrm{GIC}(k')-\mathrm{GIC}(k)\ge\frac{n}{2}(k-k')\left(\epsilon_k/2-2pC_n/n\right),
$$
and for $k'>k$, it holds with probability at least $1-c_2/C_n$ that
$$
\mathrm{GIC}(k')-\mathrm{GIC}(k)\ge\frac{(k'-k)C_n}{4},
$$
where $c_1$ and $c_2$ are some positive constants.
\end{theorem}

The above result shows that GIC in \eqref{eq:gic} is consistent if $C_n\to\infty$ and $C_n/n\to 0$, which further relaxes the conditions $C_n/\log\log n\to\infty$ and $C_n/n\to0$ in \cite{Zhao:Krishnaiah:Bai:1986}, broadening the class of consistent estimation criteria.
For example, $C_n=\log\log n$ and $C_n=(\log\log n)^{1/2}$ both lead to consistent estimation criteria under Theorem \ref{thm4}.

As pointed out in \citet{Kritchman:Nadler:2009} and \citet{Nadler:2010}, although BIC (i.e., GIC with $C_n=\log n/2$) is asymptotically consistent, the penalty term $C_n=\log n/2$ can be too large in finite sample and may, especially when the SNR is low, lead to underestimation. Our result thus offers theoretical guarantee for GIC with a smaller penalty term compared to that in BIC and under general distributions beyond Gaussianity. In Section \ref{section:experiments}, we demonstrate that GIC with $C_n=\log\log n$ and $C_n=(\log\log n)^{1/2}$ have good finite sample performances, and both outperform BIC when the SNR is low.


\subsection{Standard spiked covariance model with a diverging $p$}\label{sec:GIC}
In this section, we propose a general information criterion (GIC) that consistently estimates the number of spikes $k$ as $n,p\to\infty$.

Assuming that $n,p\to\infty$, $p/n\rightarrow c>0$ and $k$ is fixed, \citet*{Bai:Choi:Fujikoshi:2018} (referred to as ``BCF" herein after), proposed the following method for estimating the number of principal components $k$.
When $0<c<1$, BCF estimates $k$ by minimizing
\begin{equation}
\label{eq:aic1}
\ell_1(k')=-\sum_{i=k'+1}^p\log d_i+(p-k')\log \bar d_{k'+1}-(p-k'-1)(p-k'+2)/n,
\end{equation}
where $\bar{d}_{k'+1}=\sum_{i=k'+1}^pd_i/(p-k')$, and when $c>1$, BCF estimates $k$ by minimizing
\begin{equation}
\label{eq:aic2}
\ell_2(k')=-\sum_{i=k'+1}^{n-1}\log d_i+(n-1-k')\log \check{d}_{k'+1}-(n-k'-2)(n-k'+1)/p,
\end{equation}
where $ \check{d}_{k'+1}=\sum_{i=k'+1}^{n-1}d_i/(n-1-k')$.
The criterion in \eqref{eq:aic1} is equivalent to AIC and the criterion in \eqref{eq:aic2} is
referred to as quasi-AIC in \citet{Bai:Choi:Fujikoshi:2018}.
Assuming $k$ is fixed, the consistency of BCF is established by assuming the following gap condition when $0<c<1$,
\begin{equation}
\label{eq:consistency3}
\psi(\lambda_k)-1-\log \psi(\lambda_k)>2c,
\end{equation}
and the following gap condition when $c>1$,
\begin{equation}
\label{eq:consistency4}
\psi(\lambda_k)/c-1-\log(\psi(\lambda_k)/c)>2/c.
\end{equation}
Note that BCF essentially defines two different criteria for cases $0<c<1$ and $c>1$, respectively.

Now with a possibly divergent $k$, we consider GIC that estimates $k$ by minimizing
\begin{equation}
\label{eq:aic}
\textrm{GIC}(k')=\frac{n}{2}\left\{\sum_{i=1}^{k'}\log d_i+(p-k')\log\bar{d}_{k'+1}\right\}+\gamma k'(p-k'/2+1/2),
\end{equation}
where $\bar{d}_{k'+1}=\sum_{i=k'+1}^pd_i/(p-k')$ and $\gamma>0$. It is seen that when $\gamma=1$, (\ref{eq:aic}) is equivalent to AIC; when $\gamma=\log n/2$, (\ref{eq:aic}) is equivalent to BIC.
Unlike the finite $p$ case, GIC does not rely on the quadratic expansion of $\log L_{k'}$, as the residuals from the expansion no longer tend to zero when $p,n\rightarrow\infty$. Alternatively, we exploit the fact that $\log L_{k'}$ can be written as a function of eigenvalues of the sample covariance matrix so that we may use the results derived in Section \ref{section:asy}.

Let
\[
\varphi(x)=\frac{1}{2c}[\psi(x)-1-\log \psi(x)],
\]
where $\psi(x)=x+cx/(x-1)$. The definition of $\varphi(x)$ is related to the formula in \eqref{eq:aic}. Specifically, the definition of $\varphi(x)$ is related to the expansion of $\log L_k-\log L_{k'}$ in GIC$(k')-$GIC$(k)$ for all $k'<k$.
Noting that
\[
\varphi((1+\sqrt{c})\lambda_{k+1})=\varphi(1+\sqrt{c})=1/2+\sqrt{1/c}-\log(1+\sqrt{c})/c,
\]
we have $\varphi(1+\sqrt{c})\in(1/2,1)$ for $c>0$.
We assume the following gap condition.
\renewcommand{\theassumption}{\arabic{assumption}}
\setcounter{assumption}{4}
\begin{assumption}\label{ass6}
$\varphi(1+\sqrt{c})<\gamma<\varphi(\lambda_{k})$.
\end{assumption}

The parameter $\gamma$ in Assumption \ref{ass6} is bounded from below and above. To establish the consistency of GIC, we mainly use the results of Theorem \ref{thm1} and the fact that for all $k<i\leq k'=o(p)$, it holds that $d_{i}\stackrel{a.s.}{\longrightarrow} (1+\sqrt{c})^2$. As a result, we need
$\gamma<\varphi(\lambda_{k})$ and  $\gamma>\varphi(1+\sqrt{c})$ for all $k'<k$ and $k<k'=o(p)$, respectively.
If $\lambda_k\rightarrow\infty$, Assumption \ref{ass6} holds trivially for any constant $\gamma>\varphi(1+\sqrt{c})$. Note that this condition is more relaxed than \eqref{eq:consistency3} in \citet{Bai:Choi:Fujikoshi:2018} as $\gamma$ is allowed to be less than 1.

\begin{theorem}
\label{thm5}
Suppose Assumptions \ref{ass1}, \ref{ass22}, \ref{ass3}-\ref{ass6} and one of the following two conditions holds:
\begin{enumerate}
\item[(i)] $Y_{ij}$ is Gaussian and $k=o(n^{1/3})$,
\item[(ii)] $\mathbb{E}Y_{ij}^4<\infty$ and $k=o(n^{1/4})$.
\end{enumerate}
Let $\emph{GIC}(k')$ be as defined in (\ref{eq:aic}).
Then for all $k'<k$, we have
\[
\mathbb{P}\left[\emph{GIC}(k)<\emph{GIC}(k')\right]\rightarrow1,
\]
and for all $k<k'=o(p)$, we have
\[
\mathbb{P}\left[\emph{GIC}(k)<\emph{GIC}(k')\right]\rightarrow1.
\]
\end{theorem}


As both AIC and BIC are special cases of GIC in \eqref{eq:aic}, the result in Theorem \ref{thm5} can be used to establish the consistency of AIC and BIC. Specifically, Theorem \ref{thm5} implies that AIC (i.e., $\gamma=1$) is consistent under condition \eqref{eq:consistency3}.
Since $\gamma=1$ for AIC, Assumption \ref{ass6} holds if \eqref{eq:consistency3} is satisfied, as $\varphi(1+\sqrt{c})<1$ and $\gamma=1$.
Thus, our result shows that AIC is consistent for any $c>0$ as long as \eqref{eq:consistency3} is satisfied. This generalizes the result in \citet{Bai:Choi:Fujikoshi:2018}, which only establishes the consistency of AIC when $0<c<1$.
Moreover, our result also allows $k$ to be fixed or diverge at $k=o(n^{1/3})$.

GIC adapts to the SNR via the parameter $\gamma$ in the penalty term. \citet{Bai:Choi:Fujikoshi:2018} showed that AIC has a positive probability for underestimating when \eqref{eq:consistency3} is not satisfied (i.e., low SNR). In this case, GIC can still achieve consistency as long as Assumption \ref{ass6} holds. When $\lambda_k\rightarrow\infty$, it is seen that $\varphi'(\lambda_k)>0$ and thus $\varphi(\lambda_k)$ increases with $\lambda_k$. When $\lambda_k/\log n>c$, we have $\varphi(\lambda_k)>\gamma=\log n/2$. Thus, Theorem \ref{thm5} also implies that BIC, which is GIC by setting $\gamma=\log n/2$, is consistent if $\lambda_k/\log n>c$.

In practice, we need to specify the value of $\gamma$. The upper bound for $\gamma$, i.e., $\varphi(\lambda_{k})$, is difficult to calculate since $\lambda_{k}$ is unknown. However, the lower bound $\varphi(1+\sqrt{c})$ can be calculated once $p$ and $n$ are given. As a practical choice, we may set $\gamma$ to be slightly larger than $\varphi(1+\sqrt{c})$, such as $\gamma=1.1\varphi(1+\sqrt{c})$, which performs well in our simulations. We now give a further rationale. By Taylor's expansion,
 \[
\lim_{c\rightarrow0+}\varphi(1+\sqrt{c})=1/2+\sqrt{1/c}-(\sqrt{c}-c/2)/c=1.
\]
This implies that, when $\gamma=\varphi(1+\sqrt{c})$, (\ref{eq:aic}) can be seen as a natural extension of  AIC, which was first introduced in the literature for fixed $p$ (i.e., $c\rightarrow0+$).
We also note that $\varphi'(1+\sqrt{c})<0$ and  $\varphi(1+\sqrt{c})\in(1/2,1)$ for $c>0$. Compared with AIC, $\gamma=\varphi(1+\sqrt{c})$ decreases  with the increase of $c$. Thus, when $\gamma=\varphi(1+\sqrt{c})$, (\ref{eq:aic}) may efficiently avoid the problem that AIC suffers from, that is, the fixed penalty $\gamma=1$  may be too large with the increase of $c$. However, as a natural extension of AIC for $c\rightarrow0+$, when $\gamma=\varphi(1+\sqrt{c})$, (\ref{eq:aic}) may also tend to slightly overestimate $k$. Thus, a bit larger $\gamma$ is a good choice.

\subsection{General spiked covariance model with a diverging $p$}\label{sec:AGIC}

We now consider the general spiked covariance model, where the bulk eigenvalues $\lambda_{k+1},\ldots,\lambda_{p}$ need not be the same. We start out investigation of the information criterion from the view of the penalized likelihood function, although our final result does not rely on the Gaussian assumption. For the general spiked covariance model , 
it holds that \citep{Anderson:2003}
\begin{equation}\label{eq:L2}
\log L_{k}=-\frac{n}{2}\left\{\sum_{i=1}^{k}\log d_i+\sum_{i=k+1}^{p}\log d_i\right\}=-\frac{n}{2}\left\{\sum_{i=1}^{p}\log d_i\right\}.
\end{equation}
It is seen that the $\log L_{k}$ in \eqref{eq:L2} actually does not vary with $k$ and this is different from the standard spiked covariance model in \eqref{eq:L}. As a result, the $\log L_{k}$ in \eqref{eq:L2} is not a suitable choice in an information criterion that selects $k$ as all $k$ gives the same $\log L_{k}$ value.

Under the general spiked covariance model, \citet{Fan:Guo:Zheng:2022} proposed an adjusted correlation thresholding (ACT) method that estimates $k$ by $\hat{k}=\max\{i: \hat{\lambda}^C_i>1+\sqrt{c}\}$,
where $\hat{\lambda}^C_i=-\underline{m}_{n,i}^{-1}(d_i)$ and
\[
m_{n,i}(z)=\frac{1}{p-i}\left[\sum_{j=i+1}^p\frac{1}{d_j-z}+\frac{1}{\frac{3}{4}d_i+\frac{1}{4}d_{i+1}-z}\right],
\]
\[
\underline{m}_{n,i}(z)=-\left(1-\frac{p-i}{n-1}\right)z^{-1}+\frac{p-i}{n-1}m_{n,i}(z),
\]
with $\{d_i\}_{1\le i\le p}$ being the eigenvalues of the sample correlation matrix of $\x_1,\ldots,\x_p$. Assuming that $n,p\to\infty$, $p/n\rightarrow c>0$ and $k$ is fixed, \citet{Fan:Guo:Zheng:2022} proved that ACT can estimate $k$ consistently under the condition that $\lambda_1$ is bounded.

Inspired by \citet{Fan:Guo:Zheng:2022}, we consider instead the correlation matrix, in which case the sample eigenvalues sum to $p$.
Specifically, define $\D=\Diag(\bSigma)$ and $\D_n=\Diag(\S_n)$ . 
Let $\R=\D^{-1/2}\bSigma\D^{-1/2}$ and $\R_n=\D_n^{-1/2}\S_n\D_n^{-1/2}$ be the population and sample correlation matrices, respectively.
To simplify notation, we still use $\{\lambda_i\}_{1\leq i\leq p}$ and $\{d_i\}_{1\leq i\leq p}$ to denote the eigenvalues of $\R$ and $\R_n$, respectively.
Since $\sum_{i=1}^pd_i=\tr(\R_n)=\tr(\S_n\D_n^{-1})=p$, we have
\[
\sum_{i=1}^{k}\log d_i+\sum_{i=k+1}^{p}\log d_i=\log(\prod_{i=1}^pd_i) < p\log (\sum_{i=1}^pd_i/p)=0.
\]
That is, $\sum_{i=k+1}^{p}\log d_i<-\sum_{i=1}^{k}\log d_i$, which implies that $\sum_{i=k+1}^{p}\log d_i=-\gamma_1\sum_{i=1}^{k}\log d_i$ for some $\gamma_1>1$.
Thus, we have
\begin{equation*}
\log L_{k}=-\frac{n}{2}(1-\gamma_1)\sum_{i=1}^{k}\log d_i.
\end{equation*}
Compared with \eqref{eq:L}, there are additional $p-k-1$ parameters that need to be estimated, that is, $\lambda_{k+2},\dots, \lambda_{p}$ and $k(p-k/2+1/2)+p-k-1=k(p-k/2-1/2)+p-1$.
Hence, one natural option is to estimate $k$ by minimizing
\begin{equation}
\label{eq:aic4}
-\frac{n}{2}\sum_{i=1}^{k'}\log d_i+\gamma k'(p-k'/2-1/2),
\end{equation}
for some $\gamma>0$. However, our analysis shows that establishing the consistency of \eqref{eq:aic4} requires $\frac{1}{2c}\log\psi((1+\sqrt{c})\lambda_{k+1})<\gamma<\frac{1}{2c}\log\psi(\lambda_{k})$, which makes selecting an appropriate $\gamma$ infeasible as $\lambda_{k}, \lambda_{k+1}, \dots, \lambda_{p}$ are unknown in practice.

To overcome this disadvantage, we consider the adjusted GIC that estimates $k$ by minimizing
\begin{equation}
\label{eq:aic3}
\mathrm{AGIC}(k')=-\frac{n}{2}\sum_{i=1}^{k'}\log \frac{d_i}{1+\delta_i}+\gamma k'(p-k'/2-1/2),
\end{equation}
where $\delta_i=\frac{1}{n}\sum_{j=i+1}^p\frac{d_j}{d_i-d_j}$ and $\gamma>0$.
We assume the following assumptions and show that AGIC in (\ref{eq:aic3}) is consistent under general distributions with finite moments of order greater than 4.
\renewcommand{\theassumption}{\arabic{assumption}$^*$}
\setcounter{assumption}{3}
\begin{assumption}\label{ass44}
$\psi'(\lambda_k)>0$ and $\lambda_{k+1}\leq 1$.
\end{assumption}

\renewcommand{\theassumption}{\arabic{assumption}$^*$}
\setcounter{assumption}{4}
\begin{assumption}\label{ass666}
$\frac{1}{2c}\log (1+\sqrt{c})<\gamma<\frac{1}{2c}\log \lambda_k$.
\end{assumption}

\begin{theorem}
\label{thm6}
Suppose Assumptions \ref{ass1}-\ref{ass3}, \ref{ass44}, \ref{ass666} hold, and there exists a positive constant $c_0$ not depending on $n$ such that $\lambda_{i}/\lambda_{i+1}\geq c_0>1$ for all $1\le i\le k-1$.
Assume one of the following two conditions holds:
\begin{enumerate}
\item[(i)] $Y_{ij}$ is Gaussian and $k=o(n^{1/3})$,
\item[(ii)] $\mathbb{E}|Y_{ij}|^{4+\zeta}<\infty$ for some constant $\zeta>0$ and $k=o(n^{1/4})$.
\end{enumerate}
Let $\mathrm{AGIC}(k')$ be defined as in (\ref{eq:aic3}).
Then for all $k'<k$, we have
\[
\mathbb{P}\left[\mathrm{AGIC}(k)<\mathrm{AGIC}(k')\right]\rightarrow1,
\]
and for all $k<k'=o(p)$, we have
\[
\mathbb{P}\left[\mathrm{AGIC}(k)<\mathrm{AGIC}(k')\right]\rightarrow1.
\]
\end{theorem}

We also note that AGIC in (\ref{eq:aic3}) can be further simplified as
\begin{equation}
\label{eq:aic5}
\hat{k}=\max\{i: \log \frac{d_i}{1+\delta_i}>\gamma\},
\end{equation}
where $\delta_i=\frac{1}{n}\sum_{j=i+1}^p\frac{d_j}{d_i-d_j}$ and $\gamma>0$.

To show the consistency of AGIC in (\ref{eq:aic5}), we make the following assumption.
\renewcommand{\theassumption}{\arabic{assumption}$^{**}$}
\setcounter{assumption}{4}
\begin{assumption}\label{ass777}
$\log (1+\sqrt{c})<\gamma<\log \lambda_k$.
\end{assumption}

\begin{theorem}
\label{thm7}
Suppose Assumptions \ref{ass1}-\ref{ass3}, \ref{ass44}, \ref{ass777} hold and one of the following two conditions holds:
\begin{enumerate}
\item[(i)] $Y_{ij}$ is Gaussian and $k=o(n^{1/3})$,
\item[(ii)] $\mathbb{E}|Y_{ij}|^{4+\zeta}<\infty$ for some constant $\zeta>0$ and $k=o(n^{1/4})$.
\end{enumerate}
Let $\mathrm{AGIC}$ be defined as in (\ref{eq:aic5}), where $i=o(p)$.
Then, we have
\[
\mathbb{P}(\hat{k}=k)\rightarrow1.
\]
\end{theorem}

As a practical choice, we set $\gamma=1.1\times\frac{1}{2c}\log (1+\sqrt{c})$ and $\gamma=1.1 \log (1+\sqrt{c})$ for  AGIC in (\ref{eq:aic3}) and (\ref{eq:aic5}), respectively. Both versions yield similar empirical results. For simplicity, we recommend using $\mathrm{AGIC}$ in (\ref{eq:aic5}) with $\gamma=1.1\log (1+\sqrt{c})$.

\section{Simulation studies}
\label{section:experiments}
In this section, we evaluate the finite sample performance of the proposed estimation criteria, and compare that with existing methods, including BIC, AIC, the modified AIC \citep[mAIC;][]{Nadler:2010}, the PC$_3$ estimator proposed in \cite{Bai:Ng:2002}, the ON$_2$ estimator proposed in \cite{Onatski:2009}, the KN estimator proposed in \cite{Kritchman:Nadler:2008}, the BCF estimator proposed in \cite{Bai:Choi:Fujikoshi:2018}, the DDPA estimator proposed in \cite{Dobriban:Owen:2019} and \cite{Dobriban:2020}, the ACT estimator proposed in \cite{Fan:Guo:Zheng:2022} and the BEMA estimator proposed in \cite{Ke:Ma:Lin:2023}. We set $\gamma=1.1\varphi(1+\sqrt{c})$ and $\gamma=1.1\log (1+\sqrt{c})$ for GIC in \eqref{eq:aic} and AGIC in \eqref{eq:aic5}, respectively.

We consider two different settings including the case of small $p$, corresponding to the method proposed in Section \ref{section:ILP}, and the case of large $p$, corresponding to the methods proposed in Sections \ref{sec:GIC} and \ref{sec:AGIC}. To evaluate the estimation accuracy, we report the fraction of times of successful recovery, i.e., $\widehat{\mathbb{P}}(\hat{k}=k)$, and the average selected number of spikes, i.e., $\widehat{\mathbb{E}}(\hat{k})$. All results are based on 200 replications.

\medskip
\textbf{Simulation 1 (the standard spiked covariance model, small $p$).}  
In this simulation, we evaluate the performance of our proposal in \eqref{eq:gic} with iterated logarithm penalties $C_n=\log\log n$ (referred to as ILP) and $C_n=(\log\log n)^{1/2}$ (referred to as ILP$_{1/2}$). We also consider BIC, AIC, mAIC, KN, ACT and AGIC. Note that PC$_3$, ON$_2$, DDPA and BEMA are not included in this study as they are not suited for the small $p$, large $n$ scenario. In this simulation, we set $\x\sim\mathcal{N}_p(\0,\bSigma)$, where $\bSigma=\rho \Q\Q^\top+\I_p$. The matrix $\Q$ is obtained by first generating a $\D=(a_{ij})_{p\times k}$ matrix with independent $N(0,1)$ entries and then considering the QR decomposition of $\D=\Q_{p\times k}\S_{k\times k}$.
We set $p=12$, $k=3$, and vary $n$ from 100 to 500.  We set $\rho=2\delta((p-k/2+1/2)\log\log n/n)^{1/2}$, where $\delta$ increases from $1.5$ to $2.5$. This choice for $\rho$ is to keep the signal-to-noise ratio small, so to highlight the difference in performance when signal is weak.
We restrict the candidate spike size in the range of $k'\in\{0,1,\ldots,p-1\}$.

It is seen from Table \ref{tab:1} that, most methods indeed perform well when $n$ is large and the signal is strong, that is, $\delta$ is large. However, when $n$ is small or when $\delta$ is small, ILP and ILP$_{1/2}$ perform better than other methods in almost all settings. One exception is when $n=100$ and $\delta=1.5$. In this small $n$ and small $\delta$ setting, most methods underestimate whereas AIC, which tends to overestimate $k$, performs better in comparison.

\begin{table}[!t]
\setlength{\tabcolsep}{2pt}
\centering
{\renewcommand{\arraystretch}{0.75}
\begin{tabular}{c|ccc|cc|cc|cc|ccccccccc}
\hline
&&                           $\widehat{\mathbb{P}}(\hat{k}=k)$ & $\widehat{\mathbb{E}}(\hat{k})$ &      $\widehat{\mathbb{P}}(\hat{k}=k)$ & $\widehat{\mathbb{E}}(\hat{k})$ & $\widehat{\mathbb{P}}(\hat{k}=k)$ & $\widehat{\mathbb{E}}(\hat{k})$  &    $\widehat{\mathbb{P}}(\hat{k}=k)$ & $\widehat{\mathbb{E}}(\hat{k})$ &   $\widehat{\mathbb{P}}(\hat{k}=k)$ & $\widehat{\mathbb{E}}(\hat{k})$ &\\
\hline
&&\multicolumn{ 2}{c|}{$\delta=1.5$}&\multicolumn{ 2}{c|}{$\delta=1.75$}&\multicolumn{ 2}{c|}{$\delta=2$}& \multicolumn{ 2}{c|}{$\delta=2.25$}&\multicolumn{ 2}{c}{$\delta=2.5$}\\
\hline
\multirow{6}{*}{$n=100$}
&ILP        &   0.39  & 2.10&     0.69& 2.59&     0.83  & 2.80    &0.92  & 2.93 &     \textbf{0.98}& \textbf{2.98}\\
&ILP$_{1/2}$&   0.65  & 2.61&     \textbf{0.83}& \textbf{2.89}&     \textbf{0.91}  & \textbf{2.98}    &\textbf{0.94}  & \textbf{3.02} &     0.94& 3.02\\
&BIC        &   0.04  & 1.14&     0.15& 1.49&     0.32  & 1.91    &0.56  & 2.38 &     0.75& 2.67\\
&AIC        &   \textbf{0.73}  & \textbf{3.01}&     0.82& 3.14&     0.82  & 3.18    &0.84  & 3.18 &     0.83& 3.20\\
&mAIC       &   0.09  & 1.37&     0.28& 1.86&     0.53  & 2.31    &0.76  & 2.69 &     0.87& 2.85\\
&KN         &   0.41  & 2.41&     0.60& 2.59&     0.73  & 2.73    &0.86  & 2.89 &     0.92& 2.92\\
&ACT        &   0.04  & 1.40&     0.09& 1.64&     0.21  & 1.92    &0.31  & 2.09 &     0.43& 2.26\\
&AGIC       &   0.21  & 1.94&     0.33& 2.19&     0.49  & 2.45    &0.61  & 2.58 &     0.73& 2.72\\
\hline
\multirow{6}{*}{$n=200$}
&ILP        &   0.60  & 2.47&     0.83& 2.81&     0.95  & 2.95    &\textbf{0.99}  & \textbf{2.99} &     \textbf{1.00}& \textbf{3.00}\\
&ILP$_{1/2}$&   \textbf{0.84}  & \textbf{2.86}&     \textbf{0.93}& \textbf{2.98}&     \textbf{0.97}  & \textbf{3.03}    &0.97  & 3.03 &     0.97& 3.04\\
&BIC        &   0.03  & 1.41&     0.20& 1.62&     0.48  & 2.19    &0.77  & 2.70 &     0.93& 2.92\\
&AIC        &   0.83  & 3.11&     0.85& 3.15&     0.85 & 3.16     &0.84  & 3.19 &     0.84& 3.19\\
&mAIC       &   0.30  & 1.88&     0.64& 2.52&     0.86  & 2.83    &0.97  & 2.96 &     0.99& 2.99\\
&KN         &   0.59  & 2.59&     0.84& 2.84&     0.96  & 2.96    &\textbf{0.99}  & \textbf{2.99} &     \textbf{1.00}& \textbf{3.00}\\
&ACT        &   0.15  & 1.76&     0.33& 2.11&     0.51  & 2.38    &0.69  & 2.64 &     0.79& 2.77\\
&AGIC       &   0.41  & 2.31&     0.63& 2.61&     0.76  & 2.76    &0.87  & 2.87 &     0.93& 2.93\\
\hline
\multirow{6}{*}{$n=500$}
&ILP        &   0.79  & 2.75&     0.95& 2.95&     \textbf{1.00}  & \textbf{3.00}    &\textbf{1.00}  & \textbf{3.00} &     \textbf{1.00}& \textbf{3.00}\\
&ILP$_{1/2}$&   \textbf{0.95}  & \textbf{2.99}&     \textbf{0.98}& \textbf{3.02}&     0.98  & 3.02    &0.98  & 3.02 &     0.98& 3.02\\
&BIC        &   0.05  & 1.20&     0.34& 1.90&     0.94  & 2.94    &0.95  & 2.94 &     0.99& 2.99\\
&AIC        &   0.86  & 3.15&     0.81& 3.22&     0.85  & 3.16    &0.84  & 3.19 &     0.84& 3.19\\
&mAIC       &   0.68  & 2.56&     0.92& 2.91&     \textbf{1.00}  & \textbf{3.00}    &\textbf{1.00}  & \textbf{3.00} &     \textbf{1.00}& \textbf{3.00}\\
&KN         &   0.81  & 2.81&     0.95& 2.95&     0.99  & 2.99    &0.99  & 2.99 &     \textbf{1.00}& \textbf{3.00}\\
&ACT        &   0.41  & 2.26&     0.69& 2.64&     0.84  & 2.83    &0.92  & 2.92 &     0.96& 2.96\\
&AGIC       &   0.68  & 2.66&     0.88& 2.88&     0.95  & 2.95    &0.97  & 2.97 &     0.99& 2.99\\
\hline
\end{tabular}
}
\caption{Performances of ILP, ILP$_{1/2}$, BIC, AIC, mAIC, KN, ACT and AGIC when $p=12$ and $k=3$ in Simulation 1. Marked in boldface are those achieving the best evaluation criteria in each setting.}
\label{tab:1}
\end{table}

\medskip
\textbf{Simulation 2 (the standard spiked covariance model, large $p$).} 
In this simulation, we set $\x\sim\mathcal{N}_p(\0,\bSigma)$, where $\bSigma=\rho \Q\Q^\top+\I_p$. The matrix $\Q$ is obtained by first generating a $\D=(a_{ij})_{p\times k}$ matrix with independent $N(0,1)$ entries and then considering the QR decomposition of $\D=\Q_{p\times k}\S_{k\times k}$.
We set $k=10$ and restrict the candidate spike size in the range of $k'\in\{0,1,\ldots,20\}$.
We consider several different large $p$ settings including $p<n$, $p>n$ and $p=n$.
The performances of GIC in \eqref{eq:aic} and AGIC in \eqref{eq:aic5} are evaluated,  along with AIC, BCF, PC$_3$, ON$_2$  KN, BEMA, DPPA, ACT and AGIC.

It is seen from Table \ref{tab:2} that, except for a few settings,  GIC and AGIC have better performances than BEMA, DDPA and ACT when $\rho$ is relatively low, while all methods perform well when $\rho$ is large.  Even for these few settings, the performances of GIC and AGIC are also similar to those of competing methods. Moreover, DDPA tends to slightly over-select the number of components and this is consistent with the observation in \cite{Dobriban:Owen:2019}, which proposed a more conservative DDPA+ algorithm. ACT and BEMA's performances improve with $\rho$ but do not perform as well as GIC and AGIC, especially when $n$ is small.

We also note that, for the standard spiked covariance model, GIC outperforms AGIC. If we are sure that data are generated from a standard spiked covariance model, we recommend using GIC. In practice, we usually do not know whether data are generated from a standard spiked covariance model or a general spiked covariance model. In such cases, we recommend using AGIC. An additional simulation under the standard spiked covariance model with comparisons with BCF and ACT can be found in the Supplementary Materials.

\begin{table}[!t]
\setlength{\tabcolsep}{2.5pt}
\centering
{\renewcommand{\arraystretch}{0.75}
\begin{tabular}{c|c|cc|cc|cc|cc|ccccccccc}
\hline
\multirow{10}{*}{\begin{tabular}[c]{@{}c@{}}$n=200$\\$p=200$ \end{tabular}}
&&                           $\widehat{\mathbb{P}}(\hat{k}=k)$ & $\widehat{\mathbb{E}}(\hat{k})$ &      $\widehat{\mathbb{P}}(\hat{k}=k)$ & $\widehat{\mathbb{E}}(\hat{k})$ & $\widehat{\mathbb{P}}(\hat{k}=k)$ & $\widehat{\mathbb{E}}(\hat{k})$  &    $\widehat{\mathbb{P}}(\hat{k}=k)$ & $\widehat{\mathbb{E}}(\hat{k})$ &   $\widehat{\mathbb{P}}(\hat{k}=k)$ & $\widehat{\mathbb{E}}(\hat{k})$ &\\\cline{2-12}
&& \multicolumn{ 2}{c|}{$\rho=2$}&\multicolumn{ 2}{c|}{$\rho=3$}&\multicolumn{ 2}{c|}{$\rho=4$}&\multicolumn{ 2}{c|}{$\rho=5$}&\multicolumn{ 2}{c}{$\rho=6$}\\\hline
&AIC                       &0.00  & 3.65 &     0.25& 8.85&   0.96  & 9.96&     \textbf{1.00}& \textbf{10.0}&     \textbf{1.00}& \textbf{10.0}&\\
&GIC                      &0.00   & 6.80 &     \textbf{0.79}& \textbf{9.79}&   0.96  & 10.1&     0.96& 10.1&     0.96& 10.1&\\
&BCF                      &0.00  & 3.65 &     0.25& 8.85&   0.96  & 9.96&     \textbf{1.00}& \textbf{10.0}&     \textbf{1.00}& \textbf{10.0}&\\
&PC$_3$                   &0.00  & 7.77 &     0.76& 9.78&   0.99  & 10.0&     \textbf{1.00}  & \textbf{10.0}&     0.99  & 10.0&\\
&ON$_2$                   &0.05  & 11.0 &     0.11& 10.4&   0.21  & 10.3&     0.26  & 10.4&     0.31  & 9.86&\\
&KN                       &0.00  & 5.44 &     0.22& 9.04&   0.81  & 9.79&   \textbf{1.00}& \textbf{10.0}&     \textbf{1.00}& \textbf{10.0}&\\
&BEMA                     &0.00  & 4.78 &     0.00& 7.40&   0.70  & 9.71&   \textbf{1.00}& \textbf{10.0}&     \textbf{1.00}& \textbf{10.0}&\\
&DDPA                     &\textbf{0.27}  & \textbf{9.21} &     0.50& 10.7&   0.42  & 10.9&   0.38& 10.9&     0.37  & 11.0&\\
&ACT                      &0.00  & 4.32 &     0.36& 8.85&   0.97  & 9.97&     \textbf{1.00}& \textbf{10.0}&     \textbf{1.00}& \textbf{10.0}&\\
&AGIC                     &0.01  & 6.54 &     0.54& 9.46&   \textbf{1.00}  & \textbf{10.0}&     \textbf{1.00}& \textbf{10.0}&     \textbf{1.00}& \textbf{10.0}&\\
\hline
\multirow{10}{*}{\begin{tabular}[c]{@{}c@{}}$n=500$\\$p=200$ \end{tabular}}
& & \multicolumn{ 2}{c|}{$\rho=1$}&\multicolumn{ 2}{c|}{$\rho=1.5$}&\multicolumn{ 2}{c|}{$\rho=2$}&\multicolumn{ 2}{c|}{$\rho=2.5$}&\multicolumn{ 2}{c}{$\rho=3$}\\\hline
&AIC                       &0.00  & 3.45 &     0.10& 8.60&   0.99  & 9.99&     \textbf{1.00}& \textbf{10.0}&     \textbf{1.00}  & \textbf{10.0}\\
&GIC                       &0.00  & 4.75 &     \textbf{0.65}& \textbf{9.60}&   \textbf{1.00}  & \textbf{10.0}&   \textbf{1.00}& \textbf{10.0}&     \textbf{1.00}  & \textbf{10.0}&\\
&BCF                       &0.00  & 3.45 &     0.10& 8.60&   0.99  & 9.99&     \textbf{1.00}& \textbf{10.0}&     \textbf{1.00}  & \textbf{10.0}&\\
&PC$_3$                    &0.00  & 1.00 &     0.00& 1.00&   0.00  & 1.00&     0.00  & 2.64&     0.00  & 5.46&\\
&ON$_2$                    &0.05  & 10.9 &     0.11& 10.4&   0.22  & 10.7&     0.40  & 10.3&     0.56  & 10.1&\\
&KN                        &0.00  & 4.04 &     0.10& 8.66&   0.78  & 9.78&   \textbf{1.00}& \textbf{10.0}&     \textbf{1.00}& \textbf{10.0}&\\
&BEMA                      &0.00  & 3.10 &     0.00& 7.25&   0.40  & 9.41&     \textbf{1.00}  & \textbf{10.0}&     \textbf{1.00}& \textbf{10.0}&\\
&DDPA                      &\textbf{0.07}  & \textbf{7.79} &     0.40& 11.0&   0.31  & 11.2&     0.25& 11.3&     0.25  & 11.4&\\
&ACT                       &0.00  & 3.21 &     0.32& 8.65&   \textbf{1.00}  & \textbf{10.0}&     \textbf{1.00}  & \textbf{10.0}&     \textbf{1.00}  & \textbf{10.0}&\\
&AGIC                      &0.00  & 5.12 &     0.54& 9.44&   \textbf{1.00}  & \textbf{10.0}&     \textbf{1.00}  & \textbf{10.0}&     \textbf{1.00}  & \textbf{10.0}&\\
\hline
\multirow{10}{*}{\begin{tabular}[c]{@{}c@{}}$n=200$\\$p=500$ \end{tabular}}
&& \multicolumn{ 2}{c|}{$\rho=3$}&\multicolumn{ 2}{c|}{$\rho=4$}&\multicolumn{ 2}{c|}{$\rho=5$}&\multicolumn{ 2}{c|}{$\rho=6$}&\multicolumn{ 2}{c}{$\rho=7$}\\\hline
&AIC                       &0.00  & 1.36 &     0.00& 4.87&   0.00  & 7.80&   0.49  & 9.41&     \textbf{1.00}& \textbf{10.0}&      \\
&GIC                       &0.00  & 6.20 &     0.53& 9.43&   0.94  &9.94&   \textbf{1.00}& \textbf{10.0}&     \textbf{1.00}& \textbf{10.0}&     \\
&BCF                       &0.00  & 5.70 &     0.20& 9.10&   0.84  & 9.84&   \textbf{1.00}& \textbf{10.0}&      \textbf{1.00}& \textbf{10.0}& \\
&PC$_3$                    &0.00  & 1.00 &     0.00& 1.00&   0.00  & 1.10&     0.00  & 2.70&     0.00  & 4.47&\\
&ON$_2$                    &0.05  & 10.4 &     0.07& 10.2&   0.15  & 10.5&     0.20  & 10.0&     0.26  & 10.2&\\
&KN                        &0.00  & 5.74 &     0.05& 8.46&   0.29  & 9.12&   0.86& 9.86&     \textbf{1.00}& \textbf{10.0}&\\
&BEMA                      &0.00  & 5.50 &     0.10& 8.10&   0.20  & 9.10&     0.90& 9.90&     \textbf{1.00}& \textbf{10.0}&\\
&DDPA                      &\textbf{0.25}  & \textbf{9.50} &     0.42& 10.9&   0.29  & 11.2&     0.15& 11.3&     0.27  & 11.3&\\
&ACT                       &0.00  & 5.60 &     0.34& 8.86&   0.91  & 9.93&     \textbf{1.00}& \textbf{10.0}&     \textbf{1.00}& \textbf{10.0}&\\
&AGIC                      &0.02  & 7.27 &     \textbf{0.53}& \textbf{9.46}&   \textbf{0.95}  & \textbf{9.99}& \textbf{1.00}& \textbf{10.0}&     \textbf{1.00}& \textbf{10.0}&\\
\hline
\end{tabular}}
\caption{Performances of AIC, GIC, BCF, PC$_3$, ON$_2$, KN, BEMA, DDPA, ACT and AGIC in Simulation 2 with $k=10$. Marked in boldface are those achieving the best evaluation criteria in each setting.}
\label{tab:2}
\end{table}

\medskip

\textbf{Simulation 3 (the general spiked covariance model, large $p$).} 
In this simulation, we set $\x\sim\mathcal{N}_p(\0,\bSigma)$, where $\bSigma=\rho \Q\Q^\top+\Diag\{\nu_1^2,\ldots, \nu_p^2\}$, and $\nu_1^2,\ldots, \nu_p^2$ are i.i.d. from Uniform (0, 5). The matrix $\Q$ is obtained by first generating a $\D=(a_{ij})_{p\times k}$ matrix with independent $N(0,1)$ entries and then considering the QR decomposition of $\D=\Q_{p\times k}\S_{k\times k}$.
We set $k=10$ and restrict the candidate spike size in the range of $k'\in\{0,1,\ldots,20\}$.
We consider several different large $p$ settings including $p<n$, $p>n$ and $p=n$. For example, when $n=500$, $p=200$, $\rho=3$, after generating a $\bSigma$, we calculate its eigenvalues,  $\lambda_k=5.618$, $\lambda_{k+1}=4.843$, $\ldots$, $\lambda_{p-1}=0.1146$, $\lambda_{p}=0.04099$. Thus, $\text{SNR}=\lambda_k/\lambda_{k+1}-1=0.1601$.

The performance of AGIC in \eqref{eq:aic5} is evaluated,  along with PC$_3$, ON$_2$  KN, BEMA, DPPA and ACT. Methods such as AIC, GIC, and BCF are not included in this study as they are designed for the standard spiked covariance model.
It is seen from Table \ref{tab:3} that, in most settings, AGIC outperforms DDPA, ACT and BEMA.
Similar to our observations in Simulation 2, DDPA tends to slightly over-select the number of components; ACT and BEMA's performances improve with $\rho$ but do not perform as well as AGIC, especially when $n$ is small. 

\begin{table}[!t]
\setlength{\tabcolsep}{2.5pt}
\centering
{\renewcommand{\arraystretch}{0.75}
\begin{tabular}{c|c|cc|cc|cc|cc|ccccccccc}
\hline
\multirow{10}{*}{\begin{tabular}[c]{@{}c@{}}$n=200$\\$p=200$ \end{tabular}}
&&                           $\widehat{\mathbb{P}}(\hat{k}=k)$ & $\widehat{\mathbb{E}}(\hat{k})$ &      $\widehat{\mathbb{P}}(\hat{k}=k)$ & $\widehat{\mathbb{E}}(\hat{k})$ & $\widehat{\mathbb{P}}(\hat{k}=k)$ & $\widehat{\mathbb{E}}(\hat{k})$  &    $\widehat{\mathbb{P}}(\hat{k}=k)$ & $\widehat{\mathbb{E}}(\hat{k})$ &   $\widehat{\mathbb{P}}(\hat{k}=k)$ & $\widehat{\mathbb{E}}(\hat{k})$ &\\\cline{2-12}
&& \multicolumn{ 2}{c|}{$\rho=3$}&\multicolumn{ 2}{c|}{$\rho=4$}&\multicolumn{ 2}{c|}{$\rho=5$}&\multicolumn{ 2}{c|}{$\rho=6$}&\multicolumn{ 2}{c}{$\rho=7$}\\\hline
&PC$_3$                    &0.02  & 8.06 &     \textbf{0.33}& \textbf{9.28}&   \textbf{0.82}  & \textbf{9.89}&     0.90& 9.92&     0.96  & 10.0&\\
&ON$_2$                    &0.04  & 10.6 &     0.06& 10.9&   0.13  & 10.8&     0.18  & 10.8&     0.17  & 10.3&\\
&KN                        &\textbf{0.12}  & \textbf{8.18} &     0.27& 9.94&   0.26  & 10.6&     0.22& 11.4&     0.07& 12.1&\\
&BEMA                      &0.00  & 4.62 &     0.00& 6.05&   0.00  & 7.62&     0.05& 8.52&     0.21& 8.96&\\
&DDPA                      &0.00  & 0.97 &     0.00& 3.02&   0.04  & 6.39&     0.15& 13.1&     0.23& 15.7&\\
&ACT                       &0.00  & 5.39 &     0.05& 7.30&   0.41  & 9.01&     0.83& 9.79&     0.95& 9.95&\\
&AGIC                      &0.01  & 6.90 &     0.17& 8.53&   0.60  & 9.52&     \textbf{0.92}& \textbf{9.92}&     \textbf{0.97}& \textbf{9.98}&\\
\hline
\multirow{10}{*}{\begin{tabular}[c]{@{}c@{}}$n=500$\\$p=200$ \end{tabular}}
& & \multicolumn{ 2}{c|}{$\rho=1.5$}&\multicolumn{ 2}{c|}{$\rho=2$}&\multicolumn{ 2}{c|}{$\rho=2.5$}&\multicolumn{ 2}{c|}{$\rho=3$}&\multicolumn{ 2}{c}{$\rho=3.5$}\\\hline
&PC$_3$                    &0.00  & 1.00 &     0.00& 1.03&   0.00  & 1.34&     0.00  & 2.25&     0.00  & 3.08&\\
&ON$_2$                    &\textbf{0.07}  & \textbf{11.0} &     0.06& 11.1&   0.14  & 10.8&     0.17  & 11.0&     0.23  & 10.5&\\
&KN                        &0.00  & 18.9 &     0.00& 18.9&   0.00  & 18.9&     0.00& 19.0&       0.00& 19.0&\\
&BEMA                      &0.00  & 4.43 &     0.00& 6.00&   0.00  & 7.43&     0.14& 8.71&       0.65& 9.60&\\
&DDPA                      &0.00  & 0.39 &     0.00& 0.85&   0.00  & 10.2&     0.03& 74.5&       0.02& 142&\\
&ACT                       &0.00  & 5.01 &     0.06& 7.49&   0.47  & 9.16&     0.88  & 9.86&     0.99  & 9.99&\\
&AGIC                      &0.00  & 6.32 &     \textbf{0.13}& \textbf{8.44}&   \textbf{0.64}  & \textbf{9.56}&     \textbf{0.94}  & \textbf{9.94}&     \textbf{1.00}  & \textbf{10.0}&\\
\hline
\multirow{10}{*}{\begin{tabular}[c]{@{}c@{}}$n=200$\\$p=500$ \end{tabular}}
&& \multicolumn{ 2}{c|}{$\rho=4$}&\multicolumn{ 2}{c|}{$\rho=5$}&\multicolumn{ 2}{c|}{$\rho=6$}&\multicolumn{ 2}{c|}{$\rho=7$}&\multicolumn{ 2}{c}{$\rho=8$}\\\hline
&PC$_3$                    &0.00  & 1.00 &     0.00& 1.00&   0.00  & 1.23&     0.00  & 2.00&     0.00  & 3.03&\\
&ON$_2$                    &\textbf{0.06}  & \textbf{10.7} &     0.05& 10.5&   0.08  & 10.8&     0.10  & 10.6&     0.16  & 10.7&\\
&KN                        &0.03  & 7.34 &     0.10& 8.36&   0.27  & 9.24&     0.33& 10.1&       0.22& 11.0&\\
&BEMA                      &0.00  & 5.67 &     0.00& 6.90&   0.05  & 8.19&     0.22& 8.83&       0.50& 9.44&\\
&DDPA                      &0.00  & 0.90 &     0.00& 2.09&   0.00  & 3.81&     0.00& 5.61&       0.08& 7.73&\\
&ACT                       &0.01  & 6.48 &     0.21& 8.44&   0.56  & 9.42&     0.91& 9.90&     0.97& 9.98&\\
&AGIC                      &0.05  & 7.83 &     \textbf{0.37}& \textbf{9.10}&   \textbf{0.74}  & \textbf{9.75}&     \textbf{0.95}& \textbf{10.0}&     \textbf{0.99}& \textbf{10.0}&\\
\hline
\end{tabular}}
\caption{Performances of PC$_3$, ON$_2$, KN, BEMA, DDPA, ACT and AGIC in Simulation 3 with $k=10$. Marked in boldface are those achieving the best evaluation criteria in each setting.}
\label{tab:3}
\end{table}

\medskip
\textbf{Simulation 4  (the general spiked covariance model, large $p$, comparison with ACT).}
In this simulation, we set $\x\sim\mathcal{N}_p(\0,\bSigma)$, where $\bSigma=\rho \Q\Q^\top+\Diag\{\nu_1^2,\ldots, \nu_p^2\}$, $\rho=5\sqrt{p/n}$, and $\nu_1^2,\ldots, \nu_p^2$ are i.i.d. from Uniform (0, 5). The matrix $\Q$ is obtained by first generating a $\D=(a_{ij})_{p\times k}$ matrix with independent $N(0,1)$ entries and then considering the QR decomposition of $\D=\Q_{p\times k}\S_{k\times k}$.
We set $k=10$ and restrict the candidate spike size in the range of $k'\in\{0,1,\ldots,20\}$. In this simulation, for example, when $n=500$, $p=200$, after generating a $\bSigma$, we calculate its eigenvalues,  $\lambda_k=5.737$, $\lambda_{k+1}=4.920$, ..., $\lambda_{p-1}=0.07826$, $\lambda_{p}=0.05324$. Thus, $\text{SNR}=\lambda_k/\lambda_{k+1}-1=0.1660$.

\begin{figure}[!t]
\centering
\includegraphics[width=\textwidth]{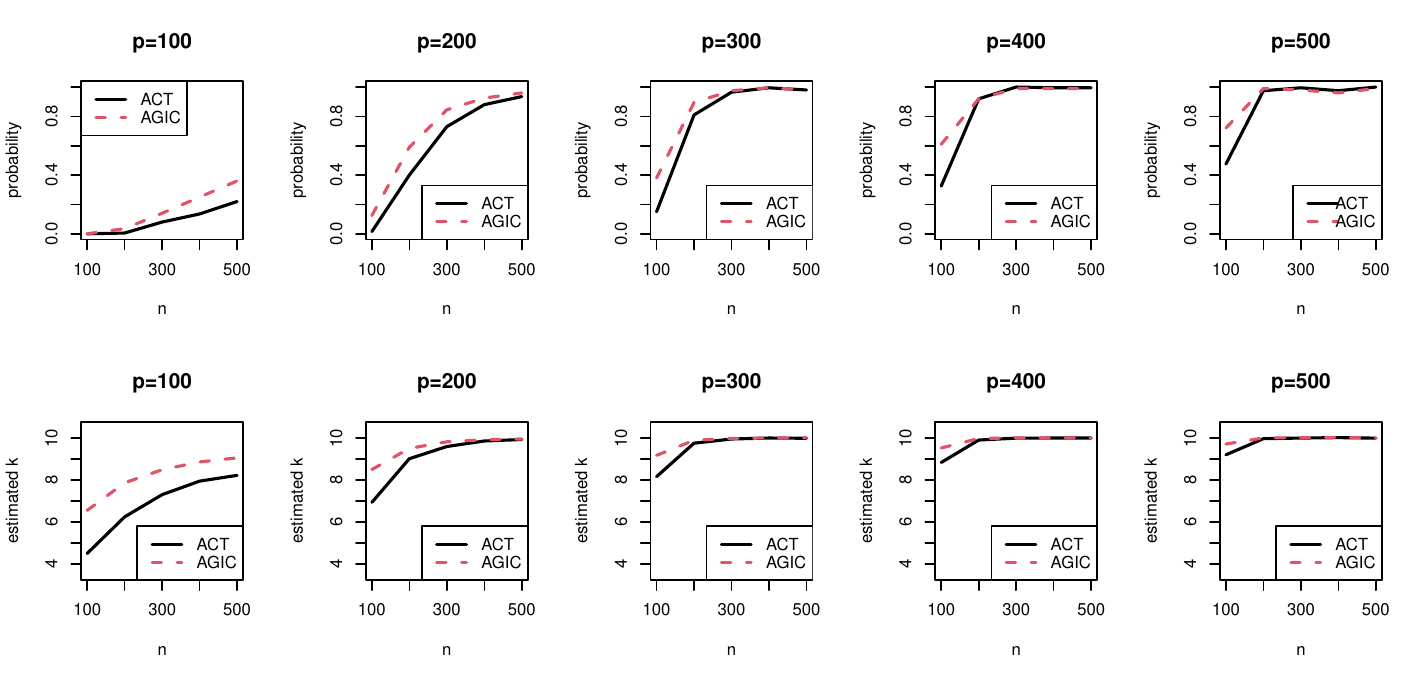}
\caption{Performances of ACT and AGIC in Simulation 4 under varying $p,n$ with $k=10$.}
\label{fig:1}
\end{figure}

 We compare AGIC with ACT. It is seen from Figure \ref{fig:1} that,   AGIC performs better than   ACT, especially when $p$ or $n$ is relatively small (i.e., $p\leq 300$ or $n\leq 300$). An additional simulation under the general spiked covariance model with comparison with ACT can be found in the Supplementary Materials.

\section{Real-world data examples}
\label{section:empirical}
In this section, we consider three real-world data examples ranging from finance to biology, and compare with existing methods including the PC$_3$ estimator proposed in \cite{Bai:Ng:2002}, the ON$_2$ estimator proposed in \cite{Onatski:2009}, the KN estimator proposed in \cite{Kritchman:Nadler:2008}, the DDPA estimator proposed in \cite{Dobriban:Owen:2019}, the ACT estimator proposed in \cite{Fan:Guo:Zheng:2022} and the BEMA estimator proposed in \cite{Ke:Ma:Lin:2023}. It is known that real-world data examples are more complicated than simulated data and may not be strictly generated from the spiked model. However, there is a tuning parameter $\gamma$ in GIC and AGIC, which makes GIC and AGIC very flexible. We may set, for example, $\gamma=2\varphi(1+\sqrt{c})$ and $\gamma=2\log (1+\sqrt{c})$ for GIC in \eqref{eq:aic} and AGIC in \eqref{eq:aic5}, respectively.

\subsection{The Fama-French 100 portfolios}
We estimate the number of factors using the excess returns of Fama-French 100 portfolios. The well-known risk factors for equity markets are well explained by $k=3$ Fama-French factors: the market factor, the size factor and the value factor. We use the daily returns of 100 industrial Portfolios from January 2, 2010 to April 30, 2019. The dimension and the sample size of the data set are  $p=100$ and $n=2346$, respectively.

When estimating $k$, AGIC selects $\hat{k}=3$, while GIC, PC$_3$, ON$_2$, KN, DDPA,  BEMA and ACT estimate the number of factors as 47, 26, 1, 57, 100, 6 and 3, respectively. It is seen that both AGIC and ACT correctly estimate the correct number of factors, while ON$_2$ and BEMA slightly underestimate and overestimate the number of factors, respectively,  and GIC, PC$_3$, KN and DDPA more notably overestimate the number of factors. Compared with AGIC, GIC does not correctly estimate the correct number of factors. One possible reason is that this data set follows more closely to a general spiked covariance model than to a standard spiked covariance model.

In this data set, the 10 largest eigenvalues of the sample covariance matrix are 130.90, 5.48, 3.12, 1.53, 1.31, 0.99, 0.70, 0.61, 0.55, 0.51. The variance explained by 3 factors is $85.64\%$ due to a large spike top eigenvalue. The 10 largest eigenvalues of the sample correlation matrix are 80.62, 3.22, 2.06, 0.82,  0.75, 0.56, 0.37, 0.35, 0.31, 0.27 and 3 factors explain $85.90\%$ total variation in 100 portfolios, while $\hat{k}=1$ and 6 explain $80.62\%$ and $88.02\%$ total variation, respectively.

\subsection{The 1000 Genomes project genotypes}
Next, we illustrate our method on a genotype data from the 1000 Genomes Project (Phase III), publicly available at \url{https://www.internationalgenome.org}. We used Plink to retain common variants with minor allele frequency greater than 0.1, and generated a set of variants in approximate linkage disequilibrium. The data pre-processing steps are performed following \citet{zhong2022empirical}. The dimension and the sample size of the data set are  $p=1000$ and $n=2504$, respectively.

Amongst these subjects, there are five self-reported ethnicity groups, including African, Caucasian, East Asian, Hispanic and South Asian. As this data matrix follows a clustering model, the covariance matrix has $k=k_0-1$ spiked eigenvalues, where $k_0$ is the number of clusters.  We treat these five self-reported groups as the ground truth, which gives $k=4$. When estimating $k$, both GIC and AGIC select $\hat{k}=4$, while PC$_3$, ON$_2$, KN, DDPA, BEMA and ACT estimate $k$ as 4, 4, 11, 12, 8 and 8, respectively. It is seen that GIC, AGIC, PC$_3$ and ON$_2$ correctly estimate the correct number of spikes, while KN, DDPA, BEMA and ACT overestimate the number of spikes.

In this data set, the 10 largest eigenvalues of the sample covariance matrix are
70.81, 32.68, 10.81, 8.82, 2.79, 2.63, 2.53, 2.45, 2.43, 2.39. The variance explained by $\hat{k}=4$ is $12.24\%$ in 1000 variables. In contrast, the 10 largest eigenvalues of the sample correlation matrix are 64.25, 30.34, 10.40, 8.43, 2.91, 2.66, 2.56, 2.52, 2.42, 2.40 and $\hat{k}=4$ explains $12.20\%$ total variation in 1000 variables, while $\hat{k}=8$, 11 and 12 explain $12.40\%$, $13.12\%$ and $13.36\%$ total variation, respectively.

\subsection{The Peripheral blood mononuclear cells}
\label{section:empirical3}
Finally, we illustrate our method on the single cell RNA-sequencing data of Peripheral blood mononuclear cells (PBMC) from 10X Genomics\footnote{\url{https://cf.10xgenomics.com/samples/cell/pbmc3k/pbmc3k_filtered_gene_bc_matrices.tar.gz.}}. The dimension and the sample size of the data set are  $p=13711$ and $n=2626$, respectively.

For the PBMC cells, there are 6 major cell types, including T cells, B cells, Monocytes, Natural Killer cells, Dendritic cells, and Platelet cells.
When estimating $k$, both GIC and AGIC select $\hat{k}=6$, while PC$_3$, ON$_2$, KN, DDPA, BEMA and ACT estimate the number of factors as 2, 9, 217, 2626, 7 and 20, respectively. For BEMA, we sub-sampled 1000 genes due to its extreme computational demand.
It is seen that both GIC and AGIC correctly estimate the correct number of spikes, while ON$_2$, KN, DDPA, BEMA and ACT overestimate the number of spikes and PC$_3$ underestimates the number of spikes.

In this data set, the 10 largest eigenvalues of the sample covariance matrix are 93.47, 47.27, 33.21, 20.00, 19.10, 19.01, 16.32, 16.05, 15.84, 15.26. The variance explained by $\hat{k}=6$  is $1.692\%$ in 13711 variables. In contrast, the 10 largest eigenvalues of the sample correlation matrix are 93.47, 47.27, 33.21, 20.00, 19.10, 19.01, 16.32, 16.05, 15.84, 15.26 and $\hat{k}=6$ explains $1.692\%$ total variation in 13711 variables, while $\hat{k}=2$, 7 and 9 explain $1.026\%$, $1.811\%$ and $2.044\%$ total variation, respectively.

\section*{Supplementary Materials} The supplementary materials include all technical proofs and two additional simulation studies.

\section*{Acknowledgments} The authors thank the editor, the associate editor, and the referees for their valuable comments and helpful advice.

\section*{Funding} Hu is partially supported by the National Natural Science Foundation of China (Grant Nos. 12171187, 12371261). Guo was partially supported by the Key Program of the National Natural Science Foundation of China (Grant No. 12431009) and the National Key Research and Development Program of China (Grant No. 2020YFA0714100).

\section*{Disclosure Statement} The authors report there are no competing interests to declare.

\newpage
\renewcommand{\thesection}{S\arabic{section}}
\renewcommand{\thelemma}{S\arabic{lemma}}
\renewcommand{\theequation}{S\arabic{equation}}
\setcounter{section}{0}
\setcounter{lemma}{0}
\setcounter{equation}{0}
\setcounter{page}{1}

\begin{center}
{\Large\bf Supplementary Materials for ``Limiting laws and consistent estimation criteria for fixed and diverging number of spiked eigenvalues"}
\\
\medskip
{\large Jianwei Hu, Jingfei Zhang, Jianhua Guo and Ji Zhu}
\end{center}

\section{Auxiliary Lemmas}
\begin{lemma}[\cite{bickel2008regularized}]
\label{lemma:spectral}
Suppose a mean zero random vector $\Y=(Y_1,\cdots, Y_p)^\top$ satisfies $\mathbb{E}|Y_{j}|^\beta\le C$, for $\beta>2$ and some constant $C>0$, and has a covariance matrix $\bSigma$. Let $\Y^{(1)},\cdots, \Y^{(n)}$ be $n$ independent copies of $\Y$. Then there exists some universal constant $c>0$ such that the sample covariance matrix of $\Y^{(1)},\cdots, \Y^{(n)}$, denoted as $\S_n$, satisfies
\[
\mathbb{P}(\|\S_n-\bSigma\|>t)\leq cn^{-\beta/4}(2p+1)p(t/p)^{-\beta/2}.
\]
Here $\|\cdot\|$ denotes the spectral norm.
\end{lemma}

\begin{lemma}[Weyl's inequality in \cite{Weyl:1912}]
\label{lemma:weyl}
Let $\A$ and $\B$ be $p\times p$ Hermitian matrices with eigenvalues ordered as $\lambda_1\geq \cdots\geq\lambda_p$ and $\mu_1\geq\cdots\geq\mu_p$, respectively. It holds that
\[
\sup_{1\leq j\leq p}\vert\mu_j-\lambda_j\vert\leq\|\A-\B\|,
\]
where $\|\cdot\|$ denotes the spectral norm.
\end{lemma}

\begin{lemma}
\label{lemma:1}
Suppose that Assumptions \ref{ass1}, \ref{ass22} hold, $\mathbb{E}Y_{ij}^4<\infty$ and $p$ is fixed while $n\rightarrow\infty$. It holds that
$$(d_i-\lambda_i)^2=\mathcal{O}_p(1/n),\,\,1\le i\le k,$$
$$ (p-k)(\bar{d}_{k+1}-\lambda)^2=\mathcal{O}_p(1/n).$$
\end{lemma}

\begin{lemma}\label{lemma:4}
Consider $\log L_k$ and $\log\tilde L_k$ as defined in \eqref{eq:L} and \eqref{eq:tildeL}, respectively. It holds that
\begin{equation*}
\log L_k=\log\tilde  L_k+\mathcal{O}_p(1).
\end{equation*}
\end{lemma}

Let $\X=[\x_1,\ldots,\x_n]^\top\in\mathbb{R}^{n\times p}$ denote $n$ samples from a multivariate distribution, where $\x_{i}=\bSigma^{\frac{1}{2}}\y_{i}$ and $\Y=[\y_1,\ldots,\y_n]^\top\in\mathbb{R}^{n\times p}$. By spectral decomposition, we may express $\bSigma$ as
\begin{equation*}
\bSigma=\bGamma\bLambda\bGamma^\top,\quad \bLambda= \begin{pmatrix}
 \bLambda_1 \,\,& \0 \\
  \0 \,\,&  \bLambda_2
\end{pmatrix},
\end{equation*}
where $\bLambda_1=\Diag\{\lambda_1,\ldots,\lambda_k\}$, $\bLambda_2=\Diag\{\lambda_{k+1},\ldots,\lambda_p\}$, and $\bGamma=(\bGamma_1, \bGamma_2)$ is a $p\times p$ orthogonal matrix collecting the eigenvectors of $\bSigma$. Define $\Z=[\z_1,\ldots,\z_n]^\top\in\mathbb{R}^{n\times p}$, where $\z_i=\bGamma^\top\x_i$.  Write $\Z=[\Z_1,\Z_2]$ with $\Z_1\in\mathbb{R}^{n\times k}$ and $\Z_2\in\mathbb{R}^{n\times (p-k)}$. It is seen that $\Z=\Y\bGamma\bLambda^{\frac{1}{2}}$, $\Z_1=\Y\bGamma_1\bLambda_1^{\frac{1}{2}}$ and $\Z_2=\Y\bGamma_2\bLambda_2^{\frac{1}{2}}$.
Under this setting, we state two lemmas.

\begin{lemma}\label{lemma:large deviation1}
For any positive semi-definite matrix $\A\in\mathbb{R}^{n\times n}$, denote by $l_1\geq\cdots\geq l_n$ the eigenvalues of $\A$.
Suppose that $Y_{ij}$'s are independent and identically distributed normal random variables with mean zero and unit variance, i.e.,  $Y_{ij}\sim \mathcal{N}(0,1)$.
There exist some universal constant $C>0$ and some constant $\rho$ such that the following inequality holds
\begin{eqnarray*}
&&\mathbb{P}\left(\left\Vert\frac{1}{n}\bLambda_{1}^{-\frac{1}{2}}\Z_1^\top(\I+\A)\Z_1\bLambda_{1}^{-\frac{1}{2}}-\frac{1}{n}\tr(\I+\A)\I\right\Vert>t\right)\\
&\leq&2\exp\left\{Ck-\frac{nt^2}{\rho\sum_{s=1}^n(1+l_s)^2/n}\right\},
\end{eqnarray*}
for all $0<t<\rho\sum_{s=1}^n(1+l_s)^2/\{4C_1n(1+l_1)\}$, where $C_1$ is a constant depending on $C$. Consequently, the following also holds
$$
\mathbb{P}\left(\left\Vert\frac{1}{n}\bLambda_{1}^{-\frac{1}{2}}\Z_1^\top\A\Z_1\bLambda_{1}^{-\frac{1}{2}}-\frac{1}{n}\tr(\A)\I\right\Vert>t\right)
\leq2\exp\left\{Ck-\frac{nt^2}{\rho\sum_{s=1}^nl_s^2/n}\right\},
$$
for all $0<t<\rho\sum_{s=1}^nl_s^2/(4C_1nl_1)$.
\end{lemma}

\begin{lemma}\label{corollary:large deviation2}
For any positive semi-definite matrix $\A$, denote by $l_1\geq\cdots\geq l_n$ the eigenvalues of $\A$. Suppose that $Y_{ij}$'s are independent and identically distributed random variables with $\mathbb{E}(Y_{ij})=0$, $\mathrm{Var}(Y_{ij})=1$, and $\mathbb{E}Y_{ij}^4=\eta<\infty$, then we have
\[
\left\Vert\frac{1}{n}\bLambda_{1}^{-\frac{1}{2}}\Z_1^\top(\I+\A)\Z_1\bLambda_{1}^{-\frac{1}{2}}-\frac{1}{n}\tr(\I+\A)\I\right\Vert=\mathcal{O}_p(\bar lk/\sqrt{n}),
\]
and
\[
\tr\left[\frac{1}{n}\bLambda_{1}^{-\frac{1}{2}}\Z_1^\top(\I+\A)\Z_1\bLambda_{1}^{-\frac{1}{2}}-\frac{1}{n}\tr(\I+\A)\I\right]^2=\mathcal{O}_p(\bar l^2k^2/n),
\]
where $\bar l=\{\sum_{s=1}^n(1+l_s)^2/n\}^{1/2}$.
Consequently, the following also holds
\[
\left\Vert\frac{1}{n}\bLambda_{1}^{-\frac{1}{2}}\Z_1^\top\A\Z_1\bLambda_{1}^{-\frac{1}{2}}-\frac{1}{n}\tr(\A)\I\right\Vert=\mathcal{O}_p(\bar lk/\sqrt{n}),
\]
and
\[
\tr\left[\frac{1}{n}\bLambda_{1}^{-\frac{1}{2}}\Z_1^\top\A\Z_1\bLambda_{1}^{-\frac{1}{2}}-\frac{1}{n}\tr(\A)\I\right]^2=\mathcal{O}_p(\bar l^2k^2/n).
\]
\end{lemma}

The proofs of Lemmas \ref{lemma:1}-\ref{corollary:large deviation2} are deferred to Section \ref{sec:proofaux}.

\section{Proofs of main results}\label{sec:proof}

\subsection{Proof of Theorem \ref{thm1}}
\label{sec:proof1}
In what follows, we give the proof for the case of $\mathbb{E}Y_{ij}^4<\infty$.  That is, by Lemma \ref{corollary:large deviation2},  we can show that Theorem \ref{thm1} holds for $k=o(n^{1/4})$. For the case of $Y_{ij}$ is Gaussian,  instead of using Lemma \ref{corollary:large deviation2},  we can show Theorem \ref{thm1} holds for $k=o(n^{1/3})$ using Lemma \ref{lemma:large deviation1}. The proof is very similar to the case of $\mathbb{E}Y_{ij}^4<\infty$, and the details are thus omitted.

By definition, each $d_i$ solves the equation
\[
0=|d \I-\S_n|= |d \I-\S_{22}| \cdot |d \I-\mathbf K_n(d)|,
\]
where
\[
\mathbf K_n(d)=\frac{1}{n}\Z_1^\top(\I+\H_n(d))\Z_1,
\]
\[
\H_n(d)=\frac{1}{n}\Z_2(d \I-\S_{22})^{-1}\Z_2^\top.
\]
It is known that the spectrum of $\S_{22}$ goes inside the support of $F^{c,H}$ in probability (see, for example, \cite{Bai:Yao:2012}).
Next, we demonstrate that the following determinant equation has $k$ solutions outside the support of $F^{c,H}$ in probability,
\[
|d \I-\mathbf K_n(d)|=0.
\]
Rewriting the determinant equation, we have
\[
\left| d\I-\frac{1}{n}\Z_1^\top(\I+\H_n(d))\Z_1\right|=0,
\]
or equivalently,
\begin{equation}\label{thm1:eq11}
\left| \bLambda_{1}^{-1}d\I-\frac{1}{n}\bLambda_{1}^{-\frac{1}{2}}\Z_1^\top(\I+\H_n(d))\Z_1\bLambda_{1}^{-\frac{1}{2}}\right|=0.
\end{equation}
We proceed to show that \eqref{thm1:eq11} has $k$ solutions outside the support of $F^{c,H}$ in probability.

Recall $c_{n}=(p-k)/n$ and $F^{c_{n},H_n}(x)$ is defined from $F^{c,H}(x)$ with $c$ and $H$ replaced by $c_{n}$ and $H_n$, respectively. According to \citet{Bai:Silverstein:2004}, we have
\begin{equation}\label{Gaussian:convergence:rate}
\int f(x)dF_n(x)=\int f(x)dF^{c_{n},H_n}(x)+\mathcal{O}_p\left(1/n\right),
\end{equation}
where $f(x)$  is an analytic function and $\sup_x \vert f(x)\vert<\infty$.

Recall $\{\beta_j\}_{k+1\le j\le p}$ are the eigenvalues of $\S_{22}=\frac{1}{n}\Z_2^\top\Z_2$. Then, we have
\begin{equation}\label{thm1:eq101}
\begin{aligned}
ds_n(d)&=d\frac{1}{p-k}\sum_{j=k+1}^{p}\frac{1}{\beta_j-d}=d\int\frac{1}{x-d}dF_n(x)\\
&=d\int\frac{1}{x-d}dF^{c_n,H_n}(x)+\mathcal{O}_p(1/n)\\
&=ds(d)+\mathcal{O}_p(1/n),
\end{aligned}
\end{equation}
where the last two equalities hold due to \eqref{Gaussian:convergence:rate}.
By \eqref{thm1:eq101}, we have
\begin{equation}\label{thm1:eq120}
\begin{aligned}
d\underline{s}_n(d)&=-(1-c_n)+c_nds_n(d)\\
&=-(1-c_n)+c_nds(d)+\mathcal{O}_p(1/n)\\
&=d\underline{s}(d)+\mathcal{O}_p(1/n).
\end{aligned}
\end{equation}

Recall ${m_1(d)}=\int\frac{x}{d-x}dF^{c_n,H_n}(x)$ and ${m_2(d)}=\int\frac{x^2}{(d-x)^2}dF^{c_n,H_n}(x)$. Also, we have
\begin{equation}\label{thm1:eq20}
\begin{aligned}
\frac{1}{n}\tr(\H_n(d))&=\frac{1}{n}\tr\left\{\frac{1}{n}\Z_2(d \I-\S_{22})^{-1}\Z_2^\top\right\}=\frac{1}{n}\tr\{(d \I-\S_{22})^{-1}\S_{22}\}\\
&=\frac{1}{n}\sum_{j=k+1}^{p}\frac{\beta_j}{d-\beta_j}=\frac{p-k}{n}\frac{1}{p-k}\sum_{j=k+1}^{p}\frac{\beta_j}{d-\beta_j}\\
&=c_n\int\frac{x}{d-x}dF_n(x)\\
&=c_nm_1(d)+\mathcal{O}_p(1/n),
\end{aligned}
\end{equation}
where the last equality holds due to \eqref{Gaussian:convergence:rate}.

Moreover, we have that
\begin{equation}\label{thm1:eq3}
\begin{aligned}
\frac{1}{n}\tr(\H_n^2(d))&=\frac{1}{n}\tr\left\{(d \I-\S_{22})^{-1}\S_{22}(d \I-\S_{22})^{-1}\S_{22}\right\}\\
&=\frac{1}{n}\sum_{j=k+1}^{p}\frac{\beta_j^2}{(d-\beta_j)^2}=\frac{p-k}{n}\frac{1}{p-k}\sum_{j=k+1}^{p}\frac{\beta_j^2}{(d-\beta_j)^2}\\
&=c_n\int\frac{x^2}{(d-x)^2}dF_n(x)\\
&=c_nm_2(d)+\mathcal{O}_p(1/n),
\end{aligned}
\end{equation}
where the last equality holds due to \eqref{Gaussian:convergence:rate}.

It is seen that
\begin{equation}\label{thm1:eq210}
\int\frac{x}{d-x}dF_n(x)=\mathcal{O}_p(1),\,\,\,\int\frac{x^2}{(d-x)^2}dF_n(x)=\mathcal{O}_p(1).
\end{equation}
Denote the eigenvalues of $\H_n(d)$ as $l_1\geq\cdots\geq l_n$. Combining \eqref{thm1:eq20}, \eqref{thm1:eq3} and \eqref{thm1:eq210}, we have
\[
\frac{1}{n}\sum_{s=1}^n(1+l_s)^2=\left\{1+\frac{2}{n}\tr(\H_n(d))+\frac{1}{n}\tr(\H_n^2(d))\right\}=1+\mathcal{O}_p(1).
\]
By Lemma \ref{corollary:large deviation2} ($\H_n(d)$ is positive semi-definite when $d$ is outside the support of $F^{c,H}$), we have,
\[
\left\Vert\frac{1}{n}\bLambda_{1}^{-\frac{1}{2}}\Z_1^\top(\I+\H_n(d))\Z_1\bLambda_{1}^{-\frac{1}{2}}-\frac{1}{n}\tr(\I+\H_n(d))\I\right\Vert=\mathcal{O}_p(k/\sqrt{n}).
\]
As a result, we have the following spectral decomposition
\begin{equation}\label{thm1:eq40}
\begin{aligned}
&\frac{1}{n}\tr(\I+\H_n(d))\I-\frac{1}{n}\bLambda_{1}^{-\frac{1}{2}}\Z_1^\top(\I+\H_n(d))\Z_1\bLambda_{1}^{-\frac{1}{2}}\\
=&\V_n(d)\Diag\{v_1(d),\ldots,v_k(d)\}\V_n(d)^\top,
\end{aligned}
\end{equation}
where $\V_n(d)$ is an orthogonal matrix and $\max_{1\leq j\leq k}\vert v_j(d)\vert=\mathcal{O}_p(k/\sqrt{n})$.
Moreover, we have
\begin{eqnarray}\label{thm1:eq70}
&&\bLambda_{1}^{-1}d\I-\frac{1}{n}\bLambda_{1}^{-\frac{1}{2}}\Z_1^\top(\I+\H_n(d))\Z_1\bLambda_{1}^{-\frac{1}{2}}\\\nonumber
&=&\bLambda_{1}^{-1}d\I-\frac{1}{n}\tr(\I+\H_n(d))\I+\left\{\frac{1}{n}\tr(\I+\H_n(d))\I-\frac{1}{n}\bLambda_{1}^{-\frac{1}{2}}\Z_1^\top(\I+\H_n(d))\Z_1\bLambda_{1}^{-\frac{1}{2}}\right\}\\\nonumber
&=&\bLambda_{1}^{-1}d\I-\frac{1}{n}\tr(\I+\H_n(d))\I+\M_n(d),
\end{eqnarray}
where
\begin{equation}\label{thm1:eq6}
\M_n(d)=\frac{1}{n}\tr(\I+\H_n(d))\I-\frac{1}{n}\bLambda_{1}^{-\frac{1}{2}}\Z_1^\top(\I+\H_n(d))\Z_1\bLambda_{1}^{-\frac{1}{2}}.
\end{equation}
Noting that $\bLambda_{1}^{-1}=\Diag\left\{\frac{1}{\lambda_1},\ldots,\frac{1}{\lambda_k}\right\}$
and then \eqref{thm1:eq11} becomes
\begin{eqnarray*}
&&\left|\bLambda_{1}^{-1}d\I-\frac{1}{n}\bLambda_{1}^{-\frac{1}{2}}\Z_1^\top(\I+\H_n(d))\Z_1\bLambda_{1}^{-\frac{1}{2}}\right|\\
&=&\left|\Diag\left\{\frac{d-\frac{1}{n}\tr(\I+\H_n(d))\lambda_1}{\lambda_1},\ldots,\frac{d-\frac{1}{n}\tr(\I+\H_n(d))\lambda_k}{\lambda_k}\right\}+\M_n(d)\right|=0.
\end{eqnarray*}
Let $\hat{f}_{i}\triangleq\hat{f}_{i}(d)$'s,  $1\leq i\leq k$, be eigenvalues of the following matrix,
\[
\A_n(d)=\Diag\left\{\frac{d-\frac{1}{n}\tr(\I+\H_n(d))\lambda_1}{\lambda_1},\ldots,\frac{d-\frac{1}{n}\tr(\I+\H_n(d))\lambda_k}{\lambda_k}\right\}+\M_n(d).
\]
Then, $d_1,\ldots,d_k$ are solutions of
\begin{equation}
\label{thm1:eq7}
\hat{f}_i(d)=0,\,\,\, 1\leq i\leq k.
\end{equation}
Let $f_i\triangleq f_i(d)$'s, $1\leq i\leq k$, be eigenvalues of the following matrix,
\[
\F_n(d)=\Diag\left\{\frac{d-\frac{1}{n}\tr(\I+\H_n(d))\lambda_1}{\lambda_1},\ldots,\frac{d-\frac{1}{n}\tr(\I+\H_n(d))\lambda_k}{\lambda_k}\right\}.
\]
Then, we have
\begin{equation*}
\vert \A_n(d)\vert =\vert \F_n(d)+\M_n(d)\vert,
\end{equation*}
where $\F_n(d)=\Diag\{f_1,\ldots, f_k\}.$
Combining \eqref{thm1:eq40} and \eqref{thm1:eq6}, we have
\begin{equation}\label{thm1:eq108}
\begin{aligned}
\M_n(d)=\V_n(d)\Diag\left\{v_1(d),\ldots,v_k(d))\right\}\V_n(d)^\top,
\end{aligned}
\end{equation}
where $\V_n(d)$ is an orthogonal matrix and $\max_{1\leq j\leq k}\vert v_j(d) \vert=\mathcal{O}_p(k/\sqrt{n})$.
Hence,
\begin{equation}
\label{thm1:eq8}
\tr(\M_n^2(d))=\tr\left[\Diag\left\{v_1^2(d),\ldots,v_k^2(d)\right\}\right]=\mathcal{O}_p(k^3/n).
\end{equation}
By Lemma \ref{lemma:weyl}, we have
\begin{equation}
\label{thm1:eq55}
\max_{1\leq i\leq k}|\hat{f}_i-f_i|\leq||\M_n(d)||=\mathcal{O}_p(k/\sqrt{n}).
\end{equation}
Then, we have
\begin{equation}
\tilde{\F}_n(d)\triangleq\Diag\{f_1+\tilde{\epsilon}_1,\ldots,f_k+\tilde{\epsilon}_k\}=\tilde{\V}_n(d)^\top(\F_n(d)+\M_n(d))\tilde{\V}_n(d),
\end{equation}
where $\tilde{\V}_n(d)$ is an orthogonal matrix and  $\max_{1\leq j\leq k}|\tilde{\epsilon}_j|=\mathcal{O}_p(k/\sqrt{n})$.
Thus, we have
\begin{equation}\label{eqn:de04}
\vert \A_n(d)\vert=\vert \tilde{\F}_n(d)\vert=\prod_{1\leq j\leq k}(f_j+\tilde{\epsilon}_j),
\end{equation}
where $\max_{1\leq j\leq k}|\tilde{\epsilon}_j|=\mathcal{O}_p(k/\sqrt{n})$ and $\sum_{1\leq j\leq k}|\tilde{\epsilon}_j|=\mathcal{O}_p(k^2/\sqrt{n})$.

Correspondingly, for all $d$ with the property that $\min_{1\leq t\leq k}|f_t|$ is uniformly bounded away from zero, we have
\[
|\sum_{1\leq t\leq k}\frac{1}{f_t}\tilde{\epsilon}_t|\leq\frac{1}{\min_{1\leq t\leq k}|f_t|}\sum_{1\leq t\leq k}|\tilde{\epsilon}_t|=\mathcal{O}_p(k^2/\sqrt{n}),
\]
\[
|\sum_{1\leq t\leq k}\sum_{t< s\leq k}\frac{1}{f_tf_s}\tilde{\epsilon}_t\tilde{\epsilon}_s|\leq\frac{1}{\min_{1\leq t\leq k,t< s\leq k}|f_tf_s|}\sum_{1\leq t\leq k}\sum_{t< s\leq k}|\tilde{\epsilon}_t\tilde{\epsilon}_s|=\mathcal{O}_p(k^4/n).
\]
Under the condition $k=o(n^{1/4})$, we have
\begin{equation}\label{thm1:eq9}
\begin{aligned}
\vert \A_n(d)\vert=\prod_{1\leq j\leq k}(f_j+\tilde{\epsilon}_j)&=\prod_{1\leq j\leq k}f_j+\prod_{1\leq j\leq k}f_j\sum_{1\leq t\leq k}\frac{1}{f_t}\tilde{\epsilon}_t+\prod_{1\leq j\leq k}f_j\sum_{1\leq t\leq k}\sum_{t< s\leq k}\frac{1}{f_tf_s}\tilde{\epsilon}_t\tilde{\epsilon}_s+\cdots\\
&=\prod_{1\leq j\leq k}f_j+\mathcal{O}_p(k^2/\sqrt{n})\prod_{1\leq j\leq k}f_j+\mathcal{O}_p(k^4/n)\prod_{1\leq j\leq k}f_j+\cdots\\
&=(1+\mathcal{O}_p(k^2/\sqrt{n}))\prod_{1\leq j\leq k}f_j=(1+o_p(1))\vert \F_n(d)\vert.
\end{aligned}
\end{equation}
Note that both $\A_n(d)$ and $\F_n(d)$ are rational functions of $d$. The asymptotic equivalence \eqref{thm1:eq9} holds for all $d$ in an open set $\U\subset \mathbb{R}$ (where
$\min_{1\leq t\leq k}|f_t|$ is bounded away from zero). Since $\U$ has accumulation points and the error term $o_p(1)$ is uniformly on compact subsets of $\U$,  by the uniqueness theorem for rational functions, for all $d$ outside the support of $F^{c,H}$, we have $\vert \A_n(d)\vert=(1+o_p(1))\vert \F_n(d)\vert$. By \eqref{thm1:eq7}, $d_1,\ldots,d_k$ tend to the solutions of
\begin{equation*}
\frac{d-\frac{1}{n}\tr(\I+\H_n(d))\lambda_i}{\lambda_i}=0,\,\,\, 1\leq i\leq k.
\end{equation*}

Next, consider the solution of
\begin{equation}
\label{thm1:eq103}
\frac{d-\frac{1}{n}\tr(\I+\H_n(d))\lambda_i}{\lambda_i}=0.
\end{equation}
By \eqref{thm1:eq20}, for all $1\leq i\leq k$, we have
\begin{equation}\label{thm1:eq50}
\frac{1}{n}\tr(\I+\H_n(d_i))=1+c_nm_1(d_i)+\mathcal{O}_p(k/n).
\end{equation}
By \eqref{thm1:eq103} and \eqref{thm1:eq50}, for all $1\leq i\leq k$, we have
\begin{equation}
\label{thm1:eq109}
\frac{d_i}{\lambda_i}-(1+c_nm_1(d_i))=o_p(1).
\end{equation}
By the definitions of $s(d)$ and $\underline{s}(d)$, it is seen that $m_1(d)=-1-ds(d)$. Thus, we have
\begin{equation}
\label{thm1:eq104}
1+c_nm_1(d)=-d\underline{s}(d).
\end{equation}
By \eqref{thm1:eq109} and \eqref{thm1:eq104}, for all $1\leq i\leq k$, we have
\begin{equation}\label{thm1:eq140}
\frac{d_i}{\lambda_i}+d_i\underline{s}(d_i)=o_p(1).
\end{equation}
That is, for all $1\leq i\leq k$, we have
\begin{equation}
\label{thm1:eq110}
\frac{-\underline{s}^{-1}(d_i)}{\lambda_i}-1=o_p(1).
\end{equation}
Recall $\underline{s}(d)=-\frac{1-c_n}{d}+c_ns(d)$. It is known that (see, for example, \cite{Bai:Silverstein:1998} and \cite{Bai:Silverstein:2004})
\begin{equation*}
\begin{aligned}
d&=-\frac{1}{\underline{s}(d)}+c_n\int\frac{t}{1+t\underline{s}(d)}dH_n(t)\\
&=-\underline{s}^{-1}(d)+c_n\int\frac{-\underline{s}^{-1}(d)t}{-\underline{s}^{-1}(d)-t}dH_n(t).
\end{aligned}
\end{equation*}
That is,
\begin{equation}\label{thm1:eq107}
\begin{aligned}
\frac{d_i}{\lambda_i}&=\frac{-\underline{s}^{-1}(d_i)}{\lambda_i}
+c_n\int\frac{\frac{-\underline{s}^{-1}(d_i)}{\lambda_i}t}{\frac{-\underline{s}^{-1}(d_i)}{\lambda_i}\lambda_i-t}dH_n(t).
\end{aligned}
\end{equation}
Recall
\[
\psi_n(x)=x+c_n\int\frac{xt}{x-t}dH_n(t).
\]
Combining \eqref{thm1:eq110} and \eqref{thm1:eq107}, for all $1\leq i\leq k$, we have
\begin{equation}\label{thm1:eq111}
\frac{d_i-\psi_n(\lambda_i)}{\lambda_i}=o_p(1).
\end{equation}
Thus, for all  $1\leq i\leq k$, we have
\begin{equation*}
\frac{d_i}{\psi_n(\lambda_i)}\stackrel{p}{\longrightarrow}1.
\end{equation*}

For $1\leq i\leq k$, we have
\begin{equation}\label{thm1:eq201}
\hat{s}_{n}(d_i)-s_{n}(d_i)=\frac{1}{p-k}\sum_{j=k+1}^p\frac{1}{d_j-d_i}-\frac{1}{p-k}\sum_{j=k+1}^p\frac{1}{\beta_j-d_i}=\frac{1}{p-k}\sum_{j=k+1}^p\frac{\beta_j-d_j}{(d_j-d_i)(\beta_j-d_i)}.
\end{equation}
By the interlacing inequality in \citet{Horn:Johnson:2013}, for $j=k+1,\ldots, p$, we have
\[
d_j\leq \beta_j\leq d_{j-k}.
\]
Then, we have
\[
0\leq\frac{1}{p-k}\sum_{j=k+1}^p(\beta_j-d_j)\leq \frac{1}{p-k}\sum_{j=k+1}^{2k}\beta_j-\frac{1}{p-k}\sum_{j=p-k+1}^pd_j<\frac{k\beta_{k+1}}{p-k}.
\]
Combining this result with \eqref{thm1:eq201}, we have
\[
0\leq d_i\hat{s}_{n}(d_i)-d_is_{n}(d_i)<\frac{k\beta_{k+1}d_i}{(p-k)(d_i-\beta_{k+1})(d_i-d_{k+1})}.
\]
That is, for all $1\leq i\leq k$, we have
\[
0\leq d_i\hat{s}_{n}(d_i)-d_is_{n}(d_i)<\frac{k\beta_{k+1}d_k}{(p-k)(d_k-\beta_{k+1})(d_k-d_{k+1})}=\mathcal{O}_p(k/n).
\]
By the definitions of $\underline{s}_{n}(d_i)$ and $\hat{\underline{s}}_{n}(d_i)$, for all $1\leq i\leq k$, we have
\begin{equation}\label{thm1:eq202}
d_i\hat{\underline{s}}_{n}(d_i)-d_i\underline{s}_{n}(d_i)=\mathcal{O}_p(k/n).
\end{equation}
By \eqref{thm1:eq120} and \eqref{thm1:eq202}, for all $1\leq i\leq k$, we have
\begin{equation}\label{thm1:eq203}
d_i\hat{\underline{s}}_{n}(d_i)-d_i\underline{s}(d_i)=\mathcal{O}_p(k/n).
\end{equation}
Combining this result with \eqref{thm1:eq140}, for all $1\leq i\leq k$, we have
\begin{equation*}
\frac{d_i}{\lambda_i}+d_i\hat{\underline{s}}_n(d_i)=o_p(1).
\end{equation*}
That is, for all  $1\leq i\leq k$, we have
\begin{equation*}
\frac{-\hat{\underline{s}}_n^{-1}(d_i)}{\lambda_i}-1=o_p(1).
\end{equation*}

\eop

\subsection{Proof of Theorem \ref{thm2}}
In what follows, we only give the proof for the case of $\mathbb{E}Y_{ij}^4<\infty$.  That is, by Lemma \ref{lemma:diagonal} (ii),  we can show that Theorem \ref{thm2} (ii) holds for $k=o(n^{1/4})$. For the case of $Y_{ij}$ is Gaussian,  instead of using Lemma \ref{lemma:diagonal} (ii),  we can show that Theorem \ref{thm2} (i) holds for $k=o(n^{1/3})$ using Lemma \ref{lemma:diagonal} (i). The proof is very similar to the case of $\mathbb{E}Y_{ij}^4<\infty$, and the details are thus omitted.

For $1\leq i\leq k$, according to \eqref{thm1:eq108},  we have
\begin{equation*}
\M_n(d_i)_{ii}=\sum_{l=1}^kv_l(d_i)V_{il}^2=\mathcal{O}_p(k/\sqrt{n}).
\end{equation*}
Combining this with the result in Lemma \ref{lemma:diagonal}, we have
\begin{equation}
\label{thm2:eq9}
\frac{d_i-\frac{1}{n}\tr(\I+\H_n(d_i))\lambda_i}{\lambda_i}=\mathcal{O}_p(k/\sqrt{n}).
\end{equation}
By \eqref{thm1:eq20}, \eqref{thm1:eq104} and \eqref{thm2:eq9}, we have
\begin{equation}\label{thm2:eq30}
\frac{d_i}{\lambda_i}+d_i\underline{s}(d_i)=\mathcal{O}_p(k/\sqrt{n}).
\end{equation}
That is,
\begin{equation}
\label{thm2:eq11}
\frac{-\underline{s}^{-1}(d_i)}{\lambda_i}-1=\mathcal{O}_p(k/\sqrt{n}).
\end{equation}
Combining \eqref{thm1:eq107} and \eqref{thm2:eq11}, for $1\leq i\leq k$, we have
\begin{equation*}
\frac{d_i-\psi_n(\lambda_i)}{\lambda_i}=\mathcal{O}_p(k/\sqrt{n}).
\end{equation*}
By \eqref{thm1:eq120} and \eqref{thm2:eq30}, we have
\begin{equation}
\label{thm2:eq12}
\frac{-\underline{s}_n^{-1}(d_i)-\lambda_i}{\lambda_i}=\mathcal{O}_p(k/\sqrt{n}).
\end{equation}
Finally, by \eqref{thm1:eq203} and \eqref{thm2:eq30}, we have
\[
\frac{-\hat{\underline{s}}_n^{-1}(d_i)-\lambda_i}{\lambda_i}=\mathcal{O}_p(k/\sqrt{n}).
\]

\eop

\subsection{Proof of Theorem \ref{thm222}}
By \eqref{eqn:thm22206}, we have
\[
f_i=\frac{d_i-\frac{1}{n}\tr(\I+\H_n(d_i))\lambda_i}{\lambda_i}=\mathcal{O}_p(k/\sqrt{n}).
\]
Similar to the proof of Theorem \ref{thm2}, for $1\leq i\leq k$, we have
\[
\frac{d_i-\psi_n(\lambda_i)}{\lambda_i}=\mathcal{O}_p(k/\sqrt{n}),
\]
\[
\frac{-\hat{\underline{s}}_n^{-1}(d_i)-\lambda_i}{\lambda_i}=\mathcal{O}_p(k/\sqrt{n}).
\]
\eop

\subsection{Proof of Theorem \ref{thm3}}
In what follows, we only give the proof for the case of $\mathbb{E}Y_{ij}^4<\infty$.  That is, by  Theorem \ref{thm2} (ii),  we can show that Theorem \ref{thm2} (ii) holds for $k=o(n^{1/5})$. For the case of $Y_{ij}$ is Gaussian, instead of using Theorem \ref{thm2} (ii), we can show that Theorem \ref{thm3} (i) holds for $k=o(n^{1/3})$ using Theorem \ref{thm2} (i). The proof is similar to the case of $\mathbb{E}Y_{ij}^4<\infty$, and the details are thus omitted.

Recall each $d_i$ solves the equation
\[
0=|d \I-\S_n|=|d \I-\S_{22}| \cdot |d \I-\mathbf K_n(d)|,
\]
where
\[
\mathbf K_n(d)=\frac{1}{n}\Z_1^\top(\I+\H_n(d))\Z_1,
\]
\[
\H_n(d)=\frac{1}{n}\Z_2(d \I-\S_{22})^{-1}\Z_2^\top.
\]
By Theorem \ref{thm1}, we have that $\{d_i\}_{1\le i\le k}$ go outside the support of $F^{c,H}$ in probability, and it is known that the spectrum of $\S_{22}$ goes inside the support of $F^{c,H}$ in probability.
Therefore, $d_i$'s solve the determinant equation
\[
|d \I-\mathbf K_n(d)|=0,
\]
which can be rewritten as
\[
\left| d_i\I-\frac{1}{n}\Z_1^\top(\I+\H_n(d_i))\Z_1\right|=0.
\]
Note that the random matrix $\mathbf K_n(d_i)$ can be decomposed as follows
\begin{equation}\label{thm3:eq8}
\begin{aligned}
&\mathbf K_n(d_i)=\frac{1}{n}\Z_1^\top(\I+\H_n(d_i))\Z_1\\
=&\frac{1}{n}\left\{\Z_1^\top(\I+\H_n(d_i))\Z_1-\tr(\I+\H_n(d_i))\bLambda_{1}\right\}+\frac{1}{n}\tr(\I+\H_n(d_i))\bLambda_{1}\\
=&\frac{1}{\sqrt{n}}\mathbf R_n(d_i)+\frac{1}{n}\tr(\I+\H_n(d_i))\bLambda_{1},
\end{aligned}
\end{equation}
where $\mathbf R_n(d_i)=\{\Z_1^\top(\I+\H_n(d_i))\Z_1-\tr(\I+\H_n(d_i))\bLambda_{1}\}/\sqrt{n}$.
By Theorem \ref{thm2}, we have
\[
\eta_i\triangleq\frac{\sqrt{n}(d_i-\psi_n(\lambda_i))}{\lambda_i}=\mathcal{O}_p(k).
\]
Then we have
\begin{equation}
\label{thm3:eq1}
d_i \I-\mathbf K_n(d_i)=\left\{\psi_n(\lambda_i)+\lambda_i\eta_i/\sqrt{n}\right\}\I-\mathbf K_n(\psi_n(\lambda_i))-\left\{\mathbf K_n(d_i)-\mathbf K_n(\psi_n(\lambda_i))\right\}.
\end{equation}
Further, treating $\A=\left\{\psi_n(\lambda_i)+\lambda_i\eta_i/\sqrt{n}\right\}\I-\S_{22}$, $\B=\psi_n(\lambda_i)\I-\S_{22}$, and using the fact that $\A^{-1}-\B^{-1}=\A^{-1}(\B-\A)\B^{-1}$, we have
\begin{equation}\label{thm3:eq2}
\begin{aligned}
&\bLambda_{1}^{-\frac{1}{2}}\left\{\mathbf K_n(d_i)-\mathbf K_n(\psi_n(\lambda_i))\right\}\bLambda_{1}^{-\frac{1}{2}}\\
=&\frac{1}{n^2}\bLambda_{1}^{-\frac{1}{2}}\Z_1^\top\Z_2\left[\left\{(\psi_n(\lambda_i)+\lambda_i\eta_i/\sqrt{n})\I-\S_{22}\right\}^{-1}-(\psi_n(\lambda_i)\I-\S_{22})^{-1}\right]\Z_2^\top\Z_1\bLambda_{1}^{-\frac{1}{2}}\\
=&-\frac{\eta_i\lambda_i}{n^2\sqrt{n}}\bLambda_{1}^{-\frac{1}{2}}\Z_1^\top\Z_2\left[\left\{(\psi_n(\lambda_i)+\lambda_i\eta_i/\sqrt{n})\I-\S_{22}\right\}^{-1}(\psi_n(\lambda_i)\I-\S_{22})^{-1}\right]\Z_2^\top\Z_1\bLambda_{1}^{-\frac{1}{2}}.
\end{aligned}
\end{equation}
Denote
\begin{equation}\label{thm3:eq3}
\D_n=\frac{1}{n}\lambda_i^2\Z_2\left[\left\{(\psi_n(\lambda_i)+\lambda_i\eta_i/\sqrt{n})\I-\S_{22}\right\}^{-1}(\psi_n(\lambda_i)\I-\S_{22})^{-1}\right]\Z_2^\top.
\end{equation}
Recall $m_3(d)=\int\frac{x}{(d-x)^2}dF^{c_n,H_n}$. It holds that
\begin{equation}\label{thm3:eq4}
\begin{aligned}
\frac{1}{n}\tr(\D_n)=&\frac{1}{n}\lambda_i^2\tr\left[\left\{(\psi_n(\lambda_i)+\lambda_i\eta_i/\sqrt{n})\I-\S_{22}\right\}^{-1}(\psi_n(\lambda_i)\I-\S_{22})^{-1}\S_{22}\right]\\
=&\frac{1}{n}\lambda_i^2\sum_{j=k+1}^{p}\frac{\beta_j}{(\psi_n(\lambda_i)+\lambda_i\eta_i/\sqrt{n}-\beta_j)(\psi_n(\lambda_i)-\beta_j)}\\
=&\lambda_i^2\frac{p-k}{n}\frac{1}{p-k}\sum_{j=k+1}^{p}\frac{\beta_j}{(\psi_n(\lambda_i)-\beta_j)^2}\frac{\psi_n(\lambda_i)-\beta_j}{(\psi_n(\lambda_i)+\lambda_i\eta_i/\sqrt{n}-\beta_j)}\\
=&\lambda_i^2\frac{p-k}{n}\frac{1}{p-k}\sum_{j=k+1}^{p}\frac{\beta_j}{(\psi_n(\lambda_i)-\beta_j)^2}(1+\mathcal{O}_p(\eta_i/\sqrt{n}))\\
=&\lambda_i^2\frac{p-k}{n}\int\frac{x}{(\psi_n(\lambda_i)-x)^2}dF_n(x)(1+\mathcal{O}_p(\eta_i/\sqrt{n}))\\
=&\lambda_i^2c_nm_3(\psi_n(\lambda_i))(1+\mathcal{O}_p(\eta_i/\sqrt{n}))(1+\mathcal{O}_p({1}/{n}))\\
=&\lambda_i^2c_nm_3(\psi_n(\lambda_i))(1+\mathcal{O}_p(k/\sqrt{n})),
\end{aligned}
\end{equation}
where the last two equalities hold due to \eqref{Gaussian:convergence:rate}.
Noting that
\begin{eqnarray*}
\lambda_i^2m_3(\psi_n(\lambda_i))&=&\int\frac{\lambda_i}{(\psi_n(\lambda_i)-x)}\frac{\lambda_ix}{(\psi_n(\lambda_i)-x)}dF^{c_n,H_n}(x)\\
&\leq&\frac{\lambda_i}{\psi_n(\lambda_i)-(1+\sqrt{c})^2}\int\frac{\lambda_ix}{\psi_n(\lambda_i)-x}dF^{c_n,H_n}(x)=\mathcal{O}(1),
\end{eqnarray*}
we have
$
\frac{1}{n}\tr(\D_n)=\mathcal{O}_p(1).
$

Define $m_4(d)=\int\frac{x^2}{(d-x)^4}dF^{c_n,H_n}(x)$. Using a similar argument as in \eqref{thm3:eq4}, we have
\begin{eqnarray*}
&&\frac{1}{n}\tr(\D_n^2)
=\frac{1}{n}\lambda_i^4\sum_{j=k+1}^{p}\frac{\beta_j^2}{(\psi_n(\lambda_i)+\lambda_i\eta_i/\sqrt{n}-\beta_j)^2(\psi_n(\lambda_i)-\beta_j)^2}\\
&=&\lambda_i^4\frac{p-k}{n}\frac{1}{p-k}\sum_{j=k+1}^{p}\frac{\beta_j^2}{(\psi_n(\lambda_i)+\lambda_i\eta_i/\sqrt{n}-\beta_j)^2(\psi_n(\lambda_i)-\beta_j)^2}\\
&=&\lambda_i^4\frac{p-k}{n}\int\frac{x^2}{(\psi_n(\lambda_i)-x)^4}dF_n(x)(1+\mathcal{O}_p(\eta_i/\sqrt{n}))\\
&=&\lambda_i^4 c_nm_4(\psi_n(\lambda_i))(1+\mathcal{O}_p(k/\sqrt{n})).
\end{eqnarray*}
Noting that
\begin{eqnarray*}
\lambda_i^4m_4(\psi_n(\lambda_i))
&=&\int\frac{\lambda_i^2}{(\psi_n(\lambda_i)-x)^2}\frac{\lambda_i^2x^2}{(\psi_n(\lambda_i)-x)^2}dF^{c_n,H_n}(x)\\
&\leq&\frac{\lambda_i^2}{(\psi_n(\lambda_i)-(1+\sqrt{c})^2)^2}\int\frac{\lambda_i^2x^2}{(\psi_n(\lambda_i)-x)^2}dF^{c_n,H_n}(x)=\mathcal{O}(1),
\end{eqnarray*}
we also have $\frac{1}{n}\tr(\D_n^2)=\mathcal{O}_p(1)$.

Let $\tau_1\geq\cdots\geq \tau_n$ denote the eigenvalues of $\D_n$. We have
\[
\frac{1}{n}\sum_{s=1}^n\tau_s^2=\frac{1}{n}\tr(\D_n^2)=\mathcal{O}_p(1).
\]
By Lemma \ref{corollary:large deviation2}, we have
\begin{equation}\label{thm3:eq5}
\left\Vert\frac{1}{n}\bLambda_{1}^{-\frac{1}{2}}\Z_1^\top\D_n\Z_1\bLambda_{1}^{-\frac{1}{2}}-\frac{1}{n}\tr(\D_n)\I\right\Vert=\mathcal{O}_p(k/\sqrt{n}).
\end{equation}
Combining \eqref{thm3:eq2}, \eqref{thm3:eq3}, \eqref{thm3:eq4} and \eqref{thm3:eq5}, we have the following spectral decomposition
\begin{equation}\label{thm3:eq6}
\begin{aligned}
&\bLambda_{1}^{-\frac{1}{2}}\{\mathbf K_n(d_i)-\mathbf K_n(\psi_n(\lambda_i))\}\bLambda_{1}^{-\frac{1}{2}}\\
=&-\frac{\lambda_i\eta_i}{\sqrt{n}}c_nm_3(\psi_n(\lambda_i))(1+\mathcal{O}_p(k/\sqrt{n}))\I-\frac{\eta_i}{\sqrt{n}\lambda_i}\W(d_i)\Diag\{w_1(d_i),\ldots,w_k(d_i)\}\W(d_i)^\top,
\end{aligned}
\end{equation}
where $\W(d_i)$ is orthogonal and $\max_{1\leq j\leq k}\vert w_j(d_i)\vert=\mathcal{O}_p(k/\sqrt{n})$.

By \eqref{thm1:eq20}, we have
\begin{equation}\label{thm3:eq30}
\frac{1}{n}\tr(\I+\H_n(\psi_n(\lambda_i)))=1+c_nm_1(\psi_n(\lambda_i))+\mathcal{O}_p(1/n).
\end{equation}

Recall $\underline{s}(d)=-\frac{1-c_n}{d}+c_ns(d)$. It is known that (see, for example, \cite{Bai:Silverstein:1998} and \cite{Bai:Silverstein:2004})
\begin{equation*}
\begin{aligned}
d&=-\frac{1}{\underline{s}(d)}+c_n\int\frac{t}{1+t\underline{s}(d)}dH_n(t)\\
&=-\underline{s}^{-1}(d)+c_n\int\frac{-\underline{s}^{-1}(d)t}{-\underline{s}^{-1}(d)-t}dH_n(t).
\end{aligned}
\end{equation*}
Let $d=\psi_n(\lambda_i)$. Then, we have
\begin{equation}\label{thm3:eq022}
\begin{aligned}
\psi_n(\lambda_i)&=-\underline{s}^{-1}(\psi_n(\lambda_i))+c_n\int\frac{-\underline{s}^{-1}(\psi_n(\lambda_i))t}{-\underline{s}^{-1}(\psi_n(\lambda_i))-t}dH_n(t).
\end{aligned}
\end{equation}
Recall
\begin{equation}\label{thm3:eq023}
\psi_n(\lambda_i)=\lambda_i+c_n\int\frac{\lambda_it}{\lambda_i-t}dH_n(t).
\end{equation}
Combining \eqref{thm3:eq022} and \eqref{thm3:eq023}, we have (see also, \citet{Bai:Ding:2012} and \cite{Zhang:Zheng:Pan:Zhong:2022})
\begin{equation}\label{thm3:eq10}
\lambda_i\underline{s}(\psi_n(\lambda_i))=-1.
\end{equation}
By \eqref{thm1:eq104}, we have
\[
1+c_nm_1(\psi_n(\lambda_i))=-\psi_n(\lambda_i)\underline{s}(\psi_n(\lambda_i)).
\]
That is,
\[
1+c_nm_1(\psi_n(\lambda_i))=-\frac{\psi_n(\lambda_i)}{\lambda_i}\lambda_i\underline{s}(\psi_n(\lambda_i)).
\]
Combining this result with \eqref{thm3:eq10}, we have
$m_1(\psi_n(\lambda_i))=\int\frac{t}{\lambda_i-t}dH_n(t)$. Combining this result with \eqref{thm3:eq30}, we have
\begin{equation}\label{thm3:eq7}
\frac{1}{n}\tr(\I+\H_n(\psi_n(\lambda_i)))=1+c_n\int\frac{t}{\lambda_i-t}dH_n(t)+\mathcal{O}_p(1/n).
\end{equation}
Combining \eqref{thm3:eq8}, \eqref{thm3:eq1}, \eqref{thm3:eq6} and \eqref{thm3:eq7}, we have
\begin{equation}\label{thm3:eq9}
\begin{aligned}
&\bLambda_{1}^{-\frac{1}{2}}(d_i \I-\mathbf K_n(d_i))\bLambda_{1}^{-\frac{1}{2}}\\
=&\bLambda_{1}^{-1}\psi_n(\lambda_i)\I+\frac{\lambda_i\eta_i}{\sqrt{n}}\bLambda_{1}^{-1}-[1+c_n\int\frac{t}{\lambda_i-t}dH_n(t)+\mathcal{O}_p(1/n)]\I\\
&\hspace{0.5in}-\frac{1}{\sqrt{n}}\bLambda_{1}^{-\frac{1}{2}}\mathbf R_n(\psi_n(\lambda_i))\bLambda_{1}^{-\frac{1}{2}}+\frac{\lambda_i\eta_i}{\sqrt{n}}c_nm_3(\psi_n(\lambda_i))(1+\mathcal{O}_p(k/\sqrt{n}))\I\\
&\hspace{0.5in}+\frac{\eta_i}{\sqrt{n}\lambda_i}\W(d_i)\Diag\{w_1(d_i),\ldots,w_k(d_i)\}\W(d_i)^\top.
\end{aligned}
\end{equation}

Noting that $\frac{\psi_n(\lambda_i)}{\lambda_j}-\frac{\psi_n(\lambda_i)}{\lambda_i}=\psi_n(\lambda_i)(\frac{1}{\lambda_j}-\frac{1}{\lambda_i})$
and $\bLambda_{1}^{-1}=\Diag\left\{\frac{1}{\lambda_1},\ldots,\frac{1}{\lambda_k}\right\}$, we can reorganize \eqref{thm3:eq9} as
$$
\bLambda_{1}^{-\frac{1}{2}}(d_i \I-\mathbf K_n(d_i))\bLambda_{1}^{-\frac{1}{2}}=\B_1+\B_2+\B_3+\B_4,
$$
where
\begin{eqnarray*}
\B_1&=&\Diag\left\{\psi_n(\lambda_i)(\frac{1}{\lambda_1}-\frac{1}{\lambda_i}),\ldots,\psi_n(\lambda_i)(\frac{1}{\lambda_k}-\frac{1}{\lambda_i})\right\}\\
&&+\frac{\eta_i\lambda_i}{\sqrt{n}}\Diag\left\{\frac{1}{\lambda_1},\ldots,\frac{1}{\lambda_k}\right\},\\
\B_2&=&\mathcal{O}_p(1/n)\I+\frac{\eta_i\lambda_i}{\sqrt{n}}c_nm_3(\psi_n(\lambda_i))(1+\mathcal{O}_p(k/\sqrt{n}))\I,\\
\B_3&=&-\frac{1}{\sqrt{n}}\bLambda_{1}^{-\frac{1}{2}}\mathbf R_n(\psi_n(\lambda_i))\bLambda_{1}^{-\frac{1}{2}},\\
\B_4&=&\frac{\eta_i}{\sqrt{n}\lambda_i}\W(d_i)\Diag\left\{w_1(d_i),\ldots,w_k(d_i)\right\}\W(d_i)^\top.
\end{eqnarray*}
Since $\psi_n(\lambda_i)(\frac{1}{\lambda_i}-\frac{1}{\lambda_i})=0$, it is seen that the $i$-th diagonal element of $\B_1$ is $\eta_i/\sqrt{n}$. The $i$-th diagonal element of $\B_2$ is
\[
\mathcal{O}_p(1/n)+\frac{\eta_i}{\sqrt{n}}\lambda_ic_nm_3(\psi_n(\lambda_i))(1+\mathcal{O}_p(k/\sqrt{n})).
\]
As $\W(d_i)$ is an orthogonal matrix and $\max_{1\leq j\leq k}\vert w_j(d_i)\vert=\mathcal{O}_p(k/\sqrt{n})$, the $i$-th diagonal element of $\B_4$ is $\frac{\eta_i}{\lambda_i}\mathcal{O}_p(k/n)$.
Therefore, the $i$-th diagonal element of $\bLambda_{1}^{-\frac{1}{2}}(d_i \I-\mathbf K_n(d_i))\bLambda_{1}^{-\frac{1}{2}}$ can be written as
\begin{eqnarray*}
&&\frac{\eta_i}{\sqrt{n}}\{1+c_nm_3(\psi_n(\lambda_i))\lambda_i(1+\mathcal{O}_p(k/\sqrt{n}))\}-\left\{\frac{1}{\sqrt{n}}\bLambda_{1}^{-\frac{1}{2}}\mathbf R_n(\psi_n(\lambda_i))\bLambda_{1}^{-\frac{1}{2}}\right\}_{ii}\\
&&+\frac{\eta_i}{\lambda_i}\mathcal{O}_p(k/n)+\mathcal{O}_p(1/n)=\mathcal{O}_p(k^{5/2}/n),
\end{eqnarray*}
where the equality is true by Lemma \ref{lemma:diagonal}, as $\A_n(d_i)=\bLambda_{1}^{-\frac{1}{2}}(d_i \I-\mathbf K_n(d_i))\bLambda_{1}^{-\frac{1}{2}}$.
Recall
\[
\eta_i=\frac{\sqrt{n}(d_i-\psi_n(\lambda_i))}{\lambda_i}=\mathcal{O}_p(k).
\]
Thus, under the condition $k=o(n^{1/5})$, we have
\[
\eta_i(1+c_nm_3(\psi_n(\lambda_i))\lambda_i)-[\bLambda_{1}^{-\frac{1}{2}}\mathbf R_n(\psi_n(\lambda_i))\bLambda_{1}^{-\frac{1}{2}}]_{ii}=o_p(1).
\]
That is,
\begin{equation}\label{thm3:eq14}
\frac{\eta_i}{\mathrm{T}_n(\lambda_i)}=1+o_p(1),
\end{equation}
where
\[
\mathrm{T}_n(\lambda_i)\triangleq\frac{1}{1+c_nm_3(\psi_n(\lambda_i))\lambda_i}\left[\bLambda_{1}^{-\frac{1}{2}}\mathbf R_n(\psi_n(\lambda_i))\bLambda_{1}^{-\frac{1}{2}}\right]_{ii}.
\]
Recall $m_3(d)=-m_1'(d)$. Combining \eqref{thm1:eq104} and \eqref{thm3:eq10}, we have
\begin{equation}\label{thm3:eq11}
1+c_nm_1(d)=d/\lambda_i,
\end{equation}
where $d=\psi_n(\lambda_i)$.
By \eqref{thm3:eq11}, we have
\[
c_nm_1'(d)d'=\frac{d'}{\lambda_i}-\frac{d}{\lambda_i^2}.
\]
Noting that
\[
1-c_nm_1'(d)\lambda_i=\frac{d}{\lambda_i}\frac{1}{d'},
\]
we have
\begin{equation}\label{thm3:eq13}
\mathrm{T}_n(\lambda_i)=\frac{\lambda_i\psi_n'(\lambda_i)}{\psi_n(\lambda_i)}\left[\bLambda_{1}^{-\frac{1}{2}}\mathbf R_n(\psi_n(\lambda_i))\bLambda_{1}^{-\frac{1}{2}}\right]_{ii}.
\end{equation}

Recall $\Z_1=\Y\bGamma_1\bLambda_1^{\frac{1}{2}}$ and  $\Z_2=\Y\bGamma_2\bLambda_2^{\frac{1}{2}}$. Note that $\mathbf R_n(d)=\{\Z_1^\top(\I+\H_n(d))\Z_1-\tr(\I+\H_n(d))\bLambda_{1}\}/\sqrt{n}$, then we have
\begin{equation}\label{thm3:eq15}
\bLambda_{1}^{-\frac{1}{2}}\mathbf R_n(\psi_n(\lambda_i))\bLambda_{1}^{-\frac{1}{2}}=\{\bGamma_1^\top\Y^\top(\I+\H_n(\psi_n(\lambda_i)))\Y\bGamma_1-\tr(\I+\H_n(\psi_n(\lambda_i)))\I\}/\sqrt{n}.
\end{equation}
Next, noting that
\[
\I+\H_n(d)=\I+\frac{1}{n}\Z_2(d \I-\S_{22})^{-1}\Z_2^\top=d(d \I-\frac{1}{n}\Z_2\Z_2^\top)^{-1},
\]
we have
\begin{equation}\label{thm3:eq16}
\bGamma_1^\top\Y^\top(\I+\H_n(\psi_n(\lambda_i)))\Y\bGamma_1=\psi_n(\lambda_i)\bGamma_1^\top\Y^\top(\psi_n(\lambda_i)\I-\frac{1}{n}\Y\bGamma_2\bLambda_2\bGamma_2^\top\Y^\top)^{-1}\Y\bGamma_1.
\end{equation}
 By \eqref{thm3:eq7}, \eqref{thm3:eq13}, \eqref{thm3:eq15} and \eqref{thm3:eq16}, we have
 \begin{equation}\label{thm3:eq013}
\mathrm{T}_n(\lambda_i)=\frac{\lambda_i\psi_n'(\lambda_i)}{\psi_n(\lambda_i)}\mathrm{\tilde{T}}_n(\lambda_i)+o_p(1),
\end{equation}
where
\[
\mathrm{\tilde{T}}_n(\lambda_i)=\sqrt{n}\left[\frac{1}{n}\bGamma_1^\top\Y^\top(\I
-\frac{1}{\psi_n(\lambda_i)}\frac{1}{n}\Y\bGamma_2\bLambda_2\bGamma_2^\top\Y^\top)^{-1}\Y\bGamma_1-\frac{\psi_n(\lambda_i)}{\lambda_i}\I\right]_{ii}.
\]
By Theorem 2.1 of \citet{Zhang:Zheng:Pan:Zhong:2022} (it also holds under the conditions of this paper), we have
\begin{equation}\label{thm3:eq12}
\frac{\mathrm{\tilde{T}}_n(\lambda_i)}{\tilde{\sigma}_{\lambda_i}}\stackrel{d}{\longrightarrow}N(0,1),
\end{equation}
where $\tilde{\sigma}_{\lambda_i}^2=\frac{\psi_n^2(\lambda_i)\{(\nu_4-3)\psi_n'(\lambda_i)\sum_{j=1}^p\Gamma_{ji}^4+2\}}{\lambda_i^2\psi_n'(\lambda_i)}$.
By \eqref{thm3:eq14}, \eqref{thm3:eq013} and  \eqref{thm3:eq12}, the limiting distribution of $\eta_i$ is normal, i.e.,
\begin{equation}\label{thm3:eq20}
\frac{\eta_i}{\sigma_{\lambda_i}}\stackrel{d}{\longrightarrow}N(0,1),
\end{equation}
where $\sigma_{\lambda_i}^2=(\nu_4-3)[\psi_n'(\lambda_i)]^2\sum_{j=1}^p\Gamma_{ji}^4+2\psi_n'(\lambda_i)$.

By \eqref{thm1:eq101} and \eqref{thm3:eq10}, we have
\begin{equation}\label{thm3:eq21}
\begin{aligned}
&\underline{s}_n(d_i)+\frac{1}{\lambda_i}\\
=&\underline{s}_n(d_i)-\underline{s}(\psi_n(\lambda_i))\\
=&-\frac{1-c_n}{d_i}+c_ns_n(d_i)+\frac{1-c_n}{\psi_n(\lambda_i)}-c_ns(\psi_n(\lambda_i))\\
=&\frac{(1-c_n)[d_i-\psi_n(\lambda_i)]}{d_i\psi_n(\lambda_i)}+\frac{c_n}{d_i}\left[d_is_n(d_i)-d_is(\psi_n(\lambda_i))\right]\\
=&\frac{\lambda_i(1-c_n)}{d_i\psi_n(\lambda_i)}\cdot\frac{[d_i-\psi_n(\lambda_i)]}{\lambda_i}+\frac{c_n}{d_i}\left[d_is(d_i)-d_is(\psi_n(\lambda_i))\right]+\frac{c_n}{d_i}\left[d_is_n(d_i)-d_is(d_i)\right]\\
=&\frac{\lambda_i(1-c_n)}{d_i\psi_n(\lambda_i)}\cdot\frac{d_i-\psi_n(\lambda_i)]}{\lambda_i}+c_n\left[s(d_i)-s(\psi_n(\lambda_i))\right]+\frac{c_n}{d_i}\mathcal{O}_p(1/n)\\
=&\frac{\lambda_i(1-c_n)}{d_i\psi_n(\lambda_i)}\cdot\frac{[d_i-\psi_n(\lambda_i)]}{\lambda_i}+c_n\int\frac{d_i-\psi_n(\lambda_i)}{(x-d_i)(x-\psi_n(\lambda_i))}dF^{c_n,H_n}(x)+\frac{1}{d_i}\mathcal{O}_p(1/n)\\
=&\frac{d_i-\psi_n(\lambda_i)}{\lambda_i}\left\{\frac{\lambda_i(1-c_n)}{d_i\psi_n(\lambda_i)}+c_n\int\frac{\lambda_i}{(x-d_i)(x-\psi_n(\lambda_i))}dF^{c_n,H_n}(x)\right\}+\frac{1}{d_i}\mathcal{O}_p(1/n).
\end{aligned}
\end{equation}
By \eqref{thm3:eq21}, we have
\begin{equation}\label{thm3:eq22}
\begin{aligned}
&\frac{-\underline{s}_n^{-1}(d_i)-\lambda_i}{\lambda_i}\\
=&-\underline{s}_n^{-1}(d_i)[\underline{s}_n(d_i)+\frac{1}{\lambda_i}]\\
=&\frac{-\underline{s}_n^{-1}(d_i)[d_i-\psi_n(\lambda_i)]}{\lambda_i}\left\{\frac{\lambda_i(1-c_n)}{d_i\psi_n(\lambda_i)}+c_n\int\frac{\lambda_i}{(x-d_i)(x-\psi_n(\lambda_i))}dF^{c_n,H_n}(x)\right\}\\
&+\frac{-\underline{s}_n^{-1}(d_i)}{d_i}\mathcal{O}_p(1/n)\\
=&\frac{d_i-\psi_n(\lambda_i)}{\lambda_i}\left\{\frac{-\underline{s}_n^{-1}(d_i)\lambda_i(1-c_n)}{d_i\psi_n(\lambda_i)}+c_n\int\frac{-\underline{s}_n^{-1}(d_i)\lambda_i}{(x-d_i)(x-\psi_n(\lambda_i))}dF^{c_n,H_n}(x)\right\}+\mathcal{O}_p(1/n).
\end{aligned}
\end{equation}
By \eqref{eq:thm3} and \eqref{thm2:eq12}, we have
\begin{equation}\label{thm3:eq23}
\begin{aligned}
&\frac{-\underline{s}_n^{-1}(d_i)\lambda_i(1-c_n)}{d_i\psi_n(\lambda_i)}+c_n\int\frac{-\underline{s}_n^{-1}(d_i)\lambda_i}{(x-d_i)(x-\psi_n(\lambda_i))}dF^{c_n,H_n}(x)\\
=&\lambda_i^2\left\{\frac{(1-c_n)}{\psi_n^2(\lambda_i)}+c_n\int\frac{1}{(x-\psi_n(\lambda_i))^2}dF^{c_n,H_n}(x)\right\}+\mathcal{O}_p(k/\sqrt{n})\\
=&\lambda_i^2\underline{s}'(\psi_n(\lambda_i))+\mathcal{O}_p(k/\sqrt{n}).
\end{aligned}
\end{equation}
By \eqref{thm3:eq10}, we have
\begin{equation}\label{thm3:eq24}
\underline{s}(d)=-\frac{1}{\lambda_i},
\end{equation}
where $d=\psi_n(\lambda_i)$.
By \eqref{thm3:eq24}, we have
\begin{equation}\label{thm3:eq25}
\underline{s}'(d)=\frac{1}{\lambda_i^2d'}.
\end{equation}
Combining
\eqref{thm3:eq22}, \eqref{thm3:eq23} and \eqref{thm3:eq25}, we have
\begin{equation}\label{thm3:eq26}
\begin{aligned}
&\frac{-\underline{s}_n^{-1}(d_i)-\lambda_i}{\lambda_i}\\
=&\frac{d_i-\psi_n(\lambda_i)}{\lambda_i}\left\{\frac{1}{\psi_n'(\lambda_i)}+\mathcal{O}_p(k/\sqrt{n})\right\}+\mathcal{O}_p(1/n).
\end{aligned}
\end{equation}
Since $d_i\underline{s}_n(d_i)=\mathcal{O}_p(1)$, by \eqref{eq:thm01} and \eqref{thm1:eq202}, we have
\begin{equation}\label{thm3:eq27}
\begin{aligned}
&\frac{-\hat{\underline{s}}_n^{-1}(d_i)}{\lambda_i}-\frac{-\underline{s}_n^{-1}(d_i)}{\lambda_i}\\
=&\frac{d_i\hat{\underline{s}}_n(d_i)-d_i\underline{s}_n(d_i)}{\lambda_i\hat{\underline{s}}_n(d_i)d_i\underline{s}_n(d_i)}=\mathcal{O}_p(k/n).
\end{aligned}
\end{equation}
By \eqref{thm3:eq26} and \eqref{thm3:eq27}, we have
\begin{equation}\label{thm3:eq28}
\frac{\sqrt{n}(-\hat{\underline{s}}_n^{-1}(d_i)-\lambda_i)}{\lambda_i}=\frac{\sqrt{n}(d_i-\psi_n(\lambda_i))}{\lambda_i\psi_n'(\lambda_i)}+o_p(1).
\end{equation}
By \eqref{thm3:eq20} and \eqref{thm3:eq28}, we have
\[
\frac{\sqrt{n}(-\hat{\underline{s}}_n^{-1}(d_i)-\lambda_i)}{\lambda_i\tau_{\lambda_i}}\stackrel{d}{\longrightarrow}N(0,1),
\]
where $\tau_{\lambda_i}^2=(\nu_4-3)\sum_{j=1}^p\Gamma_{ji}^4+2/\psi_n'(\lambda_i)$.
\eop

\subsection{Proof of Theorem \ref{eq:thm333}}
Similar to the proof of Theorem \ref{thm3}, \eqref{thm3:eq9} can be reorganized as
$$
\bLambda_{1}^{-\frac{1}{2}}(d_i \I-\mathbf K_n(d_i))\bLambda_{1}^{-\frac{1}{2}}=\B_1+\B_2+\B_3+\B_4,
$$
where
\begin{eqnarray*}
\B_1&=&\Diag\left\{\psi_n(\lambda_i)(\frac{1}{\lambda_1}-\frac{1}{\lambda_i}),\ldots,\psi_n(\lambda_i)(\frac{1}{\lambda_k}-\frac{1}{\lambda_i})\right\}\\
&&+\frac{\eta_i\lambda_i}{\sqrt{n}}\Diag\left\{\frac{1}{\lambda_1},\ldots,\frac{1}{\lambda_k}\right\},\\
\B_2&=&\mathcal{O}_p(1/n)\I+\frac{\eta_i\lambda_i}{\sqrt{n}}c_nm_3(\psi_n(\lambda_i))(1+\mathcal{O}_p(k/\sqrt{n}))\I,\\
\B_3&=&-\frac{1}{\sqrt{n}}\bLambda_{1}^{-\frac{1}{2}}\mathbf R_n(\psi_n(\lambda_i))\bLambda_{1}^{-\frac{1}{2}},\\
\B_4&=&\frac{\eta_i}{\sqrt{n}\lambda_i}\W(d_i)\Diag\left\{w_1(d_i),\ldots,w_k(d_i)\right\}\W(d_i)^\top.
\end{eqnarray*}
It is seen that the $S(i)\times S(i)$ diagonal block of $\B_1$ is $\frac{\eta_i}{\sqrt{n}}\I_{|S(i)|}$. The $S(i)\times S(i)$ diagonal block of $\B_2$ is
\[
\mathcal{O}_p(1/n)\I_{|S(i)|}+\frac{\eta_i}{\sqrt{n}}\lambda_ic_nm_3(\psi_n(\lambda_i))(1+\mathcal{O}_p(k/\sqrt{n}))\I_{|S(i)|}.
\]
As $\W(d_i)$ is an orthogonal matrix and $\max_{1\leq j\leq k}\vert w_j(d_i)\vert=\mathcal{O}_p(k/\sqrt{n})$, the $S(i)\times S(i)$ diagonal block  of $\B_4$ is $\frac{\eta_i}{\lambda_i}\mathcal{O}_p(k/n)\I_{|S(i)|}$.
Thus, the determinant of the $S(i)\times S(i)$ diagonal block of $\bLambda_{1}^{-\frac{1}{2}}(d_i \I-\mathbf K_n(d_i))\bLambda_{1}^{-\frac{1}{2}}$ can be written as
\begin{eqnarray*}
&&|\frac{\eta_i}{\sqrt{n}}\{1+c_nm_3(\psi_n(\lambda_i))\lambda_i(1+\mathcal{O}_p(k/\sqrt{n}))\}\I_{|S(i)|}-\left\{\frac{1}{\sqrt{n}}\bLambda_{1}^{-\frac{1}{2}}\mathbf R_n(\psi_n(\lambda_i))\bLambda_{1}^{-\frac{1}{2}}\right\}_{S(i)\times S(i)}\\
&&+\frac{\eta_i}{\lambda_i}\mathcal{O}_p(1/n)\I_{|S(i)|}+\mathcal{O}_p(k/n)\I_{|S(i)|}|=\mathcal{O}_p((\frac{k^{3}}{n})^{|S(i)|}),
\end{eqnarray*}
where the equality is true by Lemma \ref{lemma:diagonal11}, as $\A_n(d_i)=\bLambda_{1}^{-\frac{1}{2}}(d_i \I-\mathbf K_n(d_i))\bLambda_{1}^{-\frac{1}{2}}$.
By Theorem \ref{thm222}, for  $1\leq i \leq k$, we have
\[
\eta_i=\frac{\sqrt{n}(d_i-\psi_n(\lambda_i))}{\lambda_i}=\mathcal{O}_p(k).
\]
Thus, under the condition $k=o(n^{1/6})$, $\eta_i$ tends to a solution of
\begin{equation*}
|\eta_i(1+c_nm_3(\psi_n(\lambda_i))\lambda_i)\I_{|S(i)|}-[\bLambda_{1}^{-\frac{1}{2}}\mathbf R_n(\psi_n(\lambda_i))\bLambda_{1}^{-\frac{1}{2}}]_{S(i)\times S(i)}|=0.
\end{equation*}
Let $\mathbf R(\psi(\lambda_i))$ be the Gaussian matrix limit of the random matrix $\mathbf R_n(\psi_n(\lambda_i))$.
Similar to the proof of \eqref{thm3:eq13}, it is seen that $\eta_i$ tends to an eigenvalue of the matrix $\mathrm{T}(\lambda_i)$, which is defined as
\[
\mathrm{T}(\lambda_i)\triangleq\frac{\lambda_i\psi'(\lambda_i)}{\psi(\lambda_i)}\left[\bLambda_{1}^{-\frac{1}{2}}\mathbf R(\psi(\lambda_i))\bLambda_{1}^{-\frac{1}{2}}\right]_{S(i)\times S(i)}.
\]
Therefore, as the index $i$ is arbitrary, all the random variables $\eta_{i}'s, i\in S(i)$, converge in distribution to the set of eigenvalues of the above matrix.
Finally, similar to the proof of Theorem \ref{thm3}, for $1\le i\le k$, the joint distribution of
\begin{equation*}
\frac{\sqrt{n}(-\hat{\underline{s}}_n^{-1}(d_i)-\lambda_i)}{\lambda_i},\,\,\,\, i\in S(i),
\end{equation*}
converges to the joint distribution of the $|S(i)|$ eigenvalues of the Gaussian random matrix
\[
\mathrm{\tilde{T}}(\lambda_i)\triangleq\frac{\lambda_i}{\psi(\lambda_i)}\left[\bLambda_{1}^{-\frac{1}{2}}\mathbf R(\psi(\lambda_i))\bLambda_{1}^{-\frac{1}{2}}\right]_{S(i)\times S(i)}.
\]
\eop

\subsection{Proof of Theorem \ref{thm4}}
By the definition of $\mathrm{GIC}(k)$, we have
\[
\mathrm{GIC}(k')-\mathrm{GIC}(k)=\log \tilde{L}_k-\log \tilde{L}_{k'}-(k-k')\{p-(k+k')/2+1/2\}C_n.
\]
Also, recall that
$$
\log \tilde{L}_k=-\frac{n}{2}\sum_{i=1}^k\log d_i-\frac{n}{2}\sum_{i=k+1}^p(d_i-1).
$$

First, consider the case of $k'<k$. Let $h(x)=x-1-\log x$ and note that $h'(x)>0$ for $x>1$. As $\lambda_k>1$, there exists a constant $\epsilon_k>0$, such that $h(\lambda_k)\ge\epsilon_k$, which also implies $\lambda_k\geq1+\epsilon_k$.
Define the event $A=\{d_k>1+\epsilon_k/2\}$. Conditioning on event $A$, we have
\begin{eqnarray*}
&&\mathrm{GIC}(k')-\mathrm{GIC}(k)\\
&=&\log \tilde{L}_k-\log \tilde{L}_{k'}-(k-k')\{p-(k+k')/2+1/2\}C_n\\
&=&\frac{n}{2}\sum_{i=k'+1}^k(d_i-1-\log d_i)-(k-k')\{p-(k+k')/2+1/2\}C_n\\
&\geq&\frac{n}{2}(k-k')\left[d_k-1-\log d_k-2\{p-(k+k')/2+1/2\}C_n/n\right],
\end{eqnarray*}
where we have used the fact that $h'(x)>0$ for $x>1$ and $d_i\ge d_k$ for $i=k'+1,\ldots,k$.
We can further write
\begin{eqnarray*}
&&\mathrm{GIC}(k')-\mathrm{GIC}(k)\\
&\geq&\frac{n}{2}(k-k')\left[d_k-1-\log d_k-2\{p-(k+k')/2+1/2\}C_n/n\right]\\
&=&\frac{n}{2}(k-k')\left[h(\lambda_k)+d_k-\lambda_k-(\log d_k-\log\lambda_k)-2\{p-(k+k')/2+1/2\}C_n/n\right]\\
&=&\frac{n}{2}(k-k')\left[h(\lambda_k)+d_k-\lambda_k-(d_k-\lambda_k)/\tilde{d}_k-2\{p-(k+k')/2+1/2\}C_n/n\right]\\
&=&\frac{n}{2}(k-k')\left[h(\lambda_k)+(d_k-\lambda_k)(1-1/\tilde{d}_k)-2\{p-(k+k')/2+1/2\}C_n/n\right],
\end{eqnarray*}
where $\tilde{d}_k$ lies between $d_k$ and $\lambda_k$.
By Lemma \ref{lemma:weyl}, we have, for $1\le i\le p$,
\[
\vert\lambda_i-d_i\vert\leq ||\S_n-\bSigma||.
\]
Moreover, conditioning on $A$, we have
\[
0<\frac{\epsilon_k}{2+\epsilon_k}<1-1/\tilde{d}_k<1.
\]
Recall $h(\lambda_k)\ge\epsilon_k$, and we have
\begin{eqnarray*}
&&\mathrm{GIC}(k')-\mathrm{GIC}(k)\\
&\ge&\frac{n}{2}(k-k')\left[h(\lambda_k)+(d_k-\lambda_k)(1-1/\tilde{d}_k)-2\{p-(k+k')/2+1/2\}C_n/n\right]\\
&\ge&\frac{n}{2}(k-k')\left[\epsilon_k-||\S_n-\bSigma||-2\{p-(k+k')/2+1/2\}C_n/n\right]\\
&\ge&\frac{n}{2}(k-k')\left(\epsilon_k-||\S_n-\bSigma||-2pC_n/n\right).
\end{eqnarray*}
By Lemma \ref{lemma:spectral}, it holds that
\begin{eqnarray*}
\mathbb{P}(A)&\ge&\mathbb{P}(|d_k-\lambda_k|\le\epsilon_k/2)\\
&\ge&\mathbb{P}(||\S_n-\bSigma||\le\epsilon_k/2)\geq1-4c_0(2p+1)p^3/(n\epsilon_k^2),
\end{eqnarray*}
where $c_0$ is some positive constant.
Therefore, it holds that
$$
\mathrm{GIC}(k')-\mathrm{GIC}(k)\ge\frac{n}{2}(k-k')\left(\epsilon_k/2-2pC_n/n\right),
$$
with probability as least $1-c_1/n$, where $c_1$ is some positive constant.

Next, consider the case of $k'>k$. Define the event $A=\{\min_{k+1\leq i\leq k'}d_i>1/2\}$. By Taylor's expansion, we have
\begin{eqnarray*}\log \tilde{L}_k-\log \tilde{L}_{k'}&=&\frac{n}{2}\sum_{i=k+1}^{k'}(\log d_i+1-d_i)\\
&=&\frac{n}{2}\sum_{i=k+1}^{k'}\left\{(d_i-1)-\frac{1}{2}(d_i-1)^2/\tilde{d}_i^2+1-d_i\right\}\\
&=&-\frac{n}{4}\sum_{i=k+1}^{k'}(d_i-1)^2/\tilde{d}_i^2,
\end{eqnarray*}
where $\tilde{d}_i$ lies between $d_i$ and $1$. Conditioning on $A$, we have $\max_{k+1\leq i\leq k'}1/\tilde{d}_i<2$, and we also have
\begin{eqnarray*}
&&\mathrm{GIC}(k')-\mathrm{GIC}(k)\\
&=&\log \tilde{L}_k-\log \tilde{L}_{k'}-(k-k')\{p-(k+k')/2+1/2\}C_n\\
&=&-\frac{n}{4}\sum_{i=k+1}^{k'}(d_i-1)^2/\tilde{d}_i^2-(k-k')\{p-(k+k')/2+1/2\}C_n\\
&=&\frac{n(k'-k)}{4}\left[4\{p-(k+k')/2+1/2\}C_n/n-\frac{1}{k'-k}\sum_{i=k+1}^{k'}(d_i-1)^2/\tilde{d}_i^2\right]\\
&\ge&\frac{n(k'-k)}{4}\left[4\{p-(k+k')/2+1/2\}C_n/n-\max_{k+1\leq i\leq k'}(d_i-1)^2/\tilde{d}_i^2\right]\\
&\ge&n(k'-k)\left[\{p-(k+k')/2+1/2\}C_n/n-\max_{k+1\leq i\leq k'}(d_i-1)^2\right],
\end{eqnarray*}
where we have used the fact that $\max_{k+1\leq i\leq k'}1/\tilde{d}_i<2$. As $\lambda_i=1$ for $k+1\leq i\leq k'$, by Lemma \ref{lemma:weyl}, we have
$$
\mathrm{GIC}(k')-\mathrm{GIC}(k)\ge n(k'-k)\left[\{p-(k+k')/2+1/2\}C_n/n-||\S_n-\bSigma||^2\right].
$$
By Lemma \ref{lemma:spectral}, it holds that
\begin{eqnarray*}
\mathbb{P}\left[||\S_n-\bSigma||^2\le\{p-(k+k')/2+1/2\}C_n/(2n)\right]\geq1-\frac{2c_0(2p+1)p^3}{\{p-(k+k')/2+1/2\}C_n},
\end{eqnarray*}
and that
\begin{eqnarray*}
\mathbb{P}(A)&\ge&{1-\mathbb{P}(\cup_{i=k+1}^{k'}|d_i-1|\geq1/2)}\\
&\ge&{1-(k'-k)\mathbb{P}(||\S_n-\bSigma||\geq1/2)\geq1-\frac{2c_0(k'-k)(2p+1)p^3}{\{p-(k+k')/2+1/2\}C_n}.}
\end{eqnarray*}
Therefore, it holds that
$$
\mathrm{GIC}(k')-\mathrm{GIC}(k)\ge{\frac{(k'-k)C_n}{4}},
$$
with probability at least $1-c_2/C_n$, where $c_2$ is some positive constant.

\eop

\subsection{Proof of Theorem \ref{thm5}}
\label{sec:proof5}
In what follows, we only give the proof for the case of $\mathbb{E}Y_{ij}^4<\infty$.  That is, by Theorem \ref{thm1} for $\mathbb{E}Y_{ij}^4<\infty$,  we can show that Theorem \ref{thm5} holds for $k=o(n^{1/4})$. For the case of $Y_{ij}$ is Gaussian, instead of using Theorem \ref{thm1} for $\mathbb{E}Y_{ij}^4<\infty$, by Theorem \ref{thm1} for $Y_{ij}$ Gaussian, we can show that Theorem \ref{thm5} holds for $k=o(n^{1/3})$. The proof is very similar to the case of $\mathbb{E}Y_{ij}^4<\infty$, and the details are thus omitted.

Recall that $\bar{d}_{k+1}=\frac{1}{p-k}\sum_{i=k+1}^pd_i$ and $\bar{d}_{k'+1}=\frac{1}{p-k'}\sum_{i=k'+1}^pd_i$. We have
\begin{equation}\label{eq:gicdiff}
\mathrm{GIC}(k')-\mathrm{GIC}(k)=\log L_k-\log L_{k'}-\gamma(k-k')\{p-(k+k')/2+1/2\},
\end{equation}
where
\begin{eqnarray*}
\log L_{k} &=& -\frac{n}{2}\sum_{i=1}^{k}\log d_i-\frac{n}{2}(p-k)\log\bar{d}_{k+1}, \\
\log L_{k'} &=& -\frac{n}{2}\sum_{i=1}^{k'}\log d_i-\frac{n}{2}(p-k')\log\bar{d}_{k'+1}.
\end{eqnarray*}
Suppose $k'<k$. By Taylor's expansion, we have
\begin{eqnarray*}
&&\log L_k-\log L_{k'}\\
&=&-\frac{n}{2}\sum_{i=1}^{k}\log d_i-\frac{n}{2}(p-k)\log\bar{d}_{k+1}+\frac{n}{2}\sum_{i=1}^{k'}\log d_i+\frac{n}{2}(p-k')\log\bar{d}_{k'+1}\\
&=&-\frac{n}{2}(p-k')\log\bar{d}_{k+1}+\frac{n}{2}(p-k')\log\bar{d}_{k'+1}+\frac{n}{2}(k-k')\log\bar{d}_{k+1}-\frac{n}{2}\sum_{i=k'+1}^k\log d_i\\
&=&\frac{n}{2}(p-k')\log(\bar{d}_{k'+1}/\bar{d}_{k+1})+\frac{n}{2}(k-k')\log\bar{d}_{k+1}-\frac{n}{2}\sum_{i=k'+1}^k\log d_i\\
&=&\frac{n}{2}(p-k')\log\left\{\left(1-\frac{k-k'}{p-k'}\right)\frac{\sum_{i=k'+1}^pd_i}{\sum_{i=k+1}^pd_i}\right\}+\frac{n}{2}(k-k')\log\bar{d}_{k+1}-\frac{n}{2}\sum_{i=k'+1}^k\log {d_i}\\
&=&\frac{n}{2}(p-k')\left[-\frac{k-k'}{p-k'}+\mathcal{O}\left\{\frac{(k-k')^2}{(p-k')^2}\right\}+\log\frac{\sum_{i=k'+1}^pd_i}{\sum_{i=k+1}^pd_i}\right]+\frac{n}{2}(k-k')\log\bar{d}_{k+1}-\frac{n}{2}\sum_{i=k'+1}^k\log d_i\\
&=&-\frac{n}{2}(k-k')+\mathcal{O}(k-k')+\frac{n}{2}(p-k')\log\left(\frac{\frac{1}{p-k}\sum_{i=k'+1}^kd_i+\frac{1}{p-k}\sum_{i=k+1}^pd_i}{\frac{1}{p-k}\sum_{i=k+1}^pd_i}\right)\\
&&+\frac{n}{2}(k-k')\log\bar{d}_{k+1}-\frac{n}{2}\sum_{i=k'+1}^k\log d_i\\
&=&-\frac{n}{2}(k-k')+\mathcal{O}(k-k')+\frac{n}{2}(p-k')\log\left(1+\frac{\frac{1}{p-k}\sum_{i=k'+1}^kd_i}{\bar{d}_{k+1}}\right)\\
&&+\frac{n}{2}(k-k')\log\bar{d}_{k+1}-\frac{n}{2}\sum_{i=k'+1}^k\log d_i,
\end{eqnarray*}
where we have used the fact that ${\bar{d}_{k+1}=\frac{1}{p-k}\sum_{i=k+1}^pd_i}$ in the last step.

Next, we consider three different cases for $\frac{1}{p-k}\sum_{i=k'+1}^k\psi_n(\lambda_i)$. By Theorem \ref{thm1}, for all $1\leq i\leq k$, we have
\[
\frac{d_{i}}{\psi_n(\lambda_{i})}=1+o_p(1).
\]
Under the case that $\frac{1}{p-k}\sum_{i=k'+1}^k\psi_n(\lambda_i)=o(1)$, it holds that
$\frac{1}{p-k}\sum_{i=k'+1}^kd_i=o_p(1)$. Using Taylor's expansion, we have
\begin{eqnarray*}
\log\left(1+\frac{\frac{1}{p-k}\sum_{i=k'+1}^kd_i}{\bar{d}_{k+1}}\right)&=&\frac{1}{1+x_d}\times\frac{\frac{1}{p-k}\sum_{i=k'+1}^kd_i}{\bar{d}_{k+1}}\\
&=&\frac{1}{(1+x_d)\bar{d}_{k+1}}\times\frac{1}{p-k}\sum_{i=k'+1}^kd_i,
\end{eqnarray*}
where $x_d\in(0,\frac{\frac{1}{p-k}\sum_{i=k'+1}^kd_i}{\bar{d}_{k+1}})$.
According to \cite{Bai:Choi:Fujikoshi:2018}, it holds that
\[
\bar{d}_{k+1}=\frac{1}{p-k}\sum_{i=k+1}^pd_i\stackrel{a.s.}{\longrightarrow}1.
\]
Since $(1+x_d)\bar{d}_{k+1}\in (\bar{d}_{k+1}, \bar{d}_{k+1}+\frac{1}{p-k}\sum_{i=k'+1}^kd_i)$,
it holds that
\[
(1+x_d)\bar{d}_{k+1}=1+o_p(1).
\]
Hence,
\[
\log\left(1+\frac{\frac{1}{p-k}\sum_{i=k'+1}^kd_i}{\bar{d}_{k+1}}\right)=\frac{1}{p-k}\sum_{i=k'+1}^kd_i(1+o_p(1)).
\]
Thus, we have
\begin{eqnarray*}
&&\log L_k-\log L_{k'}\\
&=&-\frac{n}{2}(k-k')(1+o_p(1))+\frac{n}{2}\times\frac{p-k'}{p-k}\sum_{i=k'+1}^kd_i(1+o_p(1))-\frac{n}{2}\sum_{i=k'+1}^k\log d_i\\
&=&-\frac{n}{2}(k-k')(1+o_p(1))+\frac{n}{2}\left(1+\frac{k-k'}{p-k}\right)\sum_{i=k'+1}^kd_i(1+o_p(1))-\frac{n}{2}\sum_{i=k'+1}^k\log d_i\\
&=&\frac{n}{2}\sum_{i=k'+1}^k\{d_{i}(1+o_p(1))-1-\log d_{i}\}+\frac{n}{2}\times\frac{k-k'}{p-k}\sum_{i=k'+1}^kd_i(1+o_p(1)).
\end{eqnarray*}

Therefore, for all $k'<k$, by \eqref{eq:gicdiff}, we have
\begin{eqnarray*}
&&\mathrm{GIC}(k')-\mathrm{GIC}(k)\\
&=&\frac{n}{2}\sum_{i=k'+1}^k\{d_{i}(1+o_p(1))-1-\log d_{i}\}+\frac{n}{2}\times\frac{k-k'}{p-k}\sum_{i=k'+1}^kd_i(1+o_p(1))\\
&&\quad-\gamma(k-k')\{p-(k+k')/2+1/2\}\\
&\geq&\frac{n}{2}\sum_{i=k'+1}^k\{d_{i}(1+o_p(1))-1-\log d_{i}\}-\gamma(k-k')\{p-(k+k')/2+1/2\}\\
&=&\frac{n}{2}\sum_{i=k'+1}^k\{d_{i}(1+o_p(1))-1-\log d_{i}\}-\gamma(k-k')\{p-(k+k')/2+1/2\}\\
&\geq&\frac{n}{2}(k-k')\left\{d_{k}(1+o_p(1))-1-\log d_{k}-2\gamma \{p-(k+k')/2+1/2\}/n\right\}\\
&=&\frac{n}{2}(k-k')\left[\psi_n(\lambda_k)(1+o_p(1))-1-\log \psi_n(\lambda_k)(1+o_p(1))-2\gamma \{p-(k+k')/2+1/2\}/n\right]\\
&=&\frac{n}{2}(k-k')\left[\psi_n(\lambda_k)(1+o_p(1))-1-\log \psi_n(\lambda_k)-2\gamma \{p-(k+k')/2+1/2\}/n\right].
\end{eqnarray*}
Since $\gamma<\varphi(\lambda_{k})=\frac{1}{2c}[\psi(\lambda_k)-1-\log \psi(\lambda_k)]$, it then holds under this case that
$$
\mathbb{P}\{\mathrm{GIC}(k')>\mathrm{GIC}(k)\}\rightarrow 1.
$$

Under the case that $\frac{1}{p-k}\sum_{i=k'+1}^k\psi_n(\lambda_i)/C \rightarrow 1$, where $C$ is a constant, it holds that
$\frac{1}{p-k}\sum_{i=k'+1}^kd_i/C\stackrel{p}{\longrightarrow} 1$.
Hence, we may write
$\log\left(1+\frac{\frac{1}{p-k}\sum_{i=k'+1}^kd_i}{\bar{d}_{k+1}}\right)=\log(C+1)+o_p(1)$.
Moreover, it holds that
\[
\prod_{i=k'+1}^kd_i\leq \left(\frac{1}{k-k'}\sum_{i=k'+1}^kd_i\right)^{k-k'},
\]
\[
\sum_{i=k'+1}^k\log d_i\leq(k-k')\log\left(\frac{1}{k-k'}\sum_{i=k'+1}^kd_i\right)=(k-k')\{\log C+o_p(1)+(1+\mathcal{O}(1))\log p\}.
\]
Therefore, we have
\begin{eqnarray*}
&&\mathrm{GIC}(k')-\mathrm{GIC}(k)\\
&\geq&-\frac{n}{2}(k-k')(1+o_p(1))+\frac{n}{2}(p-k')\{\log(C+1)+o_p(1)\}\\
&&-\frac{n}{2}(k-k')\{\log C+o_p(1)+(1+\mathcal{O}(1))\log p\}-\gamma(k-k')\{p-(k+k')/2+1/2\}.
\end{eqnarray*}
It then holds under this case that
$$
\mathbb{P}\{\mathrm{GIC}(k')>\mathrm{GIC}(k)\}\rightarrow 1.
$$

Finally, for the case $\frac{1}{p-k}\sum_{i=k'+1}^k\psi_n(\lambda_i)/C_k\rightarrow1$, where $C_k\rightarrow\infty$, it holds that
$$
\frac{1}{p-k}\sum_{i=k'+1}^kd_i/C_k\stackrel{p}{\longrightarrow}1.
$$
Hence, we may write
$$
\log\left(1+\frac{\frac{1}{p-k}\sum_{i=k'+1}^kd_i}{\bar{d}_{k+1}}\right)=\log C_k+o_p(1).
$$
Moreover, it holds that
\[
\sum_{i=k'+1}^k\log d_i\leq(k-k')\log\left(\frac{1}{k-k'}\sum_{i=k'+1}^kd_i\right)=(k-k')\{\log C_k+o_p(1)+(1+\mathcal{O}(1))\log p\}.
\]
Therefore, under the case that $\frac{1}{p-k}\sum_{i=k'+1}^k\psi_n(\lambda_i)/C_k\rightarrow1$, we have
\begin{eqnarray*}
&&\mathrm{GIC}(k')-\mathrm{GIC}(k)\\
&\geq&-\frac{n}{2}(k-k')(1+o_p(1))+\frac{n}{2}(p-k')\{\log C_k+o_p(1)\}\\
&&-\frac{n}{2}(k-k')\{\log C_k+o_p(1)+(1+\mathcal{O}(1))\log p\}-\gamma(k-k')\{p-(k+k')/2+1/2\}.
\end{eqnarray*}
It then holds under this case that
$$
\mathbb{P}\{\mathrm{GIC}(k')>\mathrm{GIC}(k)\}\rightarrow 1.
$$

Suppose $k'>k$.
By Taylor's expansion, we have
\begin{eqnarray*}
&&\log L_k-\log L_{k'}\\
&=&\frac{n}{2}(p-k')\log(\bar{d}_{k'+1}/\bar{d}_{k+1})-\frac{n}{2}(k'-k)\log\bar{d}_{k+1}+\frac{n}{2}\sum_{i=k+1}^{k'}\log d_i\\
&=&\frac{n}{2}(p-k')\log\left\{\left(1+\frac{k'-k}{p-k'}\right)\frac{\sum_{i=k'+1}^pd_i}{\sum_{i=k+1}^pd_i}\right\}-\frac{n}{2}(k'-k)\log\bar{d}_{k+1}+\frac{n}{2}\sum_{i=k+1}^{k'}\log d_i\\
&=&\frac{n}{2}(k'-k)+\mathcal{O}(k'-k)+\frac{n}{2}(p-k')\log\left(\frac{\frac{1}{p-k}\sum_{i=k+1}^pd_i-\frac{1}{p-k}\sum_{i=k+1}^{k'}d_i}{\frac{1}{p-k}\sum_{i=k+1}^pd_i}\right)\\
&&\quad-\frac{n}{2}(k'-k)\log\bar{d}_{k+1}+\frac{n}{2}\sum_{i=k+1}^{k'}\log d_i\\
&=&\frac{n}{2}(k'-k)+\mathcal{O}(k'-k)+\frac{n}{2}(p-k')\log\left(1-\frac{\frac{1}{p-k}\sum_{i=k+1}^{k'}d_i}{\frac{1}{p-k}\sum_{i=k+1}^pd_i}\right)\\
&&-\frac{n}{2}(k'-k)\log\bar{d}_{k+1}+\frac{n}{2}\sum_{i=k+1}^{k'}\log d_i\\
&=&\frac{n}{2}(k'-k)(1+o_p(1))-\frac{n}{2}\times\frac{p-k'}{p-k}\sum_{i=k+1}^{k'}d_i(1+o_p(1))+\frac{n}{2}\sum_{i=k+1}^{k'}\log d_i\\
&=&\frac{n}{2}(k'-k)(1+o_p(1))-\frac{n}{2}\left(1-\frac{k'-k}{p-k}\right)\sum_{i=k+1}^{k'}d_i(1+o_p(1))+\frac{n}{2}\sum_{i=k+1}^{k'}\log d_i\\
&=&-\frac{n}{2}\sum_{i=k+1}^{k'}\{d_i-1-\log d_{i}+o_p(1)\}+\frac{n}{2}\times\frac{k'-k}{p-k}\sum_{i=k+1}^{k'}d_i(1+o_p(1)).
\end{eqnarray*}
According to \cite{Bai:Choi:Fujikoshi:2018}, for all $k<i=o(p)$, it holds that
\[
d_{i}\stackrel{a.s.}{\longrightarrow} (1+\sqrt{c})^2.
\]
Hence, for all $k<k'=o(p)$, we have
\[
\frac{1}{p-k}\sum_{i=k+1}^{k'}d_i\geq \frac{k'-k}{p-k}d_{k'}=o_p(1).
\]

Therefore, for all $k<k'=o(p)$, by \eqref{eq:gicdiff}, we have
\begin{eqnarray*}
&&\mathrm{GIC}(k')-\mathrm{GIC}(k)\\
&=&-\frac{n}{2}\sum_{i=k+1}^{k'}\{d_i-1-\log d_{i}+o_p(1)\}+\frac{n}{2}\times\frac{k'-k}{p-k}\sum_{i=k+1}^{k'}d_i(1+o_p(1))\\
&&\quad+\gamma(k'-k)\{p-(k'+k)/2-1/2\}\\
&\geq&-\frac{n}{2}\sum_{i=k+1}^{k'}\{d_{i}-1-\log d_{i}+o_p(1)\}+\gamma(k'-k)\{p-(k'+k)/2-1/2\}\\
&\geq&-\frac{n}{2}(k'-k)\{d_{k+1}-1-\log d_{k+1}+o_p(1)\}+\gamma(k'-k)\{p-(k'+k)/2-1/2\}\\
&=&\frac{n}{2}(k'-k)\left[2\gamma\{p-(k'+k)/2-1/2\}/n-\{(1+\sqrt{c})^2-1-2\log(1+\sqrt{c})\}+o_p(1)\right].
\end{eqnarray*}
Since $\gamma>\varphi(1+\sqrt{c})=1/2+\sqrt{1/c}-\log(1+\sqrt{c})/c$, it then holds under this case that
$$
\mathbb{P}\{\mathrm{GIC}(k')>\mathrm{GIC}(k)\}\rightarrow 1.
$$
\eop

\subsection{Proof of Theorem \ref{thm6}}
In what follows, we only give the proof for the case of $\mathbb{E}Y_{ij}^4<\infty$.  That is, by Theorem \ref{thm1} for $\mathbb{E}Y_{ij}^4<\infty$,  we can show that Theorem \ref{thm6} holds for $k=o(n^{1/4})$. For the case of $Y_{ij}$ is Gaussian,  instead of using Theorem \ref{thm1} for $\mathbb{E}Y_{ij}^4<\infty$, by Theorem \ref{thm1} for $Y_{ij}$ Gaussian,  we can show that Theorem \ref{thm6} holds for $k=o(n^{1/3})$. The proof is very similar to the case of $\mathbb{E}Y_{ij}^4<\infty$, and the details are thus omitted.

Let $\tilde{\R}_n=\D^{-1/2}\S_n\D^{-1/2}$. Since $\x_i=\bSigma^{\frac{1}{2}}\y_i$, we have
 \begin{equation*}
 \tilde{\R}_n=\frac{1}{n}\D^{-1/2}\bSigma^{1/2}\Y^\top\Y\bSigma^{1/2}\D^{-1/2}=\begin{pmatrix}
 \tilde{\R}_{11} \,\,& \tilde{\R}_{12} \\
  \tilde{\R}_{21} \,\,&  \tilde{\R}_{22}
\end{pmatrix}.
 \end{equation*}
Let $\hat{\sigma}_{jj}=\frac{1}{n}\v_j^\top\Y^\top\Y\v_j$,
where $\v_j=\bSigma^{1/2}\D^{-1/2}\e_j=(v_j)_{p\times 1}$ and $\e_j$ is the $j$-th column of $\I$. It is known that $\v_j^\top\v_j=1$.

When $\mathbb{E}Y_{ij}^{4}<\infty$, by Markov's inequality,  it holds that
\begin{eqnarray*}
&&P\left\{\max_{1\leq j\leq p}\left|\hat{\sigma}_{jj}-1\right|>t\right\}\\
&=&P\left\{\max_{1\leq j\leq p}\left|\frac{1}{n}\sum_{s=1}^n\sum_{1\leq i,l\leq p}v_iv_lY_{is}Y_{ls}-1\right|>t\right\}\\
&=&P\left\{\max_{1\leq j\leq p}\left|\sum_{s=1}^n\left[\sum_{1\leq i\leq p}(Y_{is}^2-1)v_i^2+2\sum_{1\leq i<l\leq p}v_iv_lY_{is}Y_{ls}\right]\right|>nt\right\}\\
&\leq&\sum_{1\leq j\leq p}P\left\{\left|\sum_{s=1}^nW_s\right|>nt\right\}\\
&\leq&\sum_{1\leq j\leq p}\frac{\mathbb{E}(\sum_{s=1}^nW_s)^2}{(nt)^2}\\
&\leq&\sum_{1\leq j\leq p}\frac{\sum_{s=1}^n\left[\sum_{1\leq i\leq p}\mathbb{E}(Y_{is}^2-1)^2v_i^4+4\sum_{1\leq i<l\leq p}v_i^2v_l^2\right]}{(nt)^{2}}\\
&\leq&\sum_{1\leq j\leq p}\frac{\sum_{s=1}^n\left[\eta\sum_{1\leq i\leq p}v_i^4+4\sum_{1\leq i<l\leq p}v_i^2v_l^2\right]}{(nt)^{2}}\\
&\leq&\sum_{1\leq j\leq p}\frac{\sum_{s=1}^n\left[(\eta-2)\sum_{1\leq i\leq p}v_i^4+2\right]}{(nt)^{2}}\\
&\leq&\left\{
\begin{array}{ll} \frac{n(2p+\eta-2)}{(nt)^{2}}, & \,\mbox{if}\,\ \eta\in [0,2], \\
\frac{np\eta}{(nt)^{2}}, & \,\mbox{if}\,\ \eta>2,
\end{array}
\right.
\end{eqnarray*}
where
\[
W_s=\sum_{1\leq i\leq p}(Y_{is}^2-1)v_i^2+2\sum_{1\leq i<l\leq p}v_iv_lY_{is}Y_{ls}
\]
and $\eta=\mathbb{E}(Y_{ij}^2-1)^2=\mathbb{E}Y_{ij}^{4}-1$.
Therefore, it follows that
$\max_{1\leq j\leq p}\left|\hat{\sigma}_{jj}-1\right|=\mathcal{O}_p(1)$.
By (S.63) in \citet{Fan:Guo:Zheng:2022}, we have $\max_{1\leq j\leq p}|\frac{\lambda_j(\tilde{\R}_n)}{\lambda_j(\R_n)}-1|=\mathcal{O}_p(1)$.
It implies that under the condition $\mathbb{E}Y_{ij}^{4}<\infty$, the eigenvalues of $\tilde{\R}_n$ and $\R_n$ may not be close.

Under the assumption that $\mathbb{E}|Y_{ij}|^{4+\zeta}<\infty$, where $\zeta>0$, 
we can assume that $Y_{ij}<\sqrt{n}\eta_n$, where $\eta_n$ satisfies $n^{1/4}\eta_n\rightarrow\infty$ (see, for example, \citet{Bai:Silverstein:2004}). Similar to the discussion in \citet{Fan:Guo:Zheng:2022}, we have $\max_{1\leq j\leq p}\left|\hat{\sigma}_{jj}-1\right|=o_p(1)$
and  $\max_{1\leq j\leq p}|\frac{\lambda_j(\tilde{\R}_n)}{\lambda_j(\R_n)}-1|=o_p(1)$.
Thus, instead of considering the eigenvalues of $\R_n$, we only consider those of  $\tilde{\R}_n$.

Recall $\x_i=\bSigma^{\frac{1}{2}}\y_i$, then we have $\tilde{\x}_i=\D^{-1/2}\bSigma^{1/2}\y_i$. That is, in the following theoretical analysis,  $\R=\D^{-1/2}\bSigma\D^{-1/2}$ is the population covariance matrix of $\tilde{\x}_i$, and $\tilde{\R}_n$ is the sample covariance matrix of $\tilde{\X}=\X\D^{-1/2}$. For simplicity, we still use $\{\lambda_i\}_{1\leq i\leq p}$ and $\{d_i\}_{1\leq i\leq p}$ to denote the eigenvalues of $\R$ and $\tilde{\R}_n$, respectively.

We need some additional notation. For $i\geq k$, we use $\{\beta_j\}_{i+1\leq j\leq p}$ to denote the eigenvalues of $\tilde{\R}_{22}$, which is a $(p-i)\times(p-i)$ matrix.  Let $F_{n,i}(x)=\frac{1}{p-i}\sum_{j=i+1}^{p}I_{(-\infty,x]}(\beta_j)$ and $H_{n,i}(x)=\frac{1}{p-i}\sum_{j=i+1}^{p}I_{(-\infty,x]}(\lambda_j)$ be the empirical spectral distributions of $\beta_{i+1},\ldots, \beta_p$ and $\lambda_{i+1},\ldots, \lambda_p$, respectively. Let $\psi_{n,i}(x)=x+c_{n,i}\int\frac{xt}{x-t}dH_{n,i}(t)$.
Define
\[
c_{n,i}=\frac{p-i}{n},\,\, s(d)=\int\frac{1}{x-d}dF^{c_{n,i},H_{n,i}}(x), \,\, \underline{s}(d)=-\frac{1-c_{n,i}}{d}+c_{n,i}s(d),
\]
\[
s_{n,i}(d)=\frac{1}{p-i}\sum_{j=i+1}^p\frac{1}{\beta_j-d}=\int\frac{1}{x-d}dF_{n,i}(x), \,\, \underline{s}_{n,i}(d)=-\frac{1-c_{n,i}}{d}+c_{n,i}s_{n,i}(d),
\]
and
\[
\hat{s}_{n,i}(d)=\frac{1}{p-i}\sum_{j=i+1}^p\frac{1}{d_j-d}, \,\, \hat{\underline{s}}_{n,i}(d)=-\frac{1-c_{n,i}}{d}+c_{n,i}\hat{s}_{n,i}(d).
\]
Let $\tilde{\delta}_i=\frac{1}{n}\sum_{j=k+1}^p\frac{d_j}{d_i-d_j}$. It is seen that $\tilde{\delta}_k=\delta_k$.

Suppose $k'<k$. Note that
\[
1+\tilde{\delta}_i=1+\frac{1}{n}\sum_{j=k+1}^p\frac{d_j-d_i+d_i}{d_i-d_j}=1-c_{n,k}+\frac{1}{n}\sum_{j=k+1}^p\frac{d_i}{d_i-d_j},
\]
and
\[
1+\delta_i=1+\frac{1}{n}\sum_{j=i+1}^p\frac{d_j-d_i+d_i}{d_i-d_j}=1-c_{n,i}+\frac{1}{n}\sum_{j=i+1}^p\frac{d_i}{d_i-d_j}.
\]
It is seen that
\begin{equation}\label{thm12:eq6}
\delta_i-\tilde{\delta}_i=\frac{1}{n}\sum_{j=i+1}^k\frac{d_j}{d_i-d_j},
\end{equation}
\begin{equation}\label{thm12:eq7}
\frac{d_i}{1+\tilde{\delta}_i}=-\hat{\underline{s}}_{n,k}^{-1}(d_i),
\end{equation}
where $\hat{\underline{s}}_{n,k}(d_i)=-\frac{1-c_{n,k}}{d_i}+c_{n,k}\hat{s}_{n,k}(d_i)$. For all $1\leq i\leq k$, similar to the proof of Theorem \ref{thm1}, we have
\begin{equation}\label{thm12:eq8}
\frac{d_i}{\psi_n(\lambda_i)}=1+o_p(1),
\end{equation}
\begin{equation}\label{thm12:eq9}
\frac{-\hat{\underline{s}}_{n,k}^{-1}(d_i)}{\lambda_i}=1+o_p(1).
\end{equation}
Therefore, for all $k'< k$, by \eqref{thm12:eq6}, \eqref{thm12:eq7}, \eqref{thm12:eq8} and \eqref{thm12:eq9}, we have
\begin{equation*}
\begin{aligned}
&\text{AGIC}(k')-\text{AGIC}(k)\\
=&\frac{n}{2}\sum_{i=k'+1}^{k}\log \frac{d_i}{1+\delta_i}+\gamma(k'-k)\{p-(k'+k)/2-1/2\}\\
=&\frac{n}{2}\sum_{i=k'+1}^{k}\{\log \frac{d_i}{1+\tilde{\delta}_i}-\log \frac{1+\delta_i}{1+\tilde{\delta}_i}\}+\gamma(k'-k)\{p-(k'+k)/2-1/2\}\\
=&\frac{n}{2}\sum_{i=k'+1}^{k}\{\log \frac{d_i}{1+\tilde{\delta}_i}-\log (1+\frac{\delta_i-\tilde{\delta}_i}{1+\tilde{\delta}_i})\}+\gamma(k'-k)\{p-(k'+k)/2-1/2\}\\
=&\frac{n}{2}\sum_{i=k'+1}^{k}\{\log \frac{d_i}{1+\tilde{\delta}_i}+o_p(1)\}+\gamma(k'-k)\{p-(k'+k)/2-1/2\}\\
=&\frac{n}{2}\sum_{i=k'+1}^{k}\{\log \lambda_i+\log \frac{-\hat{\underline{s}}_{n,k}^{-1}(d_i)}{\lambda_i}+o_p(1)\}+\gamma(k'-k)\{p-(k'+k)/2-1/2\}\\
=&\frac{n}{2}\sum_{i=k'+1}^{k}\{\log \lambda_i+o_p(1)\}+\gamma(k'-k)\{p-(k'+k)/2-1/2\}\\
\geq&\frac{n}{2}(k-k')\left[\log \lambda_k-2\gamma\{p-(k'+k)/2-1/2\}/n+o_p(1)\right].
\end{aligned}
\end{equation*}
Since $\gamma<\frac{1}{2c}\log \lambda_k$, it then holds under this case that
$$
\mathbb{P}\{\text{AGIC}(k')>\text{AGIC}(k)\}\rightarrow 1.
$$

Suppose $k'>k$.  Noting that
\[
\frac{d_i}{1+\delta_i}=\frac{d_i}{1-c_{n,i}+c_{n,i}\frac{1}{p-i}\sum_{j=i+1}^p\frac{d_i}{d_i-d_j}}=\frac{d_i}{1-c_{n,i}-c_{n,i}\hat{s}_{n,i}(d_i)d_i},
\]
we have
\begin{equation}\label{thm12:eq1}
\frac{d_i}{1+\delta_i}=-\hat{\underline{s}}_{n,i}^{-1}(d_i).
\end{equation}
For all $k<i=o(p)$, suppose that there exists an $i$ such that
\begin{equation}\label{thm12:eq2}
-\hat{\underline{s}}_{n,i}^{-1}(d_i)>(1+\sqrt{c})\lambda_{k+1}+o_p(1).
\end{equation}
Similar to the proof of Theorem \ref{thm1}, $\beta_{i+1}$ goes inside the support of $F^{c,H}$ in probability. Next, we proceed to demonstrate that $d_i$ is outside the support of $F^{c,H}$ in probability under the conditions that $d_i>\beta_{i+1}$ and  $d_i>d_{i+1}$.

Similar to the proof of \eqref{thm1:eq203}, it is seen that
\begin{equation}\label{thm12:eq0}
\frac{d_i\hat{\underline{s}}_{n,i}(d_i)}{d_i\underline{s}(d_i)}=1+o_p(1).
\end{equation}
By \eqref{thm12:eq0}, it is known that
\begin{equation*}
\begin{aligned}
d_i&=-\frac{1}{\underline{s}(d_i)}+c_{n,i}\int\frac{t}{1+t\underline{s}(d_i)}dH_{n,i}(t)\\
&=-\underline{s}^{-1}(d_i)+c_{n,i}\int\frac{-\underline{s}^{-1}(d_i)t}{-\underline{s}^{-1}(d_i)-t}dH_{n,i}(t)\\
&=-\hat{\underline{s}}_{n,i}^{-1}(d_i)\frac{d_i\hat{\underline{s}}_{n,i}(d_i)}{d_i\underline{s}(d_i)}+c_{n,i}\int\frac{-\hat{\underline{s}}_{n,i}^{-1}(d_i)\frac{d_i\hat{\underline{s}}_{n,i}(d_i)}{d_i\underline{s}(d_i)}t}{-\hat{\underline{s}}_{n,i}^{-1}(d_i)\frac{d_i\hat{\underline{s}}_{n,i}(d_i)}{d_i\underline{s}(d_i)}-t}dH_{n,i}(t)\\
&=-\hat{\underline{s}}_{n,i}^{-1}(d_i)+c_{n,i}\int\frac{-\hat{\underline{s}}_{n,i}^{-1}(d_i)t}{-\hat{\underline{s}}_{n,i}^{-1}(d_i)-t}dH_{n,i}(t)+o_p(1)\\
&=-\hat{\underline{s}}_{n,i}^{-1}(d_i)+c_{n,i}\int\frac{-\hat{\underline{s}}_{n,i}^{-1}(d_i)t}{-\hat{\underline{s}}_{n,i}^{-1}(d_i)-t}dH_{n,i}(t)+o_p(1).
\end{aligned}
\end{equation*}
By the definitions of $H_{n,i}(t)$ and $H_{n}(t)$, we have
\begin{equation*}
\begin{aligned}
d_i&=-\hat{\underline{s}}_{n,i}^{-1}(d_i)+c_{n,i}\int\frac{-\hat{\underline{s}}_{n,i}^{-1}(d_i)t}{-\hat{\underline{s}}_{n,i}^{-1}(d_i)-t}dH_{n,i}(t)+o_p(1)\\
&=-\hat{\underline{s}}_{n,i}^{-1}(d_i)+\frac{p-i}{n}\cdot\frac{1}{p-i}\sum_{j=i+1}^p\frac{-\hat{\underline{s}}_{n,i}^{-1}(d_i)\lambda_j}{-\hat{\underline{s}}_{n,i}^{-1}(d_i)-\lambda_j}\\
&=-\hat{\underline{s}}_{n,i}^{-1}(d_i)+\frac{p-k}{n}\cdot\frac{1}{p-k}(\sum_{j=k+1}^p\frac{-\hat{\underline{s}}_{n,i}^{-1}(d_i)\lambda_j}{-\hat{\underline{s}}_{n,i}^{-1}(d_i)-\lambda_j}
-\sum_{j=k+1}^i\frac{-\hat{\underline{s}}_{n,i}^{-1}(d_i)\lambda_j}{-\hat{\underline{s}}_{n,i}^{-1}(d_i)-\lambda_j})\\
&=-\hat{\underline{s}}_{n,i}^{-1}(d_i)+c_{n}\int\frac{-\hat{\underline{s}}_{n,i}^{-1}(d_i)t}{-\hat{\underline{s}}_{n,i}^{-1}(d_i)-t}dH_{n}(t)+o_p(1).
\end{aligned}
\end{equation*}
Combining this result with \eqref{thm12:eq2},  we have
\begin{equation}\label{thm12:eq3}
d_i>\psi_{n}((1+\sqrt{c})\lambda_{k+1})+o_p(1),
\end{equation}
which implies that $d_i$ is outside the support of $F^{c,H}$ in probability and $d_i>\beta_{i+1}$.

Since $\psi_{n}'((1+\sqrt{c})\lambda_{k+1})\geq 0$, we have $d_{k+1}\leq\beta_{k+1}\leq\psi_{n}((1+\sqrt{c})\lambda_{k+1})+o_p(1)$ (see, for example,  \cite{Bao:Pan:Zhou:2015} and \cite{Fan:Guo:Zheng:2022}). Thus, we have
\[
d_{i+1}\leq d_{k+1}\leq\psi_{n}((1+\sqrt{c})\lambda_{k+1})+o_p(1),
\]
which implies that $d_i>d_{i+1}$.

On the other hand,
for all $k<i=o(p)$, we have
\[
d_{i}\leq d_{k+1}\leq\psi_{n}((1+\sqrt{c})\lambda_{k+1})+o_p(1).
\]
This result is contradictory to \eqref{thm12:eq3}. Thus, for all $k<i=o(p)$, we have
\begin{equation}\label{thm12:eq4}
-\hat{\underline{s}}_{n,i}^{-1}(d_i)\leq(1+\sqrt{c})\lambda_{k+1}+o_p(1).
\end{equation}
Recall that $\lambda_{k+1}\leq 1$. By \eqref{thm12:eq4}, for all $k<i=o(p)$, we have
\begin{equation}\label{thm12:eq5}
-\hat{\underline{s}}_{n,i}^{-1}(d_i)\leq(1+\sqrt{c})\lambda_{k+1}+o_p(1)\leq (1+\sqrt{c})+o_p(1).
\end{equation}

Therefore, for all $k< k'=o(p)$, by \eqref{thm12:eq1} and \eqref{thm12:eq5}, we have
\begin{equation*}
\begin{aligned}
&\text{AGIC}(k')-\text{AGIC}(k)\\
=&-\frac{n}{2}\sum_{i=k+1}^{k'}\log \frac{d_i}{1+\delta_i}+\gamma(k'-k)\{p-(k+k')/2+1/2\}\\
=&-\frac{n}{2}\sum_{i=k+1}^{k'}\log (-\hat{\underline{s}}_{n,i}^{-1}(d_i))+\gamma(k'-k)\{p-(k+k')/2+1/2\}\\
\geq&-\frac{n}{2}(k'-k)\{\log (1+\sqrt{c})+o_p(1)\}+\gamma(k'-k)\{p-(k+k')/2+1/2\}\\
=&\frac{n}{2}(k'-k)\left[2\gamma\{p-(k+k')/2+1/2\}/n-\log (1+\sqrt{c})+o_p(1)\right].
\end{aligned}
\end{equation*}
Since $\gamma>\frac{1}{2c}\log (1+\sqrt{c})$, it then holds under this case that
$$
\mathbb{P}\{\text{AGIC}(k')>\text{AGIC}(k)\}\rightarrow 1.
$$

\eop

\subsection{Proof of Theorem \ref{thm7}}
The proof is similar to that of Theorem \ref{thm6}. For $i=k$, we have
\[
\log \frac{d_k}{1+\delta_k}=\log \frac{d_k}{1+\tilde{\delta}_k}=\log \lambda_k+o_p(1)>\gamma.
\]
For all $k<i=o(p)$, we have
\[
\log \frac{d_i}{1+\delta_i}\leq \log (1+\sqrt{c})+o_p(1).
\]
Thus, we have
\[
\mathbb{P}(\hat{k}=k)\rightarrow1.
\]

\eop

\subsection{Proof of Lemma \ref{lemma:diagonal}}
\label{sec:proof2}
In what follows, we give the proof for the case of $\mathbb{E}Y_{ij}^4<\infty$.  That is, by Lemma \ref{corollary:large deviation2},  we can show that Lemma \ref{lemma:diagonal} holds for $k=o(n^{1/4})$. For the case of $Y_{ij}$ is Gaussian,  instead of using Lemma \ref{corollary:large deviation2},  we can show Lemma \ref{lemma:diagonal} holds for $k=o(n^{1/3})$ using Lemma \ref{lemma:large deviation1}. The proof is very similar to the case of $\mathbb{E}Y_{ij}^4<\infty$, and the details are thus omitted.

By definition, each $d_i$ solves the equation
\[
0=\mid d \I-\S_n\mid=\mid d \I-\S_{22}\mid \cdot \mid d \I-\mathbf K_n(d)\mid,
\]
where
\[
\mathbf K_n(d)=\frac{1}{n}\Z_1^\top(\I+\H_n(d))\Z_1,
\]
\[
\H_n(d)=\frac{1}{n}\Z_2(d \I-\S_{22})^{-1}\Z_2^\top.
\]
Using the same argument as in \eqref{thm1:eq11} and \eqref{thm1:eq70}, we have
\begin{equation}
\label{thm2:eq10}
\left|\bLambda_{1}^{-1}d_i\I-\frac{1}{n}\bLambda_{1}^{-\frac{1}{2}}\Z_1^\top(\I+\H_n(d_i))\Z_1\bLambda_{1}^{-\frac{1}{2}}\right|=0,
\end{equation}
and
\begin{equation*}
\bLambda_{1}^{-1}d_i\I-\frac{1}{n}\bLambda_{1}^{-\frac{1}{2}}\Z_1^\top(\I+\H_n(d_i))\Z_1\bLambda_{1}^{-\frac{1}{2}}
=\bLambda_{1}^{-1}d_i\I-\frac{1}{n}\tr(\I+\H_n(d))\I+\M_n(d_i),
\end{equation*}
where $\M_n(d_i)=\frac{1}{n}\tr(\I+\H_n(d_i))\I-\frac{1}{n}\bLambda_{1}^{-\frac{1}{2}}\Z_1^\top(\I+\H_n(d_i))\Z_1\bLambda_{1}^{-\frac{1}{2}}$.
Similar to the proof of \eqref{thm1:eq40}, we have \begin{equation}\label{thm2:eq2}
\M_n(d_i)=\V(d_i)\Diag\{v_1(d_i),\ldots,v_k(d_i)\}\V(d_i)^\top,
\end{equation}
where $\V(d_i)$ is an orthogonal matrix and $\max_{1\leq j\leq k}\vert v_j(d_i)\vert=\mathcal{O}_p(k/\sqrt{n})$.
Since
 \[
\bLambda_{1}^{-1}=\Diag\left\{\frac{1}{\lambda_1},\ldots,\frac{1}{\lambda_k}\right\},
\]
\eqref{thm2:eq10} becomes
\begin{eqnarray*}
&&\left|\bLambda_{1}^{-1}d_i\I-\frac{1}{n}\bLambda_{1}^{-\frac{1}{2}}\Z_1^\top(\I+\H_n(d_i))\Z_1\bLambda_{1}^{-\frac{1}{2}}\right|\\\nonumber
&=&\left|\Diag\left\{\frac{d_i-\frac{1}{n}\tr(\I+\H_n(d_i))\lambda_1}{\lambda_1},\ldots,\frac{d_i-\frac{1}{n}\tr(\I+\H_n(d_i))\lambda_k}{\lambda_k}\right\}+\M_n(d_i)\right|=0.
\end{eqnarray*}
\if
Next, we write
\begin{equation}\label{thm2:eq2}
\begin{aligned}
\B(d_i)=&\U^\top\M_n(d_i)\U\\
=&\U^\top\V(d_i)\Diag\{\tilde{v}_1(d_i),\ldots,\tilde{v}_k(d_i)\}\V(d_i)^\top\U\\
\triangleq&\Q(d_i)\Diag\{\tilde{v}_1(d_i),\ldots,\tilde{v}_k(d_i)\}\Q(d_i)^\top,
\end{aligned}
\end{equation}
where $\Q(d_i)=\U^\top\V(d_i)$ is an orthogonal matrix.
\fi

To ease notation, we write $\A_n(d_i)=(A_{st})_{k\times k}$, $\M_n(d_i)=(M_{st})_{k\times k}$ and $\V(d_i)=(V_{st})_{k\times k}$ without emphasizing their dependence on $d_i$, when there is no ambiguity.
For any $1\leq s\neq t\leq k$, by the definitions of $\A_n(d_i)$ and $\M_n(d_i)$, we have $A_{st}=M_{st}$.

Similar to the proof of \eqref{thm1:eq8} in Theorem \ref{thm1}, we have
\begin{equation}\label{thm2:eq3}
\begin{aligned}
\sum_{s=1}^k\sum_{t=1}^kM_{st}^2=\tr(\M^2(d_i))=\mathcal{O}_p(k^3/n).
\end{aligned}
\end{equation}

By \eqref{thm2:eq2}, we also have
\begin{equation}\label{thm2:eq4}
\begin{aligned}
&\sum_{t=1}^kM_{st}^2=\sum_{t=1}^k(\sum_{l=1}^kv_l(d_i)V_{sl}V_{tl})^2\\
=&\sum_{t=1}^k\sum_{l=1}^kv_l^2(d_i)V_{sl}^2V_{tl}^2+\sum_{t=1}^k\sum_{l=1}^k\sum_{j\neq l}^kv_l(d_i)v_s(d_i)V_{sl}V_{tl}V_{sj}V_{tj}\\
=&\sum_{t=1}^k\sum_{l=1}^kv_l^2(d_i)V_{sl}^2V_{tl}^2+\sum_{l=1}^k\sum_{j\neq l}^kv_l(d_i)v_s(d_i)V_{sl}V_{sj}\sum_{t=1}^kV_{tl}V_{tj}\\
=&\sum_{t=1}^k\sum_{l=1}^kv_l^2(d_i)V_{sl}^2V_{tl}^2=\sum_{l=1}^kv_l^2(d_i)V_{sl}^2\sum_{t=1}^kV_{tl}^2\\
=&\sum_{l=1}^kv_l^2(d_i)V_{sl}^2\leq\max_{1\leq l\leq k}v_l^2(d_i)\sum_{l=1}^kV_{sl}^2=\mathcal{O}_p(k^2/n).
\end{aligned}
\end{equation}

Using the Laplace formula, we can write the determinant of $\A_n(d_i)$ as
\begin{equation}
\label{thm2:eq5}
|\A_n(d_i)|=\sum_{t=1}^k(-1)^{i+t}A_{it}h_{it}=0,
\end{equation}
where $h_{it}=|\H_{it}|$ and $\H_{it}$ is a $(k-1)\times (k-1)$ matrix that results from deleting the $i$-th row and the $t$-th column of $\A_n(d_i)$.

Next, we investigate $h_{ii}$ and $h_{it}$, $i\neq t$. Note that
$$
f_j=f_j(d_i)\triangleq\frac{d_i-\frac{1}{n}\tr(\I+\H_n(d_i))\lambda_j}{\lambda_j},\quad j\neq i,
$$
is uniformly bounded away from zero.
For the case of  $t=i$, without loss of generality, we assume that $i=1$ and $t=1$.
By the definitions of $\A_n(d_i)$ and $\H_{it}$, it holds that
\begin{equation}\label{eqn:de02}
\H_{11}=\F_{2}+\M_{22},
\end{equation}
where $\F_{2}=\Diag\left\{f_2,\dots,f_k\right\}$. Noting that
\begin{equation*}
\begin{pmatrix}
0 \,\,& \0^{\top} \\
\0 \,\,&  \M_{22}
\end{pmatrix}
  =\begin{pmatrix}
0 \,\,& \0^{\top} \\
\0 \,\,&  \I_{k-1}
\end{pmatrix}
\M_n(d_i)\begin{pmatrix}
0 \,\,& \0^{\top} \\
\0 \,\,&  \I_{k-1}
\end{pmatrix},
\end{equation*}
by \eqref{thm2:eq2}, we have $||\M_{22}||=\mathcal{O}_p(k/\sqrt{n})$. By Lemma \ref{lemma:weyl} and \eqref{eqn:de02}, we have
\begin{equation}\label{eqn:de03}
\F_2'\triangleq\Diag\{f_2+\epsilon_2',\ldots,f_k+\epsilon_k'\}=\V_2(d_i)^\top(\F_{2}+\M_{22})\V_2(d_i),
\end{equation}
where $\V_2(d_i)$ is an orthogonal matrix, $\max_{2\leq j\leq k}|\epsilon_j'|=\mathcal{O}_p(k/\sqrt{n})$, $\sum_{2\leq j\leq k}|\epsilon_j'|=\mathcal{O}_p(k^2/\sqrt{n})$ and $\sum_{2\leq j\leq k}\epsilon_j'^2=\mathcal{O}_p(k^3/n)$.
By \eqref{eqn:de02} and \eqref{eqn:de03}, we have
$h_{11}=|\H_{11}|=\prod_{j=2}^k(f_j+\epsilon_j')$.
Thus, for $i=t$, we have
\begin{equation}\label{eqn:de04}
h_{ii}=\prod_{j\neq i}(f_j+\epsilon_j'),
\end{equation}
where  $\max_{j\neq i}|\epsilon_j'|=\mathcal{O}_p(k/\sqrt{n})$, $\sum_{j\neq i}|\epsilon_j'|=\mathcal{O}_p(k^2/\sqrt{n})$ and $\sum_{j\neq i}\epsilon_j'^2=\mathcal{O}_p(k^3/n)$.

Correspondingly, we have
\[
|\sum_{t\neq i}\frac{1}{f_t}\epsilon_t'|\leq\frac{1}{\min_{t\neq i}|f_t|}\sum_{t\neq i}|\epsilon_t'|=\mathcal{O}_p(k^2/\sqrt{n}),
\]
\[
|\sum_{t\neq i}\sum_{s\neq i,t<s\leq k}\frac{1}{f_tf_s}\epsilon_t'\epsilon_s'|\leq\frac{1}{\min_{t\neq i,s\neq i,t<s\leq k}|f_tf_s|}\sum_{t\neq i}\sum_{s\neq i,t<s\leq k}|\epsilon_t'\epsilon_s'|=\mathcal{O}_p(k^4/n).
\]
Next, we have
\begin{equation}\label{thm2:eq6}
\begin{aligned}
h_{ii}=\prod_{j\neq i}(f_j+\epsilon_j')&=\prod_{j\neq i}f_j+\prod_{j\neq i}f_j\sum_{t\neq i}\frac{1}{f_t}\epsilon_t'+\prod_{j\neq i}f_j\sum_{t\neq i}\sum_{s\neq i,t<s\leq k}\frac{1}{f_tf_s}\epsilon_t'\epsilon_s'+\cdots\\
&=\prod_{j\neq i}f_j+\mathcal{O}_p(k^2/\sqrt{n})\prod_{j\neq i}f_j+\mathcal{O}_p(k^4/n)\prod_{j\neq i}f_j+\cdots\\
&=(1+\mathcal{O}_p(k^2/\sqrt{n}))\prod_{j\neq i}f_j.
\end{aligned}
\end{equation}

For the case of  $t\neq i$, without loss of generality, we assume that $i=1$ and $t=2$.
By the definitions of $\A_n(d_i)$ and $\H_{it}$, it holds that
\begin{equation}\label{eqn:de2}
\H_{12}=\begin{pmatrix}
 M_{21} \,\,& \boldsymbol{\alpha}^{\top} \\
  \boldsymbol{\beta} \,\,&  \F_{3}+\M_{33}
\end{pmatrix},
\end{equation}
where $\boldsymbol{\alpha}=(M_{23},\dots,M_{2k})^{\top}$, $\boldsymbol{\beta}=(M_{31},\dots,M_{k1})^{\top}$, and
$\F_{3}=\Diag\left\{f_3,\dots,f_k\right\}$.
Similar to the proof of \eqref{eqn:de04}, we have
\begin{equation}\label{eqn:de3}
\F_3''\triangleq\Diag\{f_3+\epsilon_3'',\ldots,f_k+\epsilon_k''\}=\V_3(d_i)^\top(\F_{3}+\M_{33})\V_3(d_i),
\end{equation}
where $\V_3(d_i)$ is an orthogonal matrix, $\max_{3\leq j\leq k}|\epsilon_j''|=\mathcal{O}_p(k/\sqrt{n})$, $\sum_{3\leq j\leq k}|\epsilon_j''|=\mathcal{O}_p(k^2/\sqrt{n})$ and $\sum_{3\leq j\leq k}\epsilon_j''^2=\mathcal{O}_p(k^3/n)$.
By \eqref{eqn:de2} and \eqref{eqn:de3}, we have
\begin{equation}\label{eqn:de4}
\begin{aligned}
h_{12}=|\H_{12}|&=|\begin{pmatrix}
 1 \,\,& \0^{\top} \\
 \0\,\,&  \V_3(d_i)^\top
\end{pmatrix}
\begin{pmatrix}
 M_{21} \,\,& \boldsymbol{\alpha}^{\top} \\
  \boldsymbol{\beta} \,\,&  \F_{3}+\M_{33}
\end{pmatrix}
\begin{pmatrix}
 1 \,\,& \0^{\top} \\
 \0\,\,&  \V_3(d_i)
\end{pmatrix}|\\
&=|\begin{pmatrix}
 M_{21} \,\,& \boldsymbol{\alpha}^{\top}\V_3(d_i) \\
  \V_3(d_i)^\top\boldsymbol{\beta} \,\,&  \F_{3}
\end{pmatrix}|\\
&=|\begin{pmatrix}
 M_{21}'' \,\,& \boldsymbol{\alpha''}^{\top}\\
  \boldsymbol{\beta}'' \,\,&  \F_{3}''
\end{pmatrix}|=(M_{21}''-\sum_{j=3}^k\frac{M_{2j}''M_{j2}''}{f_j+\epsilon_j''})\prod_{j=3}^k(f_j+\epsilon_j''),
\end{aligned}
\end{equation}
where $\boldsymbol{\alpha}''=\V_3(d_i)^\top\boldsymbol{\alpha}\triangleq(M_{23}'',\dots,M_{2k}'')^{\top}$, $\boldsymbol{\beta}''=\V_3(d_i)^\top\boldsymbol{\beta}\triangleq(M_{31}'',\dots,M_{k1}'')^{\top}$.

Since  $\V_3(d_i)$ is an orthogonal matrix, we have
\begin{equation}\label{eqn:de5}
\sum_{j=3}^kM_{j1}''^2=\sum_{j=3}^kM_{j1}^2,\,\,\,\sum_{j=3}^kM_{2j}''^2=\sum_{j=3}^kM_{2j}^2.
\end{equation}
By \eqref{thm2:eq4} and \eqref{eqn:de5}, we have
\begin{equation}\label{eqn:de6}
\begin{aligned}
\sum_{j=3}^kM_{j1}''^2=\sum_{j=3}^kM_{j1}^2\leq\sum_{j=1}^kM_{j1}^2=\sum_{j=1}^kM_{1j}^2=\mathcal{O}_p(k^2/n),
\end{aligned}
\end{equation}
and
\begin{equation}\label{eqn:de11}
\begin{aligned}
\sum_{j=3}^kM_{2j}''^2=\sum_{j=3}^kM_{2j}^2\leq\sum_{j=1}^kM_{2j}^2=\mathcal{O}_p(k^2/n).
\end{aligned}
\end{equation}
By \eqref{eqn:de6}, \eqref{eqn:de11}, and the fact that $M_{21}=M_{21}''$, we have
\begin{equation}\label{eqn:de7}
\begin{aligned}
|M_{21}''-\sum_{j=3}^k\frac{M_{2j}''M_{j2}''}{f_j+\epsilon_j''}|&\leq|M_{21}''|+|\sum_{j=3}^k\frac{M_{2j}''M_{j1}''}{f_j+\epsilon_j''}|\leq|M_{21}''|+\sum_{j=3}^k\frac{|M_{2j}''M_{j1}''|}{|f_j+\epsilon_j''|}\\
&\leq|M_{21}''|+\lambda_0\sum_{j=3}^k|M_{2j}''M_{j1}''|\leq |M_{21}''|+\lambda_0\sqrt{\sum_{j=3}^kM_{2j}''^2\sum_{j=3}^kM_{j1}''^2}\\
&=|M_{21}|+\mathcal{O}_p(k^2/n),
\end{aligned}
\end{equation}
where $\lambda_0=\max_{3\leq j\leq k}\frac{1}{|f_j+\epsilon_j''|}$.
By \eqref{thm2:eq4}, we have $\max_{1\leq s\leq k}\sum_{t=1}^kM_{st}^2=\mathcal{O}_p(k^2/n)$, which implies that
\begin{equation}\label{eqn:de8}
\begin{aligned}
|M_{21}|&\leq\max_{1\leq s\leq k}\max_{1\leq t\leq k}|M_{st}|\leq\max_{1\leq s\leq k}\sqrt{\sum_{t=1}^kM_{st}^2}=\mathcal{O}_p(k/\sqrt{n}).
\end{aligned}
\end{equation}
Combining \eqref{eqn:de7} and \eqref{eqn:de8}, we have
\begin{equation}\label{eqn:de9}
\begin{aligned}
|M_{21}''-\sum_{j=3}^k\frac{M_{2j}''M_{j2}''}{f_j+\epsilon_j''}|&=\mathcal{O}_p(k/\sqrt{n}).
\end{aligned}
\end{equation}
By \eqref{eqn:de3}, \eqref{eqn:de4} and \eqref{eqn:de9}, we have
$h_{12}=\epsilon_2''\prod_{j=3}^k(f_j+\epsilon_j'')$,
where $\epsilon_2''=\mathcal{O}_p(k/\sqrt{n})$, $\max_{3\leq j\leq k}|\epsilon_j''|=\mathcal{O}_p(k/\sqrt{n})$, $\sum_{3\leq j\leq k}|\epsilon_j''|=\mathcal{O}_p(k^2/\sqrt{n})$ and $\sum_{3\leq j\leq k}\epsilon_j''^2=\mathcal{O}_p(k^3/n)$.
Thus, for $t\neq i$, we have
\begin{equation}\label{eqn:de10}
h_{it}=\epsilon_t''\prod_{j\neq i, t}(f_j+\epsilon_j''),
\end{equation}
where $\epsilon_t''=\mathcal{O}_p(k/\sqrt{n})$, $\max_{j\neq i, t}|\epsilon_j''|=\mathcal{O}_p(k/\sqrt{n})$, $\sum_{j\neq i, t}|\epsilon_j''|=\mathcal{O}_p(k^2/\sqrt{n})$  and $\sum_{j\neq i, t}\epsilon_j''^2=\mathcal{O}_p(k^3/n)$.

Correspondingly, we have
\begin{eqnarray*}
\sum_{s\neq i,t}\frac{1}{|f_s|}|\epsilon_s''|&\leq&\frac{1}{\min_{s\neq i,t}|f_s|}\sum_{s\neq i,t}|\epsilon_s''|=\mathcal{O}_p(k^2/\sqrt{n}),
\end{eqnarray*}
and
\begin{eqnarray*}&&\sum_{s\neq i,t}\sum_{u\neq i,t, s<u\leq k}\frac{1}{|f_sf_u|}|\epsilon_s''\epsilon_u''|\\
&\leq&\frac{1}{\min_{s\neq i,t,u\neq i,t,s<u\leq k}|f_sf_u|}\sum_{s\neq i,t}|\epsilon_s''|\sum_{u\neq i,t,s<u\leq k}|\epsilon_u''|=\mathcal{O}_p(k^4/n).
\end{eqnarray*}
Next, we have
\begin{equation*}
\begin{aligned}
\prod_{j\neq i,t}(|f_j|+|\epsilon_j''|)&=\prod_{j\neq i,t}|f_j|+\prod_{j\neq i,t}|f_j|\sum_{s\neq i,t}\frac{1}{|f_s|}|\epsilon_s''|+\prod_{j\neq i,t}|f_j|\sum_{s\neq i,t}\sum_{u\neq i,t,s<u\leq k}\frac{1}{|f_sf_u|}|\epsilon_s''\epsilon_u''|+\cdots\\
&=\prod_{j\neq i,t}|f_j|+\mathcal{O}_p(k^2/\sqrt{n})\prod_{j\neq i,t}|f_j|+\mathcal{O}_p(k^4/n)\prod_{j\neq i,t}|f_j|+\cdots\\
&=(1+\mathcal{O}_p(k^2/\sqrt{n}))\prod_{j\neq i,t}|f_j|\leq\frac{1+\mathcal{O}_p(k^2/\sqrt{n})}{\min_{j\neq i}|f_j|}\prod_{j\neq i}|f_j|.
\end{aligned}
\end{equation*}
Thus, we have
\begin{equation}\label{thm2:eq7}
|h_{it}|\leq\frac{1+\mathcal{O}_p(k^2/\sqrt{n})}{\min_{j\neq i}|f_j|}\prod_{j\neq i}|f_j||\epsilon_t''|.
\end{equation}
Combining \eqref{thm2:eq4}, \eqref{eqn:de10}, \eqref{thm2:eq7}, and the fact that $A_{st}=M_{st}$, $1\leq s\neq t\leq k $, we have
\begin{equation}\label{thm2:eq8}
\begin{aligned}
\sum_{t\neq i}^k|A_{it}h_{it}|&\leq\frac{1+\mathcal{O}_p(k^2/\sqrt{n})}{\min_{j\neq i}|f_j|}\prod_{j\neq i}|f_j|\sum_{t\neq i}|A_{it}\epsilon_t''|\\
&\leq\frac{1+\mathcal{O}_p(k^2/\sqrt{n})}{\min_{j\neq i}|f_j|}\prod_{j\neq i}|f_j|\sqrt{\sum_{t\neq i}A_{it}^2\sum_{t\neq i}\epsilon_t''^2}\\
&\leq\frac{1+\mathcal{O}_p(k^2/\sqrt{n})}{\min_{j\neq i}|f_j|}\prod_{j\neq i}|f_j|\sqrt{\sum_{t\neq i}M_{it}^2\sum_{t\neq i}\epsilon_t''^2}\\
&\leq\frac{1+\mathcal{O}_p(k^2/\sqrt{n})}{\min_{j\neq i}|f_j|}\prod_{j\neq i}|f_j|\sqrt{\sum_{t\neq i}M_{it}^2\sum_{t\neq i}\epsilon_t''^2}\\
&=\frac{1+\mathcal{O}_p(k^2/\sqrt{n})}{\min_{j\neq i}|f_j|}\sqrt{\mathcal{O}_p(\frac{k^{2}}{n})\mathcal{O}_p(\frac{k^{3}}{n})}\prod_{j\neq i}|f_j|\\
&=\frac{1+\mathcal{O}_p(k^2/\sqrt{n})}{\min_{j\neq i}|f_j|}\mathcal{O}_p(\sqrt{\frac{k^5}{n^2}})\prod_{j\neq i}|f_j|.
\end{aligned}
\end{equation}
By \eqref{thm2:eq5},  \eqref{thm2:eq6}, and \eqref{thm2:eq8}, we have
\begin{eqnarray*}
|A_{ii}|&=&\left|\frac{\sum_{t\neq i}^k(-1)^{i+t}A_{it}h_{it}}{-h_{ii}}\right|\leq\frac{1+\mathcal{O}_p(\sqrt{k^4/n})}{1+\mathcal{O}_p(\sqrt{k^4/n})}\frac{1}{\min_{j\neq i}|f_j|}\mathcal{O}_p(\sqrt{\frac{k^5}{n^2}})\\
&=&\mathcal{O}_p(\frac{k^{5/2}}{n}).
\end{eqnarray*}

\eop

\subsection{Proof of Lemma \ref{lemma:diagonal11}}

Let $S(i)=\{i_1,\dots,i_{|S(i)|}\}$ and $T=\{t_1,\dots,t_{|S(i)|}\}$. Using the Laplace formula, we can write the determinant of $\A_n(d_i)$ as
\begin{equation}
\label{lemma11:eq1}
|\A_n(d_i)|=\sum_{T=\{t_1,\dots,t_{|S(i)|}\}}(-1)^{i_1+\dots+i_{|S(i)|}+t_1+\dots+t_{|S(i)|}}|\A_{S(i)\times T}|h_{S(i)\times T}=0,
\end{equation}
where $h_{S(i)\times T}=|\H_{S(i)\times T}|$ and $\H_{S(i)\times T}$ is a $(k-|S(i)|)\times (k-|S(i)|)$ matrix that results from deleting the $S(i)$-th rows and the $T$-th columns of $\A_n(d_i)$.

Next, we investigate $\H_{S(i)\times S(i)}$ and $\H_{S(i)\times T}$, $T\neq S(i)$. Note that
$$
f_j\triangleq\frac{d_i-\frac{1}{n}\tr(\I+\H_n(d_i))\lambda_j}{\lambda_j},\quad j\notin S(i),
$$
is uniformly bounded away from zero.
For the case of  $T=S(i)$, similar to the proof of \eqref{eqn:de04} in Lemma \ref{lemma:diagonal}, we have
\begin{equation}\label{eqn:de001}
h_{S(i)\times S(i)}=\prod_{j\notin S(i)}(f_j+\epsilon_j'),
\end{equation}
where $\max_{j\notin S(i)}|\epsilon_j'|=\mathcal{O}_p(k/\sqrt{n})$, $\sum_{j\notin S(i)}|\epsilon_j'|=\mathcal{O}_p(k^2/\sqrt{n})$ and $\sum_{j\notin S(i)}\epsilon_j'^2=\mathcal{O}_p(k^3/n)$.
Correspondingly, we have
$$
|\sum_{t\notin S(i)}\frac{1}{f_t}\epsilon_t'|\leq\frac{1}{\min_{t\notin S(i)}|f_t|}\sum_{t\notin S(i)}|\epsilon_t'|=\mathcal{O}_p(k^2/\sqrt{n}),
$$
$$
|\sum_{t\notin S(i)}\sum_{s\notin S(i),t< s\leq k}\frac{1}{f_tf_s}\epsilon_t'\epsilon_s'|\leq\frac{1}{\min_{t\notin S(i),s\notin S(i),t< s\leq k}|f_tf_s|}\sum_{t\notin S(i)}\sum_{s\notin S(i),t< s\leq k}|\epsilon_t'\epsilon_s'|=\mathcal{O}_p(k^4/n).
$$
Next, we have
\begin{equation}\label{lemma11:eq2}
\begin{aligned}
h_{S(i)\times S(i)}=\prod_{j\notin S(i)}(f_j+\epsilon_j')&=\prod_{j\notin S(i)}f_j+\prod_{j\notin S(i)}f_j\sum_{t\notin S(i)}\frac{1}{f_t}\epsilon_t'+\prod_{j\notin S(i)}f_j\sum_{t\notin S(i)}\sum_{s\notin S(i),t< s\leq k}\frac{1}{f_tf_s}\epsilon_t'\epsilon_s'+\cdots\\
&=\prod_{j\notin S(i)}f_j+\mathcal{O}_p(k^2/\sqrt{n})\prod_{j\notin S(i)}f_j+\mathcal{O}_p(k^4/n)\prod_{j\notin S(i)}f_j+\cdots\\
&=(1+\mathcal{O}_p(k^2/\sqrt{n}))\prod_{j\notin S(i)}f_j.
\end{aligned}
\end{equation}
For the case of  $T\neq S(i)$, let $T_1=S(i)\cap T$ and $T_2=T/T_1$. Similar to the proof of \eqref{eqn:de10} in Lemma \ref{lemma:diagonal}, we have
\begin{equation}\label{eqn:de002}
\begin{aligned}
h_{S(i)\times T}=\prod_{t\in T_2}\epsilon_t''\prod_{j\notin S(i)\cup T_2}(f_j+\epsilon_j''),
\end{aligned}
\end{equation}
where $\max_{t\in T_2}|\epsilon_t''|=\mathcal{O}_p(k/\sqrt{n})$, $\max_{j\notin S(i)\cup T_2}|\epsilon_j''|=\mathcal{O}_p(k/\sqrt{n})$, $\sum_{j\notin S(i)\cup T_2}|\epsilon_j''|=\mathcal{O}_p(k^2/\sqrt{n})$  and $\sum_{j\notin S(i)\cup T_2}\epsilon_j''^2=\mathcal{O}_p(k^3/n)$.

Correspondingly, we have
\[
\sum_{s\notin S(i)\cup T_2}\frac{1}{|f_s|}|\epsilon_s''|\leq\frac{1}{\min_{s\notin S(i)\cup T_2}|f_s|}\sum_{s\notin S(i)\cup T_2}|\epsilon_s''|=\mathcal{O}_p(k^2/\sqrt{n}),
\]
and
\begin{eqnarray*}&&\sum_{s\notin S(i)\cup T_2}\sum_{u\notin S(i)\cup T_2,s< u\leq k}\frac{1}{|f_sf_u|}|\epsilon_s''\epsilon_u''|\\
&\leq&\frac{1}{\min_{s\notin S(i)\cup T_2, u\notin S(i)\cup T_2, s< u\leq k}|f_sf_u|}\sum_{s\notin S(i)\cup T_2}|\epsilon_s''|\sum_{u\notin S(i)\cup T_2,s< u\leq k}|\epsilon_u''|=\mathcal{O}_p(k^4/n).
\end{eqnarray*}
Next, we have
\begin{eqnarray*}
&&\prod_{j\notin S(i)\cup T_2}(|f_j|+|\epsilon_j''|)\\
&=&\prod_{j\notin S(i)\cup T_2}|f_j|+\prod_{j\notin S(i)\cup T_2}|f_j|\sum_{s\notin S(i)\cup T_2}\frac{1}{|f_s|}|\epsilon_s''|\\
&&\hspace{1in}+\prod_{j\notin S(i)\cup T_2}|f_j|\sum_{s\notin S(i)\cup T_2}\sum_{u\notin S(i)\cup T_2,1\leq u<s}\frac{1}{|f_sf_u|}|\epsilon_s''||\epsilon_u''|+\cdots\\
&=&\prod_{j\notin S(i)\cup T_2}|f_j|+\mathcal{O}_p(k^2/\sqrt{n})\prod_{j\notin S(i)\cup T_2}|f_j|+\mathcal{O}_p(k^4/n)\prod_{j\notin S(i)\cup T_2}|f_j|+\cdots\\
&=&(1+\mathcal{O}_p(k^2/\sqrt{n}))\prod_{j\notin S(i)\cup T_2}|f_j|\leq\frac{1+\mathcal{O}_p(k^2/\sqrt{n})}{(\min_{j\notin S(i)}|f_j|)^{|T_2|}}\prod_{j\notin S(i)}|f_j|.
\end{eqnarray*}
Hence,
\begin{equation}\label{lemma11:eq4}
|h_{S(i)\times T}|\leq\frac{1+\mathcal{O}_p(k^2/\sqrt{n})}{(\min_{j\notin S(i)}|f_j|)^{|T_2|}}\prod_{j\notin S(i)}|f_j|\prod_{t\in T_2}|\epsilon_{t}''|.
\end{equation}
Recall that $|S(i)|$ is the cardinality of $S(i)$ and $|S(i)|=|T|$, then we have
\begin{equation*}
\begin{aligned}
&\sum_{T\neq S(i)}\Big||\A_{S(i)\times T}|h_{S(i)\times T}\Big|\\
&\leq\frac{1+\mathcal{O}_p(k^2/\sqrt{n})}{(\min_{j\notin S(i)}|f_j|)^{|T_2|}}\prod_{j\notin S(i)}|f_j|\sum_{T\neq S(i)}(\max_{u\in S(i)}\sum_{v\in T}|A_{uv}|)^{|S(i)|}\prod_{t\in T_2}|\epsilon_{t}''|\\
&\leq\frac{1+\mathcal{O}_p(k^2/\sqrt{n})}{(\min_{j\notin S(i)}|f_j|)^{|T_2|}}\prod_{j\notin S(i)}|f_j|\sum_{T\neq S(i)}(\max_{u\in S(i)}\sqrt{|T|\sum_{v\in T}A_{uv}^2})^{|S(i)|}\prod_{t\in T_2}|\epsilon_{t}''|\\
&\leq\frac{1+\mathcal{O}_p(k^2/\sqrt{n})}{(\min_{j\notin S(i)}|f_j|)^{|S(i)|}}\prod_{j\notin S(i)}|f_j|\sum_{T\neq S(i)}(\sqrt{|S(i)|\max_{u\in S(i)}\sum_{v\in T}A_{uv}^2})^{|S(i)|}\prod_{t\in T_2}|\epsilon_{t}''|\\
&\leq\frac{1+\mathcal{O}_p(k^2/\sqrt{n})}{(\min_{j\notin S(i)}|f_j|)^{|S(i)|}}\prod_{j\notin S(i)}|f_j|\sum_{T\neq S(i)}(\sqrt{|S(i)|\max_{u\in S(i)}(\sum_{v\in T_1}A_{uv}^2+\sum_{v\in T_2}A_{uv}^2)})^{|S(i)|}\prod_{t\in T_2}|\epsilon_{t}''|.
\end{aligned}
\end{equation*}
Combining \eqref{thm2:eq4}, \eqref{lemma11:eq1}, \eqref{lemma11:eq4}, and the fact that $A_{st}=M_{st}$, $1\leq s\neq t\leq k $, we have
\begin{equation*}
\begin{aligned}
&\sum_{T\neq S(i)}\Big||\A_{S(i)\times T}|h_{S(i)\times T}\Big|\\
&\leq\frac{1+\mathcal{O}_p(k^2/\sqrt{n})}{(\min_{j\notin S(i)}|f_j|)^{|S(i)|}}\prod_{j\notin S(i)}|f_j|\sum_{T\neq S(i)}(\sqrt{|S(i)|\max_{u\in S(i)}(\sum_{v\in T_1}A_{uv}^2+\sum_{v\in T_2}M_{uv}^2)})^{|S(i)|}\prod_{t\in T_2}|\epsilon_{t}''|\\
&\leq\frac{1+\mathcal{O}_p(k^2/\sqrt{n})}{(\min_{j\notin S(i)}|f_j|)^{|S(i)|}}\prod_{j\notin S(i)}|f_j|\\
&\times\sum_{T\neq S(i)}(\sqrt{|S(i)|(\sum_{u\in S(i)}\sum_{v\in T_1}A_{uv}^2+\mathcal{O}_p(k^2/n))})^{|S(i)|}(\mathcal{O}_p(k/\sqrt{n}))^{|S(i)|}\\
&\leq\frac{1+\mathcal{O}_p(k^2/\sqrt{n})}{(\min_{j\notin S(i)}|f_j|)^{|S(i)|}}\prod_{j\notin S(i)}|f_j|\\
&\times\sum_{T\neq S(i)}(\sqrt{|S(i)|(\sum_{u\in T_1}A_{uu}^2+\sum_{u\in S(i)}\sum_{v\in T_1, u\neq v}A_{uv}^2+\mathcal{O}_p(k^2/n))})^{|S(i)|}(\mathcal{O}_p(k/\sqrt{n}))^{|S(i)|}\\
&\leq\frac{1+\mathcal{O}_p(k^2/\sqrt{n})}{(\min_{j\notin S(i)}|f_j|)^{|S(i)|}}\prod_{j\notin S(i)}|f_j|\\
&\times\sum_{T\neq S(i)}(\sqrt{|S(i)|(\sum_{u\in T_1}A_{uu}^2+\sum_{u\in S(i)}\sum_{v\in T_1, u\neq v}M_{uv}^2+\mathcal{O}_p(k^2/n))})^{|S(i)|}(\mathcal{O}_p(k/\sqrt{n}))^{|S(i)|}
\end{aligned}
\end{equation*}
Recall $S(i)=\{1\leq j\leq k: \lambda_j=\lambda_i\}$.  For $j\in S(i)$, we have
\[
f_j=\frac{d_i-\frac{1}{n}\tr(\I+\H_n(d_i))\lambda_i}{\lambda_i}=f_i.
\]
Thus, we have
\begin{equation}\label{lemma11:eq6}
\sum_{u\in T_1}A_{uu}^2+\sum_{u\in S(i)}\sum_{v\in T_1, u\neq v}M_{uv}^2+\mathcal{O}_p(k^2/n)\leq(f_i^2+\mathcal{O}_p(k^2/n))|S(i)|.
\end{equation}
By  \eqref{lemma11:eq6}, it holds that
\begin{equation}\label{lemma11:eq5}
\begin{aligned}
&\sum_{T\neq S(i)}\Big||\A_{S(i)\times T}|h_{S(i)\times T}\Big |\\
&\leq\frac{1+\mathcal{O}_p(k^2/\sqrt{n})}{(\min_{j\notin S(i)}|f_j|)^{|S(i)|}}\prod_{j\notin S(i)}|f_j|\sum_{T\neq S(i)}(\sqrt{(f_i^2+\mathcal{O}_p(k^2/n))|S(i)|^2})^{|S(i)|}(\mathcal{O}_p(k/\sqrt{n}))^{|S(i)|}\\
&\leq\frac{1+\mathcal{O}_p(k^2/\sqrt{n})}{(\min_{j\notin S(i)}|f_j|)^{|S(i)|}}\prod_{j\notin S(i)}|f_j|\sum_{T\neq S(i)}(\sqrt{(f_i^2+\mathcal{O}_p(k^2/n))|S(i)|^2}\mathcal{O}_p(k/\sqrt{n}))^{|S(i)|}\\
&\leq\frac{1+\mathcal{O}_p(k^2/\sqrt{n})}{(\min_{j\notin S(i)}|f_j|)^{|S(i)|}}\prod_{j\notin S(i)}|f_j|k^{|S(i)|}(\sqrt{(f_i^2+\mathcal{O}_p(k^2/n))|S(i)|^2}\mathcal{O}_p(k/\sqrt{n}))^{|S(i)|}\\
&=\frac{1+\mathcal{O}_p(k^2/\sqrt{n})}{(\min_{j\notin S(i)}|f_j|)^{|S(i)|}}\prod_{j\notin S(i)}|f_j|\mathcal{O}_p(((f_i+k/\sqrt{n})k^{2}|S(i)|/n^{1/2})^{|S(i)|}).
\end{aligned}
\end{equation}
Combining \eqref{lemma11:eq1}, \eqref{lemma11:eq2}, and \eqref{lemma11:eq5}, we have
\begin{equation}\label{lemma11:eq7}
\begin{aligned}
|\A_n(d_i)_{S(i)\times S(i)}|&=\left|\frac{\sum_{T\neq S(i)}(-1)^{i_1+\dots+i_{|S(i)|}+t_1+\dots+t_{|S(i)|}}|\A_{S(i)\times T}|h_{S(i)\times T}}{-h_{S(i)\times S(i)}}\right|\\
&\leq\frac{1+\mathcal{O}_p(k^2/\sqrt{n})}{1+\mathcal{O}_p(k^2/\sqrt{n})}\frac{1}{(\min_{j\notin S(i)}|f_j|)^{|T_2|}}\mathcal{O}_p(((f_i+k/\sqrt{n})k^{2}|S(i)|/n^{1/2})^{|S(i)|})\\
&=\mathcal{O}_p(((f_i+k/\sqrt{n})k^{2}|S(i)|/n^{1/2})^{|S(i)|}).
\end{aligned}
\end{equation}
Let $\F_1=\Diag\{f_i,\dots,f_i\}$. Then, we have
\begin{equation}\label{eqn:thm22201}
\A_n(d_i)_{{S(i)\times S(i)}}=\F_1+\M_{11}.
\end{equation}
where $\M_{11}=\M_n(d_i)_{{S(i)\times S(i)}}$.
Similar to the proof of \eqref{eqn:de03}, we have
\begin{equation}\label{eqn:thm22202}
\F_1'\triangleq\Diag\{f_i+\epsilon_1,\ldots,f_i+\epsilon_{|S(i)|}\}=\V_1(d_i)^\top(\F_1+\M_{11})\V_1(d_i),
\end{equation}
where $\V_1(d_i)$ is an orthogonal matrix and $\max_{1\leq j\leq |S(i)|}|\epsilon_j|=\mathcal{O}_p(k/\sqrt{n})$.
By \eqref{eqn:thm22201} and \eqref{eqn:thm22202}, we have
\begin{equation}\label{eqn:thm22203}
|\A_n(d_i)_{{S(i)\times S(i)}}|=\prod_{j=1}^{|S(i)|}(f_i+\epsilon_j)=\prod_{j=1}^{|S(i)|}(f_i+\mathcal{O}_p(k/\sqrt{n})).
\end{equation}
By \eqref{lemma11:eq7} and \eqref{eqn:thm22203}, we have
\begin{equation}\label{eqn:thm22204}
f_i+\mathcal{O}_p(k/\sqrt{n})=\mathcal{O}_p((f_i+k/\sqrt{n})k^{2}|S(i)|/n^{1/2}).
\end{equation}
Since $k=o(n^{1/4})$, and $|S(i)|<\infty$, we have
\begin{equation}\label{eqn:thm22205}
k^{2}|S(i)|/n^{1/2}=o(1).
\end{equation}
Combining \eqref{eqn:thm22204} and \eqref{eqn:thm22205}, we have
\begin{equation}\label{eqn:thm22206}
f_i=\frac{d_i-\frac{1}{n}\tr(\I+\H_n(d_i))\lambda_i}{\lambda_i}=\mathcal{O}_p(k/\sqrt{n}).
\end{equation}
Therefore, by  \eqref{lemma11:eq7} and \eqref{eqn:thm22206}, we have
\[
|\A_n(d_i)_{S(i)\times S(i)}|=\mathcal{O}_p((\frac{k^{3}}{n})^{|S(i)|}).
\]

\eop

\section{Proofs of auxiliary lemmas}\label{sec:proofaux}

\subsection{Proof of Lemma \ref{lemma:1}}
By Lemmas \ref{lemma:spectral} and \ref{lemma:weyl}, for $1\le i\le p$, we have
\[
\vert d_i-\lambda_i\vert\leq \max_{1\leq i\leq p}\vert\phi_i\vert=||\S_n-\bSigma||=\mathcal{O}_p(\sqrt{1/n}).
\]
Thus, for $1 \leq i\leq p$, we have
\[
(d_i-\lambda_i)^2=\mathcal{O}_p(1/n).
\]
Moreover, we have
\begin{eqnarray*}(p-k)(\bar{d}_{k+1}-\lambda)^2&=&(p-k)\left\{\frac{1}{p-k}\sum_{i=k+1}^pd_i-\lambda\right\}^2\\
&=&\frac{\left\{\sum_{i=k+1}^p(d_i-\lambda)\right\}^2}{p-k}\leq p\max_{1\leq i\leq p}\vert\phi_i\vert^2=\mathcal{O}_p(1/n).
\end{eqnarray*}

\eop

\subsection{Proof of Lemma \ref{lemma:4}}

We have
\begin{eqnarray*}\log L_k&=&-\frac{n}{2}\left\{\sum_{i=1}^k\log d_i+(p-k)\log\bar{d}_{k+1}\right\}\\
&=&-\frac{n}{2}\sum_{i=1}^k\log d_i-\frac{n(p-k)}{2}\log\left[\left\{\frac{1}{p-k}\sum_{i=k+1}^pd_i-1\right\}+1\right]\\
&=&-\frac{n}{2}\sum_{i=1}^k\log d_i-\frac{n}{2}\sum_{i=k+1}^pd_i+\frac{n}{2}(p-k)+\mathcal{O}_p(1)\\
&=&-\frac{n}{2}\sum_{i=1}^k\log d_i-\frac{n}{2}\sum_{i=k+1}^p(d_i-1)+\mathcal{O}_p(1),
\end{eqnarray*}
where the third equality uses Taylor's expansion and the fact that $(p-k)(\lambda-\bar{d}_{k+1})^2=\mathcal{O}_p(1/n)$ from Lemma \ref{lemma:1}.

\eop

\subsection{Proof of Lemma \ref{lemma:large deviation1}}
Let $\L=\Diag\{l_1,\cdots,l_n\}$ and
\begin{eqnarray*}
\A=\U\L\U^\top,
\end{eqnarray*}
where $\U$ is orthogonal and its column $\U_s$ is the eigenvector of $\A$  with eigenvalue $l_s$.

Let $\T=\bLambda_{1}^{-\frac{1}{2}}\Z_1^\top\U$ and $\S=\bLambda_{1}^{-\frac{1}{2}}\Z_1^\top=\Gamma_1^\top\Y^\top=(S_{im})_{k\times n}$. Then, we have $S_{im}=\sum_{1\leq r\leq p}\Gamma_{ri}Y_{mr}$ and  $T_{is}=\sum_{1\leq m\leq n}S_{im}U_{ms}$. It holds that $\mathbb{E}(T_{is})=0$, $\mathrm{Var}(T_{is})=1$, and $\mathbb{E}(T_{is}T_{js})=0$ for $i\neq j$.
For $1\leq s\leq n$, the $s$-th column of $\T$ is
\[
\T_s=(\bLambda_{1}^{-\frac{1}{2}}\Z_1^\top\U)_s\sim \mathcal{N}_k(\0, \I),
\]
and the columns are independent and identically distributed.
Note that
\begin{eqnarray*}
\left\{\frac{1}{n}\bLambda_{1}^{-\frac{1}{2}}\Z_1^\top(\I+\A)\Z_1\bLambda_{1}^{-\frac{1}{2}}\right\}_{ij}&=&\left\{\frac{1}{n}\T(\I+\L)\T^\top\right\}_{ij}=\frac{1}{n}\sum_{s=1}^n(1+l_s)T_{is}T_{js}.
\end{eqnarray*}
Therefore, we have
$$
\mathbb{E}\left\{\frac{1}{n}\bLambda_{1}^{-\frac{1}{2}}\Z_1^\top(\I+\A)\Z_1\bLambda_{1}^{-\frac{1}{2}}\right\}=\frac{1}{n}\tr(\I+\A)\I.
$$
Moreover, for $\v\in\mathbb{R}^{k}$, we have
\[
\frac{1}{n}\v^\top\left\{\frac{1}{n}\bLambda_{1}^{-\frac{1}{2}}\Z_1^\top(\I+\A)\Z_1\bLambda_{1}^{-\frac{1}{2}}\right\}\v=\frac{1}{n}\sum_{s=1}^n(1+l_s)\left(\sum_{i,j=1}^kv_iv_jT_{is}T_{js}\right).
\]
Consequently, for $\v^\top\v=1$, it holds that
\begin{eqnarray*}
&&\frac{1}{n}\v^\top\left\{\frac{1}{n}\bLambda_{1}^{-\frac{1}{2}}\Z_1^\top(\I+\A)\Z_1\bLambda_{1}^{-\frac{1}{2}}-\sum_{s=1}^n(1+l_s)\I\right\}\v\\\nonumber
&=&\frac{1}{n}\sum_{s=1}^n(1+l_s)\left(\sum_{i,j=1}^kv_iv_jT_{is}T_{js}\right)-\frac{1}{n}\sum_{s=1}^n(1+l_s)\v^\top\I\v\\\nonumber
&=&\frac{1}{n}\sum_{s=1}^n(1+l_s)\left(\sum_{i,j=1}^kv_iv_jT_{is}T_{js}\right)-\frac{1}{n}\sum_{s=1}^n(1+l_s).
\end{eqnarray*}
Define
$$
\mathcal{S}(\T_1,\ldots,\T_n;t)=\left\{\sum_{s=1}^n(1+l_s)(\sum_{i,j=1}^kv_iv_jT_{is}T_{js})-\sum_{s=1}^n(1+l_s)>nt\right\}.
$$
Next, for any $a>0$, we have
\begin{eqnarray*}
&&\mathbb{P}\left[\frac{1}{n}\v^\top\left\{\frac{1}{n}\bLambda_{1}^{-\frac{1}{2}}\Z_1^\top(\I+\A)\Z_1\bLambda_{1}^{-\frac{1}{2}}-\tr(\I+\A)\I\right\}\v>t\right]\\\nonumber
&=&\mathbb{P}\left[\v^\top\left\{\frac{1}{n}\bLambda_{1}^{-\frac{1}{2}}\Z_1^\top(\I+\A)\Z_1\bLambda_{1}^{-\frac{1}{2}}-\sum_{s=1}^n(1+l_s)\I\right\}\v>nt\right]\\
&=&\int_{\mathcal{S}(\T_1,\ldots,\T_n;t)}\prod_{s=1}^ndF(T_{1s},\ldots,T_{ks})\\
&\leq&\int_{\mathcal{S}(\T_1,\ldots,\T_n;t)}\frac{\exp\left\{a\sum_{s=1}^n(1+l_s)\left(\sum_{i,j=1}^kv_iv_jT_{is}T_{js}\right)-a\sum_{s=1}^n(1+l_s)\right\}}{\exp(ant)}\prod_{s=1}^ndF(T_{1s},\ldots,T_{ks})\\\nonumber
&\leq&\int\frac{\exp\left\{a\sum_{s=1}^n(1+l_s)(\sum_{i,j=1}^kv_iv_jT_{is}T_{js})-a\sum_{s=1}^n(1+l_s)\right\}}{\exp(ant)}\prod_{s=1}^ndF(T_{1s},\ldots,T_{ks})\\\nonumber
&=&\exp\left\{-ant-a\sum_{s=1}^n(1+l_s)\right\}\prod_{s=1}^n\mathbb{E}\left[\exp\left\{a(1+l_s)\sum_{i,j=1}^kv_iv_jT_{is}T_{js}\right\}\right]\\
&=&\exp\left\{-ant-a\sum_{s=1}^n(1+l_s)\right\}\prod_{s=1}^n\mathbb{E}\left[\exp\left\{a(1+l_s)W_s\right\}\right],
\end{eqnarray*}
where  $F(T_{1s},\ldots,T_{ks})$ is the probability density function of $T_{1s},\ldots,T_{ks}$, and $W_s=\sum_{i,j=1}^kv_iv_jT_{is}T_{js}\sim\chi^2(1)$.
Next, we have
\begin{eqnarray*}
\mathbb{E}\left[\exp\left\{a(1+l_s)W_s\right\}\right]&=&\sum_{i=0}^\infty\int\frac{(a(1+l_s)w)^i}{i!}dF_1(w)\\
&=&1+a(1+l_s)+\sum_{i=2}^\infty\int\frac{(a(1+l_s)w)^i}{i!}dF_1(w)\\
&=&1+a(1+l_s)+\sum_{i=2}^\infty\frac{(a(1+l_s))^i\mathbb{E}(W_s^i)}{i!}\\
&=&1+a(1+l_s)+\sum_{i=2}^\infty\frac{(2a(1+l_s))^i}{i!}\frac{\Gamma(\frac{1}{2}+i)}{\Gamma(\frac{1}{2})}\\
&\leq&1+a(1+l_s)+\sum_{i=2}^\infty(2a(1+l_s))^i\\
&\leq&1+a(1+l_s)+4\rho_1a^2(1+l_s)^2,
\end{eqnarray*}
where $F_1(\cdot)$ is the probability density function of $\chi^2(1)$, $2a(1+l_1)<1$ and $\rho_1\geq\frac{1}{1-2a(1+l_1)}$.
Consequently, we have
\begin{eqnarray*}
\prod_{s=1}^n\mathbb{E}\left[\exp\left\{a(1+l_s)W_s\right\}\right]&\leq&\exp\left[\sum_{s=1}^n\log\left\{1+a(1+l_s)+4\rho_1a^2(1+l_s)^2\right\}\right]\\
&\leq&\exp\left\{a\sum_{s=1}^n(1+l_s)+4\rho_1a^2\sum_{s=1}^n(1+l_s)^2\right\}.
\end{eqnarray*}
Putting the above results together, for any $0<a<\frac{1}{2(1+l_1)}$, we have
\begin{eqnarray*}
&&\mathbb{P}\left[\frac{1}{n}\v^\top\left\{\frac{1}{n}\bLambda_{1}^{-\frac{1}{2}}\Z_1^\top(\I+\A)\Z_1\bLambda_{1}^{-\frac{1}{2}}-\tr(\I+\A)\I\right\}\v>t\right]\\\nonumber
&\leq&\exp\left\{-ant+4\rho_1a^2\sum_{s=1}^n(1+l_s)^2\right\}.
\end{eqnarray*}
Let $f(a)=-ant+4\rho_1a^2\sum_{s=1}^n(1+l_s)^2$.
Setting the derivative $f'(a)$ to zero gives $a=\frac{nt}{8\rho_1\sum_{s=1}^n(1+l_s)^2}$. Thus, we have
\begin{eqnarray*}
&&\mathbb{P}\left[\frac{1}{n}\v^\top\left\{\frac{1}{n}\bLambda_{1}^{-\frac{1}{2}}\Z_1^\top(\I+\A)\Z_1\bLambda_{1}^{-\frac{1}{2}}-\tr(\I+\A)\I\right\}\v>t\right]\\\nonumber
&\leq&\exp\left\{-\frac{nt^2}{16\rho_1\sum_{s=1}^n(1+l_s)^2/n}\right\}.
\end{eqnarray*}

Following a standard covering argument (see, for example, Lemma 3 in \cite{Cai:Zhang:Zhou:2010} or Lemma 1 in \cite{Hu:Keita:Fu:2023}), we have

\begin{eqnarray*}
&&\mathbb{P}\left\{\left\Vert\frac{1}{n}\bLambda_{1}^{-\frac{1}{2}}\Z_1^\top(\I+\A)\Z_1\bLambda_{1}^{-\frac{1}{2}}-\frac{1}{n}\tr(\I+\A)\I\right\Vert>t\right\}\\\nonumber
&\leq&2\times\mathbb{P}\left[C_1\sup_{1\leq j\leq5^k}\frac{1}{n}\v_j^\top\left\{\frac{1}{n}\bLambda_{1}^{-\frac{1}{2}}\Z_1^\top(\I+\A)\Z_1\bLambda_{1}^{-\frac{1}{2}}-\sum_{s=1}^n(1+l_s)\I\right\}\v_j>t\right]\\\nonumber
&\leq&2\sum_{j\leq5^k}\sup_{\v_j^\top\v_j=1}\mathbb{P}\left[\frac{1}{n}\v_j^\top\left\{\frac{1}{n}\bLambda_{1}^{-\frac{1}{2}}\Z_1^\top(\I+\A)\Z_1\bLambda_{1}^{-\frac{1}{2}}-\sum_{s=1}^n(1+l_s)\I\right\}\v_j>\frac{t}{C_1}\right]\\\nonumber
&\leq& 2\times5^k\exp\left\{-\frac{nt^2}{16\rho_1C_1\sum_{s=1}^n(1+l_s)^2/n}\right\}=2\exp\left\{Ck-\frac{nt^2}{\rho\sum_{s=1}^n(1+l_s)^2/n}\right\},
\end{eqnarray*}
where $\v_j^\top\v_j=1$ and $\rho=16\rho_1C_1$ with $C_1$ being a constant depending on $C$.
Noting that $2a(1+l_1)<1$, we have
\[
t<4\rho_1\frac{\sum_{s=1}^n(1+l_s)^2/n}{1+l_1}=\rho\frac{\sum_{s=1}^n(1+l_s)^2/n}{4C_1(1+l_1)}.
\]

\eop

\subsection{Proof of Lemma \ref{corollary:large deviation2}}
The proof is similar to that of Lemma \ref{lemma:large deviation1}.

By Markov's inequality, it holds that
\begin{equation}
\label{lemma:large deviation2-equation-1}
\begin{aligned}
&\mathbb{P}\left\{\left\Vert\frac{1}{n}\bLambda_{1}^{-\frac{1}{2}}\Z_1^\top(\I+\A)\Z_1\bLambda_{1}^{-\frac{1}{2}}-\frac{1}{n}\tr(\I+\A)\I\right\Vert>t\right\}\\
\leq&\mathbb{P}\left\{\tr\left[\frac{1}{n}\bLambda_{1}^{-\frac{1}{2}}\Z_1^\top(\I+\A)\Z_1\bLambda_{1}^{-\frac{1}{2}}-\frac{1}{n}\tr(\I+\A)\I\right]^2>t^2\right\}\\
\leq&\mathbb{P}\left\{\sum_{1\leq i,j\leq k}\left[\frac{1}{n}\sum_{s=1}^n(1+l_s)T_{is}T_{js}-\frac{1}{n}\sum_{s=1}^n(1+l_s)I_{ij}\right]^2>t^2\right\}\\
\leq&\mathbb{P}\left\{\sum_{1\leq i,j\leq k}\left[\sum_{s=1}^n(1+l_s)(T_{is}T_{js}-I_{ij})\right]^2>n^2t^2\right\}\\
\leq&\frac{\sum_{1\leq i,j\leq k}\mathbb{E}\left[\sum_{s=1}^n(1+l_s)(T_{is}T_{js}-I_{ij})\right]^2}{n^2t^2}\\
\leq&\frac{\sum_{1\leq i,j\leq k}\mathbb{E}Q_{ij}^2}{n^2t^2}.
\end{aligned}
\end{equation}
where
$Q_{ij}=\sum_{s=1}^n(1+l_s)(T_{is}T_{js}-I_{ij}).$

Let $R_s=T_{is}T_{js}-I_{ij}$. Then
\begin{equation}
\label{lemma:large deviation2-equation-2}
\begin{aligned}
\mathbb{E}Q_{ij}^2=&\sum_{1\leq s,t\leq n}(1+l_s)(1+l_t)\mathbb{E}(R_sR_t).
\end{aligned}
\end{equation}
Recall that $S_{im}=\sum_{1\leq r\leq p}\Gamma_{ri}Y_{mr}$ and  $T_{is}=\sum_{1\leq m\leq n}S_{im}U_{ms}$. Note that for $m\neq m'$, $S_{im}$ and $S_{jm'}$ are independent and for $i\neq j$, $S_{im}$ and $S_{jm}$ are uncorrelated. It holds that
\begin{equation}
\label{lemma:large deviation2-equation-0}
\begin{aligned}
\mathbb{E}(T_{is}T_{jt})=&\sum_{m_1,m_2}U_{m_1s}U_{m_2t}\mathbb{E}(S_{im_1}S_{jm_2})=\sum_{m}^nU_{ms}U_{mt}\mathbb{E}(S_{im}S_{jm})\\
=&\sum_{m}^nU_{ms}U_{mt}\sum_{r_1,r_2}\Gamma_{r_1i}\Gamma_{r_2j}\mathbb{E}(Y_{mr_1}Y_{mr_2})\\
=&\sum_{m}^nU_{ms}U_{mt}\sum_{r}^p\Gamma_{ri}\Gamma_{rj}=\sum_{r}^p\Gamma_{ri}\Gamma_{rj}\sum_{m}^nU_{ms}U_{mt}\\
=&I_{ij}I_{st}.
\end{aligned}
\end{equation}
By \eqref{lemma:large deviation2-equation-0}, we have
\begin{equation}
\label{lemma:large deviation2-equation-3}
\begin{aligned}
\mathbb{E}(R_sR_t)=&\mathbb{E}\left[(T_{is}T_{js}-I_{ij})(T_{it}T_{jt}-I_{ij})\right]\\
=&\mathbb{E}(T_{is}T_{js}T_{it}T_{jt})-I_{ij}\mathbb{E}(T_{is}T_{js})-I_{ij}\mathbb{E}(T_{it}T_{jt})+I_{ij}^2\\
=&\mathbb{E}(T_{is}T_{js}T_{it}T_{jt})-I_{ij}.
\end{aligned}
\end{equation}
By \eqref{lemma:large deviation2-equation-2} and \eqref{lemma:large deviation2-equation-3}, we have
\begin{equation}
\label{lemma:large deviation2-equation-4}
\begin{aligned}
\mathbb{E}Q_{ij}^2=&\sum_{1\leq s,t\leq n}(1+l_s)(1+l_t)\left[\mathbb{E}(T_{is}T_{js}T_{it}T_{jt})-I_{ij}\right].
\end{aligned}
\end{equation}

It holds that
\begin{equation}
\label{lemma:large deviation2-equation-5-00}
\begin{aligned}
\mathbb{E}(S_{im}S_{jm}S_{um}S_{vm})=&\sum_{r_1,r_2,r_3,r_4}\Gamma_{r_1i}\Gamma_{r_2j}\Gamma_{r_3u}\Gamma_{r_4v}\mathbb{E}(Y_{mr_1}Y_{mr_2}Y_{mr_3}Y_{mr_4}).
\end{aligned}
\end{equation}
For $r_1=r_2=r_3=r_4=r$, it holds that
\begin{equation}
\label{lemma:large deviation2-equation-5-01}
\begin{aligned}
&\sum_{r=1}^p\Gamma_{r_1i}\Gamma_{r_2j}\Gamma_{r_3u}\Gamma_{r_4v}\mathbb{E}(Y_{mr_1}Y_{mr_2}Y_{mr_3}Y_{mr_4})\\
=&\eta\sum_{r=1}^p\Gamma_{ri}\Gamma_{rj}\Gamma_{ru}\Gamma_{rv}\\
=&3\sum_{r=1}^p\Gamma_{ri}\Gamma_{rj}\Gamma_{ru}\Gamma_{rv}+(\eta-3)\sum_{r=1}^p\Gamma_{ri}\Gamma_{rj}\Gamma_{ru}\Gamma_{rv}.
\end{aligned}
\end{equation}
For $r_1=r_2=r$, $r_3=r_4=r'$, $r\neq r'$, it holds that
\begin{equation}
\label{lemma:large deviation2-equation-5-02}
\begin{aligned}
&\sum_{r\neq r'}\Gamma_{r_1i}\Gamma_{r_2j}\Gamma_{r_3u}\Gamma_{r_4v}\mathbb{E}(Y_{mr_1}Y_{mr_2}Y_{mr_3}Y_{mr_4})\\
=&\sum_{r\neq r'}\Gamma_{ri}\Gamma_{rj}\Gamma_{r'u}\Gamma_{r'v}.
\end{aligned}
\end{equation}
Similarly, for $r_1=r_3=r$, $r_2=r_4=r'$, $r\neq r'$, it holds that
\begin{equation}
\label{lemma:large deviation2-equation-5-03}
\begin{aligned}
&\sum_{r\neq r'}\Gamma_{r_1i}\Gamma_{r_2j}\Gamma_{r_3u}\Gamma_{r_4v}\mathbb{E}(Y_{mr_1}Y_{mr_2}Y_{mr_3}Y_{mr_4})\\
=&\sum_{r\neq r'}\Gamma_{ri}\Gamma_{ru}\Gamma_{r'j}\Gamma_{r'v}.
\end{aligned}
\end{equation}
For $r_1=r_4=r$, $r_2=r_3=r'$, $r\neq r'$, it holds that
\begin{equation}
\label{lemma:large deviation2-equation-5-04}
\begin{aligned}
&\sum_{r\neq r'}\Gamma_{r_1i}\Gamma_{r_2j}\Gamma_{r_3u}\Gamma_{r_4v}\mathbb{E}(Y_{mr_1}Y_{mr_2}Y_{mr_3}Y_{mr_4})\\
=&\sum_{r\neq r'}\Gamma_{ri}\Gamma_{rv}\Gamma_{r'j}\Gamma_{r'u}.
\end{aligned}
\end{equation}
It holds that
\begin{equation}
\label{lemma:large deviation2-equation-5-05}
\begin{aligned}
&\sum_{r=1}^p\Gamma_{ri}\Gamma_{rj}\Gamma_{ru}\Gamma_{rv}+\sum_{r\neq r'}\Gamma_{ri}\Gamma_{rj}\Gamma_{r'u}\Gamma_{r'v}\\
=&(\sum_{r=1}^p\Gamma_{ri}\Gamma_{rj})(\sum_{r'=1}^p\Gamma_{r'u}\Gamma_{r'v})=I_{ij}I_{uv}.
\end{aligned}
\end{equation}
Similarly, it holds that
\begin{equation}
\label{lemma:large deviation2-equation-5-06}
\begin{aligned}
\sum_{r=1}^p\Gamma_{ri}\Gamma_{rj}\Gamma_{ru}\Gamma_{rv}+\sum_{r\neq r'}\Gamma_{ri}\Gamma_{ru}\Gamma_{r'j}\Gamma_{r'v}=&I_{iu}I_{jv}.
\end{aligned}
\end{equation}
It also holds that
\begin{equation}
\label{lemma:large deviation2-equation-5-07}
\begin{aligned}
\sum_{r=1}^p\Gamma_{ri}\Gamma_{rj}\Gamma_{ru}\Gamma_{rv}+\sum_{r\neq r'}\Gamma_{ri}\Gamma_{rv}\Gamma_{r'j}\Gamma_{r'u}=&I_{iv}I_{ju}.
\end{aligned}
\end{equation}
Combining \eqref{lemma:large deviation2-equation-5-00}-\eqref{lemma:large deviation2-equation-5-07}, we have
\begin{equation}
\label{lemma:large deviation2-equation-5}
\begin{aligned}
\mathbb{E}(S_{im}S_{jm}S_{um}S_{vm})=&\sum_{r_1,r_2,r_3,r_4}\Gamma_{r_1i}\Gamma_{r_2j}\Gamma_{r_3u}\Gamma_{r_4v}\mathbb{E}(Y_{mr_1}Y_{mr_2}Y_{mr_3}Y_{mr_4}).\\
=&I_{ij}I_{uv}+I_{iu}I_{jv}+I_{iv}I_{ju}+(\eta-3)\sum_{r=1}^p\Gamma_{ri}\Gamma_{rj}\Gamma_{ru}\Gamma_{rv}.
\end{aligned}
\end{equation}

For $s=t$, it holds that
\begin{equation}
\label{lemma:large deviation2-equation-6}
\begin{aligned}
\mathbb{E}(T_{is}T_{js}T_{us}T_{vs})=&\sum_{m_1,m_2,m_3,m_4}U_{m_1s}U_{m_2s}U_{m_3s}U_{m_4s}\mathbb{E}(S_{im_1}S_{jm_2}S_{um_3}S_{vm_4}).
\end{aligned}
\end{equation}

For $m_1=m_2=m_3=m_4=m$,  by \eqref{lemma:large deviation2-equation-5}, we have
\begin{equation}
\label{lemma:large deviation2-equation-7}
\begin{aligned}
&\sum_{m=1}^nU_{m_1s}U_{m_2s}U_{m_3s}U_{m_4s}\mathbb{E}(S_{im_1}S_{jm_2}S_{um_3}S_{vm_4})\\
=&\sum_{m=1}^nU_{ms}^4\left[I_{ij}I_{uv}+I_{iu}I_{jv}+I_{iv}I_{ju}+(\eta-3)\sum_{r=1}^p\Gamma_{ri}\Gamma_{rj}\Gamma_{ru}\Gamma_{rv}\right].
\end{aligned}
\end{equation}
For $m_1=m_2=m$, $m_3=m_4=m', m'\neq m$, it holds that
\begin{equation}
\label{lemma:large deviation2-equation-8}
\begin{aligned}
&\sum_{m\neq m'} U_{m_1s}U_{m_2s}U_{m_3s}U_{m_4s}\mathbb{E}(S_{im_1}S_{jm_2}S_{um_3}S_{vm_4})\\
=&\sum_{m\neq m'} U_{m_1s}U_{m_2s}U_{m_3s}U_{m_4s}\mathbb{E}(S_{im_1}S_{jm_2})\mathbb{E}(S_{um_3}S_{vm_4})\\
=&\sum_{m\neq m'} U_{ms}^2U_{m's}^2I_{ij}I_{uv}\\
=&\left[(\sum_{m} U_{ms}^2)^2-\sum_{m} U_{ms}^4\right]I_{ij}I_{uv}\\
=&\left[1-\sum_{m} U_{ms}^4\right]I_{ij}I_{uv}.
\end{aligned}
\end{equation}
Similarly, for $m_1=m_3=m$, $m_2=m_4=m'$, $m'\neq m$, it holds that
\begin{equation}
\label{lemma:large deviation2-equation-9}
\begin{aligned}
\sum_{m\neq m'} U_{m_1s}U_{m_2s}U_{m_3s}U_{m_4s}\mathbb{E}(S_{im_1}S_{jm_2}S_{um_3}S_{vm_4})=&\left[1-\sum_{m} U_{ms}^4\right]I_{iu}I_{jv}.
\end{aligned}
\end{equation}
For $m_1=m_4=m$, $m_2=m_3=m'$, $m'\neq m$, it also holds that
\begin{equation}
\label{lemma:large deviation2-equation-10}
\begin{aligned}
\sum_{m\neq m'} U_{m_1s}U_{m_2s}U_{m_3s}U_{m_4s}\mathbb{E}(S_{im_1}S_{jm_2}S_{um_3}S_{vm_4})&=&\left[1-\sum_{m} U_{ms}^4\right]I_{iv}I_{ju}.
\end{aligned}
\end{equation}
Thus, for $s=t$, combining \eqref{lemma:large deviation2-equation-6}-\eqref{lemma:large deviation2-equation-10}, we have
\begin{equation}
\label{lemma:large deviation2-equation-11}
\begin{aligned}
\mathbb{E}(T_{is}T_{js}T_{us}T_{vs})=&I_{ij}I_{uv}+I_{iu}I_{jv}+I_{iv}I_{ju}+(\eta-3)\sum_{m=1}^nU_{ms}^4\sum_{r=1}^p\Gamma_{ri}\Gamma_{rj}\Gamma_{ru}\Gamma_{rv}.
\end{aligned}
\end{equation}
By \eqref{lemma:large deviation2-equation-11}, we have
\begin{equation}
\label{lemma:large deviation2-equation-12}
\begin{aligned}
\mathbb{E}(T_{is}^2T_{js}^2)=&1+2I_{ij}+(\eta-3)\sum_{m=1}^nU_{ms}^4\sum_{r=1}^p\Gamma_{ri}^2\Gamma_{rj}^2.
\end{aligned}
\end{equation}

For $s\neq t$, it holds that
\begin{equation}
\label{lemma:large deviation2-equation-13}
\begin{aligned}
\mathbb{E}(T_{is}T_{js}T_{ut}T_{vt})=&\sum_{m_1,m_2,m_3,m_4}U_{m_1s}U_{m_2s}U_{m_3t}U_{m_4t}\mathbb{E}(S_{im_1}S_{jm_2}S_{um_3}S_{vm_4}).
\end{aligned}
\end{equation}
For $m_1=m_2=m_3=m_4=m$, by \eqref{lemma:large deviation2-equation-13}, we have
\begin{equation}
\label{lemma:large deviation2-equation-14}
\begin{aligned}
&\sum_{m=1}^nU_{m_1s}U_{m_2s}U_{m_3t}U_{m_4t}\mathbb{E}(S_{im_1}S_{jm_2}S_{um_3}S_{vm_4})\\
=&\sum_{m=1}^nU_{ms}^2U_{mt}^2\left[I_{ij}I_{uv}+I_{iu}I_{jv}+I_{iv}I_{ju}+(\eta-3)\sum_{r=1}^p\Gamma_{ri}\Gamma_{rj}\Gamma_{ru}\Gamma_{rv}\right].
\end{aligned}
\end{equation}
For $m_1=m_2=m$, $m_3=m_4=m'$, $m'\neq m$,  by \eqref{lemma:large deviation2-equation-13}, we have
\begin{equation}
\label{lemma:large deviation2-equation-15}
\begin{aligned}
&\sum_{m\neq m'} U_{m_1s}U_{m_2s}U_{m_3t}U_{m_4t}\mathbb{E}(S_{im_1}S_{jm_2}S_{um_3}S_{vm_4})\\
=&\sum_{m\neq m'} U_{ms}^2U_{m't}^2I_{ij}I_{uv}\\
=&\left[\sum_{m} U_{ms}^2\sum_{m'} U_{m't}^2-\sum_{m} U_{ms}^2U_{mt}^2\right]I_{ij}I_{uv}\\
=&\left[1-\sum_{m} U_{ms}^2U_{mt}^2\right]I_{ij}I_{uv}.
\end{aligned}
\end{equation}
For $m_1=m_3=m$, $m_2=m_4=m'$, $m'\neq m$,  by \eqref{lemma:large deviation2-equation-13}, we have
\begin{equation}
\label{lemma:large deviation2-equation-16}
\begin{aligned}
&\sum_{m\neq m'} U_{m_1s}U_{m_2s}U_{m_3t}U_{m_4t}\mathbb{E}(S_{im_1}S_{jm_2}S_{um_3}S_{vm_4})\\
=&\sum_{m\neq m'} U_{ms}U_{mt}U_{m's}U_{m't}I_{iu}I_{jv}\\
=&\left[(\sum_{m} U_{ms}U_{mt})^2-\sum_{m} U_{ms}^2U_{mt}^2\right]I_{iu}I_{jv}\\
=&-\sum_{m} U_{ms}^2U_{mt}^2I_{iu}I_{jv}.
\end{aligned}
\end{equation}
Similarly, for $m_1=m_4=m$, $m_2=m_3=m'$, $m'\neq m$, it holds that
\begin{equation}
\label{lemma:large deviation2-equation-17}
\begin{aligned}
&\sum_{m\neq m'} U_{m_1s}U_{m_2s}U_{m_3t}U_{m_4t}\mathbb{E}(S_{im_1}S_{jm_2}S_{um_3}S_{vm_4})\\
=&-\sum_{m} U_{ms}^2U_{mt}^2I_{iv}I_{ju}.
\end{aligned}
\end{equation}
Thus, for $s\neq t$, combining \eqref{lemma:large deviation2-equation-13}-\eqref{lemma:large deviation2-equation-17}, we have
\begin{equation}
\label{lemma:large deviation2-equation-18}
\begin{aligned}
\mathbb{E}(T_{is}T_{js}T_{ut}T_{vt})=&I_{ij}I_{uv}+(\eta-3)\sum_{m=1}^nU_{ms}^2U_{mt}^2\sum_{r=1}^p\Gamma_{ri}\Gamma_{rj}\Gamma_{ru}\Gamma_{rv}.
\end{aligned}
\end{equation}
By \eqref{lemma:large deviation2-equation-18}, we have
\begin{equation}
\label{lemma:large deviation2-equation-19}
\begin{aligned}
\mathbb{E}(T_{is}T_{js}T_{it}T_{jt})=&I_{ij}+(\eta-3)\sum_{m=1}^nU_{ms}^2U_{mt}^2\sum_{r=1}^p\Gamma_{ri}^2\Gamma_{rj}^2.
\end{aligned}
\end{equation}

For  $i=j$,  by \eqref{lemma:large deviation2-equation-4}, \eqref{lemma:large deviation2-equation-12} and \eqref{lemma:large deviation2-equation-19}, we have
\begin{equation*}
\begin{aligned}
\mathbb{E}Q_{ii}^2=&\sum_{s=1}^n(1+l_s)^2\left[2+(\eta-3)\sum_{m=1}^nU_{ms}^4\sum_{r=1}^p\Gamma_{ri}^4\right]\\
&+\sum_{s\neq t}^n(1+l_s)(1+l_t)(\eta-3)\sum_{m=1}^nU_{ms}^2U_{mt}^2\sum_{r=1}^p\Gamma_{ri}^4\\
=&2\sum_{s=1}^n(1+l_s)^2\\
&+(\eta-3)\sum_{r=1}^p\Gamma_{ri}^4\sum_{m=1}^n\left[\sum_{s=1}^n(1+l_s)^2U_{ms}^4+\sum_{s\neq t}^n(1+l_s)(1+l_t)U_{ms}^2U_{mt}^2\right]\\
=&2\sum_{s=1}^n(1+l_s)^2+(\eta-3)\sum_{r=1}^p\Gamma_{ri}^4\sum_{m=1}^n\left[\sum_{s=1}^n(1+l_s)U_{ms}^2\right]^2.
\end{aligned}
\end{equation*}
That is,
\begin{equation}
\label{lemma:large deviation2-equation-20}
\begin{aligned}
\mathbb{E}Q_{ii}^2=&2\sum_{s=1}^n(1+l_s)^2+(\eta-3)\sum_{r=1}^p\Gamma_{ri}^4\sum_{m=1}^n\left[\sum_{s=1}^n(1+l_s)U_{ms}^2\right]^2\\
\leq&2\sum_{s=1}^n(1+l_s)^2+(\eta-3)\sum_{m=1}^n\left[\sum_{s=1}^n\{(1+l_s)U_{ms}\}U_{ms}\right]^2\\
\leq&2\sum_{s=1}^n(1+l_s)^2+(\eta-3)\sum_{m=1}^n\sum_{s=1}^n(1+l_s)^2U_{ms}^2\sum_{s=1}^nU_{ms}^2\\
=&2\sum_{s=1}^n(1+l_s)^2+(\eta-3)\sum_{m=1}^n\sum_{s=1}^n(1+l_s)^2U_{ms}^2\\
=&2\sum_{s=1}^n(1+l_s)^2+(\eta-3)\sum_{s=1}^n(1+l_s)^2\sum_{m=1}^nU_{ms}^2\\
=&2\sum_{s=1}^n(1+l_s)^2+|\eta-3|\sum_{s=1}^n(1+l_s)^2\\
=&(2+|\eta-3|)\sum_{s=1}^n(1+l_s)^2.
\end{aligned}
\end{equation}

Similarly, for $i\neq j$, we have
\begin{equation}
\label{lemma:large deviation2-equation-21}
\begin{aligned}
\mathbb{E}Q_{ij}^2=&\sum_{s=1}^n(1+l_s)^2\left[1+(\eta-3)\sum_{m=1}^nU_{ms}^4\sum_{r=1}^p\Gamma_{ri}^2\Gamma_{rj}^2\right]\\
&+\sum_{s\neq t}^n(1+l_s)(1+l_t)(\eta-3)\sum_{m=1}^nU_{ms}^2U_{mt}^2\sum_{r=1}^p\Gamma_{ri}^2\Gamma_{rj}^2\\
=&\sum_{s=1}^n(1+l_s)^2\\
&+(\eta-3)\sum_{r=1}^p\Gamma_{ri}^2\Gamma_{rj}^2\sum_{m=1}^n\left[\sum_{s=1}^n(1+l_s)^2U_{ms}^4+\sum_{s\neq t}^n(1+l_s)(1+l_t)U_{ms}^2U_{mt}^2\right]\\
=&\sum_{s=1}^n(1+l_s)^2+(\eta-3)\sum_{r=1}^p\Gamma_{ri}^2\Gamma_{rj}^2\sum_{m=1}^n\left[\sum_{s=1}^n(1+l_s)U_{ms}^2\right]^2\\
\leq&\sum_{s=1}^n(1+l_s)^2+|\eta-3|\sum_{s=1}^n(1+l_s)^2\\
=&(1+|\eta-3|)\sum_{s=1}^n(1+l_s)^2.
\end{aligned}
\end{equation}
Combining \eqref{lemma:large deviation2-equation-1}, \eqref{lemma:large deviation2-equation-20} and \eqref{lemma:large deviation2-equation-21}, we have
\begin{equation}
\label{lemma:large deviation2-equation-22}
\begin{aligned}
&\mathbb{P}\left\{\left\Vert\frac{1}{n}\Gamma_1^\top\Y^\top(\I+\A)\Y\Gamma_1-\frac{1}{n}\tr(\I+\A)\I\right\Vert>t\right\}\\
\leq&\frac{\sum_{1\leq i,j\leq k}\mathbb{E}Q_{ij}^2}{n^2t^2}\\
\leq&\frac{k^2(\eta+5)\sum_{s=1}^n(1+l_s)^2/n}{nt^2}.
\end{aligned}
\end{equation}
We also have
\begin{equation}
\label{lemma:large deviation2-equation-23}
\begin{aligned}
&\mathbb{P}\left\{\tr\left[\frac{1}{n}\bLambda_{1}^{-\frac{1}{2}}\Z_1^\top(\I+\A)\Z_1\bLambda_{1}^{-\frac{1}{2}}-\frac{1}{n}\tr(\I+\A)\I\right]^2>t'\right\}\\
\leq&\frac{\sum_{1\leq i,j\leq k}\mathbb{E}Q_{ij}^2}{n^2t'}\\
\leq&\frac{k^2(\eta+5)\sum_{s=1}^n(1+l_s)^2/n}{nt'}.
\end{aligned}
\end{equation}
Therefore, by \eqref{lemma:large deviation2-equation-22} and \eqref{lemma:large deviation2-equation-23}, we have
\[
\left\Vert\frac{1}{n}\bLambda_{1}^{-\frac{1}{2}}\Z_1^\top(\I+\A)\Z_1\bLambda_{1}^{-\frac{1}{2}}-\frac{1}{n}\tr(\I+\A)\I\right\Vert=\sqrt{\frac{1}{n}\sum_{s=1}^n(1+l_s)^2}\mathcal{O}_p(k/\sqrt{n}).
\]
We also have
\[
\tr\left[\frac{1}{n}\bLambda_{1}^{-\frac{1}{2}}\Z_1^\top(\I+\A)\Z_1\bLambda_{1}^{-\frac{1}{2}}-\frac{1}{n}\tr(\I+\A)\I\right]^2=\frac{1}{n}\sum_{s=1}^n(1+l_s)^2\mathcal{O}_p(k^2/n).
\]
\eop

\section{Additional Simulation Results} \label{sec:addsim}
\textbf{Simulation 5 (the standard spiked covariance model, large $p$, comparison with BCF, ACT).}
Under the same model as in Simulation 2, we consider the case of large $p$ and varying $\rho$, with $\rho=3\sqrt{p/n}$ and $k=10$. We compare GIC with BCF, ACT and AGIC.

It is seen from Table \ref{tab:5} that, both BCF and ACT can be sensitive to the values of $p$ and $n$. In comparison,  except for a few settings,  GIC and AGIC perform better than BCF and ACT, especially when $p$ or $n$ is relatively small (i.e., $p\leq 300$ or $n\leq 300$).  Even for these few settings, the performances of GIC and AGIC are also similar to those of BCF and ACT. We also note that, for the standard spiked covariance model, GIC outperforms AGIC.

\begin{table}[!t]
\setlength{\tabcolsep}{3pt}
\centering
{\renewcommand{\arraystretch}{0.85}
\begin{tabular}{ccc|cc|cc|cc|ccccccccc}
\hline
& \multicolumn{ 2}{c|}{$n=100$}&\multicolumn{ 2}{c|}{$n=200$}&\multicolumn{ 2}{c|}{$n=300$}&\multicolumn{ 2}{c|}{$n=400$}&\multicolumn{ 2}{c}{$n=500$}\\
\hline
GIC &                              $\widehat{\mathbb{P}}(\hat{k}=k)$ & $\widehat{\mathbb{E}}(\hat{k})$ &      $\widehat{\mathbb{P}}(\hat{k}=k)$ & $\widehat{\mathbb{E}}(\hat{k})$ & $\widehat{\mathbb{P}}(\hat{k}=k)$ & $\widehat{\mathbb{E}}(\hat{k})$  &    $\widehat{\mathbb{P}}(\hat{k}=k)$ & $\widehat{\mathbb{E}}(\hat{k})$ &   $\widehat{\mathbb{P}}(\hat{k}=k)$ & $\widehat{\mathbb{E}}(\hat{k})$ &\\
\hline
$p=100$                       &\textbf{0.39}  & \textbf{9.40} & \textbf{0.73}& \textbf{9.77}&\textbf{0.85}  & \textbf{9.89}& \textbf{0.93}& \textbf{9.94}&\textbf{0.96}&\textbf{9.96}&\\
$p=200$                       &\textbf{0.56}  & \textbf{9.89} & \textbf{0.79}& \textbf{9.79}&\textbf{0.92}  & \textbf{9.92}&  \textbf{0.99}& \textbf{9.99}&  \textbf{0.99} &\textbf{9.99}&\\
$p=300$                       &\textbf{0.63}  & \textbf{9.83} & \textbf{0.96}& \textbf{9.96}& 0.92 & 9.92 & \textbf{1.00}& \textbf{10.0}&  \textbf{1.00}  & \textbf{10.0}&\\
$p=400$                       &\textbf{0.58}   & \textbf{9.90} &     \textbf{0.93}& \textbf{9.95}&   0.96  & 9.96&     \textbf{1.00}& \textbf{10.0}&     \textbf{1.00}  & \textbf{10.0}&\\
$p=500$                       &\textbf{0.65}   & \textbf{9.87} &     \textbf{0.91}& \textbf{9.91}&   0.94  & 9.94&     \textbf{1.00}& \textbf{10.0}&     \textbf{1.00}  & \textbf{10.0}&\\
\hline
& \multicolumn{ 2}{c|}{$n=100$}&\multicolumn{ 2}{c|}{$n=200$}&\multicolumn{ 2}{c|}{$n=300$}&\multicolumn{ 2}{c|}{$n=400$}&\multicolumn{ 2}{c}{$n=500$}\\
\hline
BCF &                              $\widehat{\mathbb{P}}(\hat{k}=k)$ & $\widehat{\mathbb{E}}(\hat{k})$ &      $\widehat{\mathbb{P}}(\hat{k}=k)$ & $\widehat{\mathbb{E}}(\hat{k})$ & $\widehat{\mathbb{P}}(\hat{k}=k)$ & $\widehat{\mathbb{E}}(\hat{k})$  &    $\widehat{\mathbb{P}}(\hat{k}=k)$ & $\widehat{\mathbb{E}}(\hat{k})$ &   $\widehat{\mathbb{P}}(\hat{k}=k)$ & $\widehat{\mathbb{E}}(\hat{k})$ &\\
\hline
$p=100$                       &0.06  & 7.64 &     0.48& 9.35&   0.73  & 9.73&     0.88& 9.88&     0.94  & 9.97\\
$p=200$                       &0.11  & 8.39 &     0.16& 8.77&   0.64  & 9.61&     0.89& 9.89&     0.97  & 9.98\\
$p=300$                       &0.16  & 8.67 &     0.36& 9.23&   0.43  & 9.32&     0.82& 9.82&     0.94  & 9.94&\\
$p=400$                       &0.21   & 8.86 &     0.57& 9.47&   0.68  & 9.65&    0.61& 9.59&     0.91  & 9.91&\\
$p=500$                       &0.32   & 9.07 &     0.64& 9.622&   0.77  & 9.78&     0.85& 9.84&    0.82  & 9.82&\\
\hline
& \multicolumn{ 2}{c|}{$n=100$}&\multicolumn{ 2}{c|}{$n=200$}&\multicolumn{ 2}{c|}{$n=300$}&\multicolumn{ 2}{c|}{$n=400$}&\multicolumn{ 2}{c}{$n=500$}\\
\hline
ACT  &                           $\widehat{\mathbb{P}}(\hat{k}=k)$ & $\widehat{\mathbb{E}}(\hat{k})$ &      $\widehat{\mathbb{P}}(\hat{k}=k)$ & $\widehat{\mathbb{E}}(\hat{k})$ & $\widehat{\mathbb{P}}(\hat{k}=k)$ & $\widehat{\mathbb{E}}(\hat{k})$  &    $\widehat{\mathbb{P}}(\hat{k}=k)$ & $\widehat{\mathbb{E}}(\hat{k})$ &   $\widehat{\mathbb{P}}(\hat{k}=k)$ & $\widehat{\mathbb{E}}(\hat{k})$ &\\
\hline
$p=100$                       &0.00  & 4.62 &0.02& 6.22&0.16  & 7.76&  0.35& 8.62&0.54  & 9.14\\
$p=200$                       &0.02  & 6.60 & 0.35& 8.85&   0.78  & 9.72&     0.92& 9.91&     0.97  & 9.97\\
$p=300$                       &0.05  & 7.50 & 0.65& 9.54&   0.92  & 9.91&     0.99& 9.99&     \textbf{1.00}  & \textbf{10.0}&\\
$p=400$                       &0.13   & 8.01 &     0.80& 9.76&   0.97  & 9.97&    \textbf{1.00}& \textbf{10.0}&     \textbf{1.00}  & \textbf{10.0}&\\
$p=500$                       &0.13  & 8.24 &     0.82& 9.80&   0.97  & 9.99&    \textbf{1.00}  & \textbf{10.0}&  \textbf{1.00}  & \textbf{10.0}&\\
\hline
& \multicolumn{ 2}{c|}{$n=100$}&\multicolumn{ 2}{c|}{$n=200$}&\multicolumn{ 2}{c|}{$n=300$}&\multicolumn{ 2}{c|}{$n=400$}&\multicolumn{ 2}{c}{$n=500$}\\
\hline
AGIC  &                           $\widehat{\mathbb{P}}(\hat{k}=k)$ & $\widehat{\mathbb{E}}(\hat{k})$ &      $\widehat{\mathbb{P}}(\hat{k}=k)$ & $\widehat{\mathbb{E}}(\hat{k})$ & $\widehat{\mathbb{P}}(\hat{k}=k)$ & $\widehat{\mathbb{E}}(\hat{k})$  &    $\widehat{\mathbb{P}}(\hat{k}=k)$ & $\widehat{\mathbb{E}}(\hat{k})$ &   $\widehat{\mathbb{P}}(\hat{k}=k)$ & $\widehat{\mathbb{E}}(\hat{k})$ &\\
\hline
$p=100$                       &0.00  & 6.67 &0.09& 8.24&0.32  & 9.00&  0.52& 9.44&0.74  & 9.69\\
$p=200$                       &0.08  & 8.21 & 0.54& 9.45&   0.91  & 9.91&     0.98& 9.98&     0.99  & 9.99\\
$p=300$                       &0.18  & 8.78 & 0.79& 9.77&   \textbf{0.96}  & \textbf{9.97}&     0.99& 10.0&     \textbf{1.00}& \textbf{10.0}&\\
$p=400$                       &0.31   & 9.06 &     0.89& 9.89&   \textbf{0.98}  & \textbf{10.0}&    \textbf{1.00}& \textbf{10.0}&     \textbf{1.00}  & \textbf{10.0}&\\
$p=500$                       &0.34   & 9.15 &     0.90& 9.92&   \textbf{0.98} & \textbf{10.0}&    \textbf{1.00}  & \textbf{10.0}&  \textbf{1.00}  & \textbf{10.0}&\\
\hline
\end{tabular}}
\caption{Performances of GIC, BCF, ACT and AGIC  in Simulation 5 with $k=10$. Marked in boldface are those achieving the best evaluation criteria in each setting.}
\label{tab:5}
\end{table}

\medskip

\textbf{Simulation 6  (the general spiked covariance model from \cite{Fan:Guo:Zheng:2022}, large $p$, comparison with ACT).}
In this simulation, we set $\x\sim\mathcal{N}_p(\0,\bSigma)$, where $\bSigma=\Q\Q^\top+\Diag\{\nu_1^2,\ldots, \nu_p^2\}$. The matrix $\Q=(q_{ij})_{pk}$ is obtained as follows. For $1\leq j\leq 2$, let $q_{ij}$ be i.i.d. from $\textrm{N}(0,25)$ for $i\in \{1,\ldots, \rho p\}$ and $q_{ij}$ be i.i.d. from Uniform(0,0.05) for $i \in \{\rho p+1,\ldots, p\}$. For $3\leq j\leq 5$, let $q_{ij}$ be i.i.d. from $\textrm{N}(0,1)$ for $i \in \{1, \ldots, p\}$.  Let $\nu_1^2,\ldots, \nu_p^2$ be i.i.d. from Uniform(0, 180). We set $\rho=0.05$, $k=5$ and restrict the candidate spike size in the range of $k'\in\{0,1,\ldots,20\}$. In this simulation, for example, when $n=500$, $p=200$, after generating a $\bSigma$, we calculate its eigenvalues,  $\lambda_k=292.9$, $\lambda_{k+1}=178.9$, ..., $\lambda_{p-1}=0.7431$, $\lambda_{p}=0.6372$. Thus, $\text{SNR}=\lambda_k/\lambda_{k+1}-1=0.6372$.

\begin{figure}[!t]
\centering
\includegraphics[width=\textwidth]{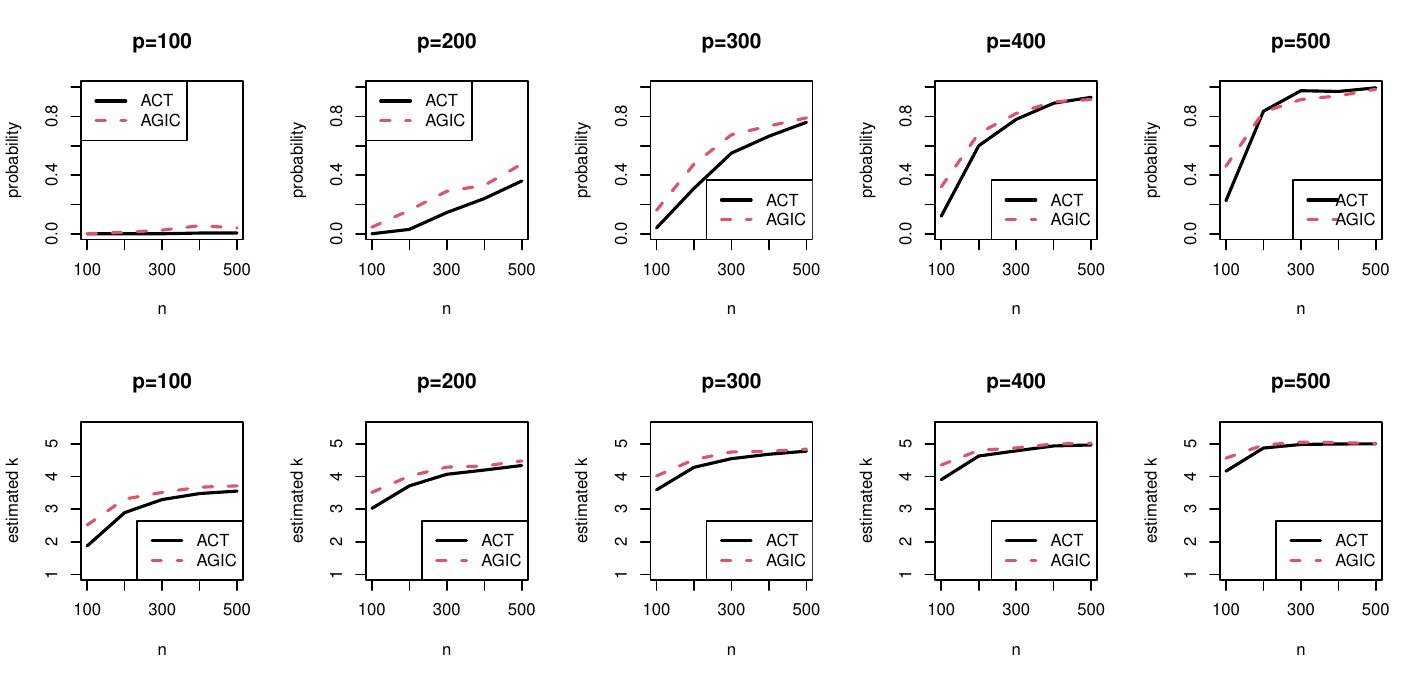}
\caption{Performances of ACT and AGIC in Simulation 6 under varying $p,n$ with $k=5$.}
\label{fig:2}
\end{figure}

For fewer spikes, i.e., $k=5$, we compare AGIC with ACT. It is seen from Figure \ref{fig:2} that, AGIC performs better than ACT, especially when $p$ or $n$ is relatively small (i.e., $p\leq 300$ or $n\leq 300$).

\begin{thebibliography}{52}
\newcommand{\enquote}[1]{``#1''}
\expandafter\ifx\csname natexlab\endcsname\relax\def\natexlab#1{#1}\fi

\bibitem[{Abdi and Williams(2010)}]{abdi2010principal}
Abdi, H. and Williams, L.~J. (2010), \enquote{Principal component analysis,}
  \textit{Wiley Interdisciplinary Reviews: Computational Statistics}, 2,
  433--459.

\bibitem[{Akaike(1974)}]{Akaike:1974}
Akaike, H. (1974), \enquote{A new look at the statistical identification
  model,} \textit{IEEE Transactions on Automatic Control}, 19, 716--723.

\bibitem[{Anderson(2003)}]{Anderson:2003}
Anderson, T.~W. (2003), \textit{An introduction to multivariate statistical
  analysis}, Wiley, New York.

\bibitem[{Bai and Ng(2002)}]{Bai:Ng:2002}
Bai, J. and Ng, S. (2002), \enquote{Determining the number of factors in
  approximate factor models,} \textit{Econometrica}, 70, 191--221.

\bibitem[{Bai et~al.(2018)Bai, Choi, and Fujikoshi}]{Bai:Choi:Fujikoshi:2018}
Bai, Z.~D., Choi, K., and Fujikoshi, Y. (2018), \enquote{Consistency of AIC and
  BIC in estimating the number of significant components in high-dimensional
  principal component analysis,} \textit{The Annals of Statistics}, 46,
  1050--1076.

\bibitem[{Bai and Ding(2012)}]{Bai:Ding:2012}
Bai, Z.~D. and Ding, X. (2012), \enquote{Estimation of spiked eigenvalues in
  spiked models,} \textit{Random Matrices: Theory and Applications}, 1,
  1150011.

\bibitem[{Bai and Silverstein(1998)}]{Bai:Silverstein:1998}
Bai, Z.~D. and Silverstein, J.~W. (1998), \enquote{No eigenvalues outside the
  support of the limiting spectral distribution of large dimensional random
  matrices,} \textit{The Annals of Probability}, 26, 316--345.

\bibitem[{Bai and Silverstein(2004)}]{Bai:Silverstein:2004}
--- (2004), \enquote{CLT for linear spectral statistics of large-dimensional
  sample covariance matrices,} \textit{The Annals of Probability}, 32,
  553--605.

\bibitem[{Bai and Yao(2008)}]{Bai:Yao:2008}
Bai, Z.~D. and Yao, J.~F. (2008), \enquote{Central limit theorems for
  eigenvalues in a spiked population model,} \textit{Annales de l'Institut
  Henri Poincare (B) Probability and Statistics}, 44, 447--474.

\bibitem[{Bai and Yao(2012)}]{Bai:Yao:2012}
--- (2012), \enquote{On sample eigenvalues in a generalized spiked population
  model,} \textit{Journal of Multivariate Analysis}, 106, 167--177.

\bibitem[{Baik and Silverstein(2006)}]{Baik:Silverstein:2006}
Baik, J. and Silverstein, J.~W. (2006), \enquote{Eigenvalues of large sample
  covariance matrices of spiked population models,} \textit{Journal of
  Multivariate Analysis}, 97, 1382--1408.

\bibitem[{Bao et~al.(2022)Bao, Ding, Wang, and Wang}]{bao2022statistical}
Bao, Z., Ding, X., Wang, J., and Wang, K. (2022), \enquote{Statistical
  inference for principal components of spiked covariance matrix,} \textit{The
  Annals of Statistics}, 50, 1144--1169.

\bibitem[{Bao et~al.(2015)Bao, Pan, and Zhou}]{Bao:Pan:Zhou:2015}
Bao, Z.~G., Pan, G.~M., and Zhou, W. (2015), \enquote{Universality for the
  largest eigenvalues of sample covariance matrices with general population,}
  \textit{The Annals of Statistics}, 43, 382--421.

\bibitem[{Bickel and Levina(2008)}]{bickel2008regularized}
Bickel, P.~J. and Levina, E. (2008), \enquote{Regularized estimation of large
  covariance matrices,} \textit{The Annals of Statistics}, 36, 199--227.

\bibitem[{Buja and Eyuboglu(1992)}]{buja1992remarks}
Buja, A. and Eyuboglu, N. (1992), \enquote{Remarks on parallel analysis,}
  \textit{Multivariate behavioral research}, 27, 509--540.

\bibitem[{Cai et~al.(2020)Cai, Han, and Pan}]{Cai:Han:Pan:2020}
Cai, T.~T., Han, X., and Pan, G. (2020), \enquote{Limiting laws for divergent
  spiked eigenvalues and largest non-spiked eigenvalues of sample covariance
  matrices,} \textit{The Annals of Statistics,}, 1255--1280.

\bibitem[{Cai et~al.(2010)Cai, Zhang, and Zhou}]{Cai:Zhang:Zhou:2010}
Cai, T.~T., Zhang, C.~H., and Zhou, H.~H. (2010), \enquote{Optimal rates of
  convergence for covariance matrix estimation,} \textit{The Annals of
  Statistics}, 38, 2118--2144.

\bibitem[{Chen and Li(2022)}]{chen2022determining}
Chen, Y. and Li, X. (2022), \enquote{Determining the Number of Factors in
  High-dimensional Generalized Latent Factor Models,} \textit{Biometrika}, 109,
  769--782.

\bibitem[{Dobriban(2020)}]{Dobriban:2020}
Dobriban, E. (2020), \enquote{Permutation methods for factor analysis and PCA,}
  \textit{The Annals of Statistics}, 48, 2824--2847.

\bibitem[{Dobriban and Owen(2019)}]{Dobriban:Owen:2019}
Dobriban, E. and Owen, A. (2019), \enquote{Deterministic parallel analysis: an
  improved method for selecting factors and principal components,}
  \textit{Journal of the Royal Statistical Society: Series B}, 81, 163--183.

\bibitem[{Fan et~al.(2022)Fan, Guo, and Zheng}]{Fan:Guo:Zheng:2022}
Fan, J., Guo, J., and Zheng, S. (2022), \enquote{Estimating number of factors
  by adjusted eigenvalues thresholding,} \textit{Journal of the American
  Statistical Association}, 117, 852--861.

\bibitem[{Fan et~al.(2018)Fan, Liu, and Wang}]{fan2018large}
Fan, J., Liu, H., and Wang, W. (2018), \enquote{Large covariance estimation
  through elliptical factor models,} \textit{Annals of statistics}, 46, 1383.

\bibitem[{Horn and Johnson(2013)}]{Horn:Johnson:2013}
Horn, R.~A. and Johnson, C.~R. (2013), \textit{Matrix Analysis}, Cambridge
  University Press, New York.

\bibitem[{Hu et~al.(2023)Hu, Keita, and Fu}]{Hu:Keita:Fu:2023}
Hu, J., Keita, S., and Fu, K. (2023), \enquote{Extreme eigenvalues of principal
  minors of random matrices with moment conditions,} \textit{Journal of the
  Korean Statistical Society}, 52, 715--735.

\bibitem[{Jiang and Bai(2021)}]{Jiang:Bai:2021}
Jiang, D. and Bai, Z. (2021), \enquote{Generalized four moment theorem and an
  application to CLT for spiked eigenvalues of high-dimensional covariance
  matrices,} \textit{Bernoulli}, 27, 274--294.

\bibitem[{Johansson(1998)}]{Johansson:1998}
Johansson, K. (1998), \enquote{On fluctuations of random Hermitian matrices,}
  \textit{Duke Mathematical Journal}, 91, 151--203.

\bibitem[{Johnstone(2001)}]{Johnstone:2001}
Johnstone, I.~M. (2001), \enquote{On the distribution of the largest eigenvalue
  in principal component analysis,} \textit{The Annals of Statistics}, 29,
  295--327.

\bibitem[{Jonsson(1982)}]{Jonsson:1982}
Jonsson, D. (1982), \enquote{Some limit theorems for the eigenvalues of a
  sample covariance matrix,} \textit{Journal of Multivariate Analysis}, 12,
  1--28.

\bibitem[{Ke et~al.(2023)Ke, Ma, and Lin}]{Ke:Ma:Lin:2023}
Ke, Z., Ma, Y., and Lin, X. (2023), \enquote{Estimation of the number of spiked
  eigenvalues in a covariance matrix by bulk eigenvalue matching analysis,}
  \textit{Journal of the American Statistical Association}, 118, 374--392.

\bibitem[{Ke(2016)}]{ke2016detecting}
Ke, Z.~T. (2016), \enquote{Detecting rare and weak spikes in large covariance
  matrices,} \textit{arXiv:1609.00883}.

\bibitem[{Kritchman and Nadler(2008)}]{Kritchman:Nadler:2008}
Kritchman, S. and Nadler, B. (2008), \enquote{Determining the number of
  components in a factor model from limited noise data,} \textit{Chemometrics
  and Intelligent Laboratory Systems}, 94, 19--32.

\bibitem[{Kritchman and Nadler(2009)}]{Kritchman:Nadler:2009}
--- (2009), \enquote{Non-parametric detection of the number of signals,
  hypothesis tests and random matrix theory,} \textit{IEEE Transactions on
  Signal Processing}, 57, 3930--3941.

\bibitem[{Li et~al.(2020)Li, Han, and Yao}]{Li:Han:Yao:2020}
Li, Z., Han, F., and Yao, J. (2020), \enquote{Asymptotic joint distribution of
  extreme eigenvalues and trace of large sample covariance matrix in a
  generalized spiked population model,} \textit{The Annals of Statistics},
  3138--3160.

\bibitem[{Ma(2012)}]{Ma:2012}
Ma, Z. (2012), \enquote{Accuracy of the Tracy-Widom limit for the extreme
  eigenvalues in the white Wishart matrices,} \textit{Bernoulli}, 18, 322--359.

\bibitem[{Nadakuditi and Edelman(2012)}]{Nadakuditi:Edelman:2008}
Nadakuditi, R.~R. and Edelman, A. (2012), \enquote{Sample Eigenvalues based
  detection of high-dimensional signals in white noise using relatively few
  samples,} \textit{IEEE Transactions on Signal Processing}, 56, 2625--2638.

\bibitem[{Nadler(2010)}]{Nadler:2010}
Nadler, B. (2010), \enquote{Nonparametric detection of signals by information
  theoretic criteria: Performance analysis and an improved estimator,}
  \textit{IEEE Transactions on Signal Processing}, 58, 2746--2756.

\bibitem[{Onatski(2009)}]{Onatski:2009}
Onatski, A. (2009), \enquote{Testing hypotheses about the number of factors in
  large factor models,} \textit{Econometrica}, 77, 1447--1479.

\bibitem[{Onatski et~al.(2013)Onatski, Moreira, and
  Hallin}]{Onatski:Moreira:Hallin:2013}
Onatski, A., Moreira, M.~J., and Hallin, M. (2013), \enquote{Asymptotic power
  of sphericity test for high-dimensional data,} \textit{The Annals of
  Statistics}, 43, 1204--1231.

\bibitem[{Passemier et~al.(2015)Passemier, Matthew, and
  Chen}]{Passemier:Mckay:Chen:2015}
Passemier, D., Matthew, and Chen, Y. (2015), \enquote{Asymptotic linear
  spectral statistics for spiked Hermitian random matrices,} \textit{Journal of
  Statistical Physics}, 160, 120--150.

\bibitem[{Passemier and Yao(2014)}]{Passemier:Yao:2014}
Passemier, D. and Yao, J. (2014), \enquote{Estimation of the number of spikes,
  possibly equal, in the high-dimensional case,} \textit{Journal of
  Multivariate Analysis}, 127, 173--183.

\bibitem[{Paul(2007)}]{Paul:2007}
Paul, D. (2007), \enquote{Asymptotics of sample eigenstructure for a large
  dimensional spiked covariance model,} \textit{Statistica Sinica}, 17,
  1617--1642.

\bibitem[{Saccenti and Timmerman(2017)}]{saccenti2017considering}
Saccenti, E. and Timmerman, M.~E. (2017), \enquote{Considering Horn' parallel
  analysis from a random matrix theory point of view,} \textit{Psychometrika},
  82, 186--209.

\bibitem[{Schwarz(1978)}]{Schwarz:1978}
Schwarz, G. (1978), \enquote{Estimating the dimension of a model,} \textit{The
  Annals of Statistics}, 6, 461--464.

\bibitem[{Silverstein(1995)}]{Silverstein:1995}
Silverstein, J.~W. (1995), \enquote{Strong convergence of the empirical
  distribution of eigenvalues of a large dimension random matrices,}
  \textit{Journal of Multivariate Analysis}, 54, 331--339.

\bibitem[{Wang et~al.(2014)Wang, Silverstein, and
  Yao}]{Wang:Silverstein:Yao:2014}
Wang, Q., Silverstein, J., and Yao, J. (2014), \enquote{A note on the CLT of
  the LSS for sample covariance matrix from a spiked population model,}
  \textit{Journal of Multivariate Analysis}, 130, 194--207.

\bibitem[{Wang and Yao(2013)}]{Wang:Yao:2013}
Wang, Q. and Yao, J. (2013), \enquote{On the sphericity test with
  large-dimensional observations,} \textit{Electronic Journal of Statistics},
  7, 2164--2192.

\bibitem[{Wang and Fan(2017)}]{Wang:Fan:2017}
Wang, W. and Fan, J. (2017), \enquote{Asymptotics of empirical eigenstructure
  for high dimensional spiked covariance,} \textit{The Annals of Statistics},
  45, 1342--1374.

\bibitem[{Wax and Kailath(1985)}]{Wax:Kailath:1985}
Wax, M. and Kailath, T. (1985), \enquote{Detection of signals by information
  theoretic criteria,} \textit{IEEE Transactions on Acoustics, Speech, and
  Signal Processing}, 33, 387--392.

\bibitem[{Weyl(1912)}]{Weyl:1912}
Weyl, H. (1912), \enquote{Der asymptotische Verteilungs gesetz der Eigenwerte
  linearer partieller Differentialgleichungen,} \textit{Mathematische Annalen},
  71, 441--479.

\bibitem[{Zhang et~al.(2022)Zhang, Zheng, Pan, and
  Zhong}]{Zhang:Zheng:Pan:Zhong:2022}
Zhang, Z., Zheng, Z., Pan, G., and Zhong, P. (2022), \enquote{Asymptotic
  independence of spiked eigenvalues and linear spectral statistics for large
  sample covariance matrices,} \textit{The Annals of Statistics}, 50,
  2205--2230.

\bibitem[{Zhao et~al.(1986)Zhao, Krishnaiah, and
  Bai}]{Zhao:Krishnaiah:Bai:1986}
Zhao, L.~C., Krishnaiah, P.~R., and Bai, Z.~D. (1986), \enquote{On detection of
  the number of signals in presence of white noise,} \textit{Journal of
  Multivariate Analysis}, 20, 1--25.

\bibitem[{Zhong et~al.(2022)Zhong, Su, and Fan}]{zhong2022empirical}
Zhong, X., Su, C., and Fan, Z. (2022), \enquote{Empirical Bayes PCA in high
  dimensions,} \textit{Journal of the Royal Statistical Society Series B:
  Statistical Methodology}, 84, 853--878.

\end{thebibliography}
\end{document}